\documentclass{article}
\usepackage[hmargin={25mm,25mm},vmargin={30mm,35mm}]{geometry}

\usepackage{amsmath,amsfonts,amssymb}
\usepackage[inline]{enumitem}
\usepackage{hyperref}
\usepackage{mathtools}
\usepackage{subfig}
\usepackage{bm}
\usepackage{float}
\usepackage{adjustbox}
\usepackage{makecell}
\usepackage{multirow}

\usepackage{todonotes} 
\newcommand*{\SavedEqref}{}
\let\SavedEqref\eqref
\renewcommand*{\eqref}[1]{%
	\begingroup
	\hypersetup{
		linkcolor=blue,
		linkbordercolor=blue,
	}%
	\SavedEqref{#1}%
	\endgroup
}

\usepackage{stmaryrd}
\usepackage{calc}
\usepackage{accents}

\newcommand{\jump}[1]{\left\llbracket{#1} \right\rrbracket}
\newcommand{\mean}[1]{ \left\{ \hspace{-0.10cm} \left\{ {#1} \right\} \hspace{-0.10cm} \right\}}
\newcommand{\meanBig}[1]{ \left\{ \hspace{-0.16cm} \left\{ {#1} \right\} \hspace{-0.16cm} \right\}_\face}
\newcommand{\tr}[1]{ \text{tr}\left({#1}\right)}
\newcommand{\tnsr}[1]{\bm{#1}}
\newcommand{\Tnsr}[1]{\mathbb{#1}}
\newcommand{\ub}[0]{\bm{u}}
\newcommand{\vb}[0]{\bm{v}}

\newcommand{\F}[0]{\tnsr{F}}
\newcommand{\V}[0]{\tnsr{V}}
\newcommand{\Fd}[0]{\boldsymbol{\mathfrak{F}}}
\newcommand{\eW}[0]{w}
\newcommand{\W}[0]{\boldsymbol{W}}
\newcommand{\Piola}[0]{\tnsr{P}}

\newcommand{\N}[0]{\bm{N}}
\newcommand{\I}[0]{\bm{1}}

\newcommand{\Grad}[0]{\nabla_{\hspace{-0.08cm}\X}}

\newcommand{\elem}[0]{T}
\newcommand{\face}[0]{F}
\newcommand{\edge}[0]{E}

\newcommand{\sumOnFEID}[0]{\sum_{\face \in \mathcal{F}_\elem^{i,\dirw}}}

\newcommand{\sumOnFDlm}[0]{\sum_{\face \in \mathcal{F}_\elem^\dirlm}}

\newcommand{\sumOnFEN}[0]{\sum_{\face \in \mathcal{F}_\elem^\neu}}
\newcommand{\sumOnElem}[0]{\sum_{\elem \in \mathcal{T}_h}}
\newcommand{\intE}[0]{\int_{\elem}}
\newcommand{\intf}[0]{\int_{\face}}
\newcommand{\X}[0]{\bm{X}}
\newcommand{\x}[0]{\bm{x}}
\newcommand{\E}[0]{\tnsr{E}}
\newcommand{\eg}[0]{\textit{e.g.}}
\newcommand{\ie}[0]{\textit{i.e.}}
\newcommand{\ea}[0]{\textit{et al.}}

\newcommand{\Th}{\mathcal{\elem}_h}
\newcommand{\Fh}{\mathcal{\face}_h}
\newcommand{\Eh}{\mathcal{\edge}_h}
\newcommand{\FT}{\mathcal{\face}_\elem}
\newcommand{\EF}{\mathcal{\edge}_\face}

\newcommand{\dir}{{\rm D}}
\newcommand{\dirw}{{\rm{D}_{\rm{N}}}}
\newcommand{\dirlm}{{\rm{D}_{\rm{L}}}}
\newcommand{\neu}{{\rm N}}
\newcommand{\internal}{{\rm i}}
\newcommand{\normal}{\boldsymbol{n}}
\newcommand{\Poly}[1]{\mathcal{P}^{#1}}
\newcommand{\Polyd}[2]{\mathbb{P}_{#1}^{#2}}
\newcommand{\uVec}{\boldsymbol{u}}
\newcommand{\vVec}{\boldsymbol{v}}
\newcommand{\wVec}{\boldsymbol{w}}
\newcommand{\zVec}{\boldsymbol{z}}
\newcommand{\yVec}{\boldsymbol{y}}
\newcommand{\sVec}{\boldsymbol{s}}

\newcommand{\Fig}[0]{Figure}
\newcommand{\Tab}[0]{Table}
\newcommand{\Sec}[0]{Section}
\newcommand{\Eq}[0]{Eq.}

\makeatletter
\newcommand{\dbtilde}[1]{{%
		\mathpalette\double@widetilde{#1}%
}}
\newcommand{\double@widetilde}[2]{%
	\sbox\z@{$\m@th#1\widetilde{#2}$}%
	\ht\z@=0.84 \ht\z@
	\widetilde{\box\z@}%
}
\makeatother

\begin{document}
\title{BR2 discontinuous Galerkin methods for finite hyperelastic deformations}
\author{BOTTI Lorenzo, VERZEROLI Luca}
\maketitle
%
	
	\begin{abstract}
In this work we introduce a dG framework for nonlinear elasticity based on a Bassi-Rebay (BR2) formulation. 
The framework encompasses compressible and incompressible hyperelastic materials and is capable of dealing with large deformations. 
In order to achieve stability, we combine higher-order lifting operators for the BR2 stabilization term with an adaptive stabilization strategy 
which relies on the BR2 Laplace operator stabilization and a penalty parameter based on the spectrum of the fourth-order elasticity tensor.
Dirichlet boundary conditions for the displacement can be imposed by means of Lagrange multipliers and Nitsche method. 
Efficiency of the solution strategy is achieved by means of state-of-the-art agglomeration based $h$-multigrid preconditioners and 
the code implementation supports distributed memory execution on modern parallel architectures. 
Several benchmark test cases are proposed in order to investigate some relevant computational aspects, 
namely the performance of the $h$-multigrid iterative solver varying the stabilization parameters and the influence of Dirichlet boundary conditions on Newton's method globalisation strategy.
	\end{abstract}
	
		


\section{Introduction}
Discontinuous Galerkin (dG) methods are widely employed in the field of Computational Fluid Dynamics (CFD) 
where they are appreciated for their turbulence modelling capabilities. 
In the last few decades, the interest in dG formulations for Computational Solid Mechanics (CSM) 
has been growing due to the following attractive features:
robustness with respect to mesh distortion, 
ability to deal with arbitrarily unstructured polytopal elements meshes, 
possibility to locally increase the accuracy by raising the polynomial degree in those regions where the solution is expected to be smooth,
availability of locking-free formulations in the incompressible and nearly-incompressible limits.
Despite those appealing properties the success of dG methods among CSM practitioners has been rather scarce, possibly because of the increased memory footprint and the lack of efficient solution strategies. An obvious downturn is related to the lack of dG modules in CSM commercial codes.

Several dG discretizations of \emph{linear elasticity} problems have been proposed and analysed in literature.
The $hp$-error analysis was first considered by Riviere \ea~\cite{Riviere2000} and 
Hansbo \ea~\cite{Hansbo2002} analysed the nearly and fully incompressible limits introducing a locking-free mixed formulation.
A Bassi-Rebay (BR2) dG method was proposed by Lew \ea~\cite{Lew2004}
and Cockburn \ea~\cite{Cockburn2006} introduced a Local Discontinuous Galerkin (LDG) method. 
Other locking-free implementations have been proposed by 
Wihler~\cite{Wihler2004,Wihler2006} and $hp$-adaptivity was considered by Houston \ea~\cite{Houston2006}.
Beam and plate modelling was tackled by Celiker~\cite{Celiker2004,Celiker2006} 
while Kirchhoff-Love linear shells were investigated by Guzey \ea~\cite{Guzey2006} and Noels~\cite{Noels2008}.
Kaufmann \ea~\cite{Kaufmann2009} exploited dG flexibility to simulate deformable bodies based on arbitrarily shaped polyhedral elements meshes.
Plasticity problems in the small deformation regime were studied by Djoko \ea~\cite{Djoko2007,Djoko2007a}.

Concerning the use of dG formulations in the context of \textit{nonlinear elasticity} problems, the following research efforts deserve to be mentioned.
In 2006 Noels and Radovitzky~\cite{Noels2006} tackled large strains of hyperelastic bodies and 
Eyck and Lew~\cite{TenEyck2006} proposed a dG formulation based on the Bassi-Rebay (BR1)~\cite{Bassi1997} gradient reconstruction. 
A novel adaptive stabilization approach for the latter formulation was proposed in Eyck \ea~\cite{Eyck2008,Eyck2008a}.
Whiteley~\cite{Whiteley2009} investigated locking phenomena in nonlinear elasticity showing the advantage of dG methods in the incompressible limit. 
Baroli \ea~\cite{Baroli2013} devised a total Lagrangian Interior Penalty (IP) dG formulation for incompressible and anisotropic soft living materials.
Challenging application oriented contributions are collected in what follows.
Becker and Noel~\cite{Becker2013} modelled cracks initiation and propagation by means of Kirchhoff-Love shell elements.
McBride and Reddy~\cite{McBride2009a} introduced a logarithmic hyper-elastoplastic model for the finite-deformation regime.
Liu \ea~\cite{Liu2013b} tackled hypo- and hyper-elastoplastic problems through an updated Lagrangian formulation.
Feistauer \ea~\cite{Feistauer2019} employed a spatial dG discretization for an elasto-dynamic system. 
To conclude, Kosis \ea~\cite{Kosik2015} considered a space-time dG formulations of the Fluid-Structure-Interaction (FSI) problem involving a compressible Newtonian fluid and a Saint Venant-Kirchhoff material. 

More recently, the introduction of Hybridizable Discontinuous Galerkin methods (HDG) has further increased the popularity of 
discontinuous Finite Element methods among CSM practitioners. Nguyen and
Peraire~\cite{Nguyen2012} proposed an HDG framework for continuum mechanics.
Kabaria \ea~\cite{Kabaria2015} proposed an HDG method for nonlinear elasticity and a suitable stabilization strategy was later proposed by Cockburn and Shen~\cite{Cockburn2019}.
Terrana \ea~\cite{Terrana2019} applied HDG methods to thin structures presenting buckling phenomena. 
Botti \ea~\cite{BottiM2017} analysed a Hybrid High-Order (HHO) methods for nonlinear elasticity with small deformations.
Abbas \ea~\cite{Abbas2018} presented a stabilized and an unstabilized HHO method for finite deformations of hyperelastic materials. 
HHO methods has been applied to incremental associative plasticity and elastoplastic deformations in Abbas \ea~\cite{Abbas2019,Abbas2019a} 
while Chouly \ea~\cite{Chouly2020} applied HHO methods to contact mechanics.\\
\bigskip


In this work we introduce an effective framework for finite deformations of elastic solids. The framework relies on the following ingredients:
\begin{enumerate*}
\item BR2 dG discretization of the Lagrangian equation of motion for hyperelastic materials with adaptive stabilization strategy;
\item dG discretization of the incompressibility constraint in Lagrangian formulation; 
\item implementation of the Lagrange multipliers method for the imposition of Dirichlet boundary conditions (BCs);
\item agglomeration based $h$-multigrid solution strategy for the fully coupled formulation.
\end{enumerate*}
Up to the author's knowledge all the aforementioned ingredients but the second are original contributions of the present manuscript.
Moreover, the numerical investigation performed on challenging 2D and 3D test cases will focus on the effectiveness of the dG framework in practice.
As a first point, we demonstrate that, thanks to the combination of adaptive stabilisation and multigrid solution strategy, 
the efficacy of the solver is maintained over a wide range of stabilization parameters values. 
Accordingly, the burden of choosing the stabilization coefficients is streamlined. 
Moreover, we show that imposing Dirichlet boundary conditions by means of Lagrange multipliers 
significantly improves the robustness of Newton's method globalisation strategy and reduces the computational cost. 

The material is organised as follows.
In \Sec~\ref{sec:nonlinear_el_prob}, we introduce the Lagrangian formulation of the nonlinear elasticity problem 
and four hyperelastic constitutive laws: two compressible and an incompressible neo-Hookean models and the Saint-Venant Kirchhoff model. 
In \Sec~\ref{sec:dg_formulation}, we present the BR2 formulation of the Lagrangian equation of motion, possibly coupled with the incompressibility constraint.
In order to deal with large deformation, an adaptive stabilization strategy and an incremental load method are required, see \Sec~\ref{sec:incremental_method}. 
In \Sec~\ref{sec:numerical_results}, after establishing the convergence rates,
challenging benchmark test cases are presented considering elastic bodies subjected to compression, torsion and stretch.

\section{Nonlinear elasticity problem} \label{sec:nonlinear_el_prob}
We consider the classical problem of seeking the static equilibrium of an elastic body 
undergoing finite deformations, see for example Odgen~\cite{Ogde}, Ciarlet and Philippe~\cite{Ciarlet1999a}, 
Gurtin \ea~\cite{Gurtin2009}, Tadmor et al.~\cite{Tadmor2011} and Bonnet \ea~\cite{Bonet2016} for additional details.

The elastic continuum body in the reference configuration occupies the bounded connected domain 
$\Omega \in \mathbb{R}^d$, $d \in {2,3}$, with Lipschitz continuous boundary $\partial \Omega$.
The material points $\X \in \Omega$ are mapped into 
spatial points $\x = \X + \ub(\X)$, where $\ub: \Omega \rightarrow \mathbb{R}^d$ is the displacement mapping.
The body in deformed configuration occupies $\Omega^{+} \coloneqq \{\X + \ub(\X), \X \in \Omega\}$.
The deformation gradient reads
\begin{equation}
\label{eq:defGrad}
\F \coloneqq \Grad \ub + \I, 
\end{equation}
where $\I$ is the second-order identity tensor 
and $\Grad$ is the gradient operator in the reference configuration.  
Introducing the density of the material of the body in reference and deformed configuration
$\rho:\Omega \rightarrow \mathbb{R}^+$ and $\rho^+:\Omega^+ \rightarrow \mathbb{R}^+$, respectively, 
mass conservation implies that, see \eg~\cite{Ogde}
\begin{equation}
\label{eq:jacMap}
J:=\det(\F(\ub)) = \frac{\rho(\X)}{\rho^+(\X + \ub(\X))}>0.
\end{equation}

In the Lagrangian framework, the elasticity problem consists of finding the displacement mapping $\ub$ such that 
\begin{subequations}
	\label{elastProb}
	\begin{alignat}{2}
	-\Grad \cdot \Piola &= \rho \, \bm{f}^+ &\qquad& \text{in } \Omega, \label{eq:linear_momentum}\\
	\ub &= \mathbf{g}_\dir &\qquad& \text{on } \partial \Omega_\dir, \label{eq:lmDbc}\\
	\Piola \N &= \mathbf{g}_\neu &\qquad& \text{on } \partial \Omega_\neu, \label{eq:lmNbc}
	\end{alignat}
\end{subequations}
where $\N$ is the unit normal vector pointing out of $\partial \Omega$,
$\bm{f}^+: \Omega^+ \rightarrow \mathbb{R}^d$ is the known body force per unit mass, 
$\bm{g}_\neu$ is the traction force per unit area imposed on the Neumann boundary $\partial \Omega_\neu$
and $\bm{g}_\dir$ is the displacement vector imposed on the Dirichlet boundary $\partial \Omega_\dir$. 
$\Piola(\F(\ub)): \mathbb{R}^{d \times d} \rightarrow \mathbb{R}^d$ is the first Piola-Kirchhoff stress tensor:
a stress measure that describes the response of the body to 
the external solicitations $\bm{f}^{+},\bm{g}_\neu,\text{ and }\bm{g}_\dir$. 
We consider \textit{hyperelastic} materials, namely materials whose mechanical properties are characterized 
by the \textit{strain-energy function} $\widehat{\eW}(\F(\ub)):\mathbb{R}^{d \times d} \rightarrow \mathbb{R}$, such that 
$\Piola \coloneqq \frac{\partial{\widehat{\eW}}}{\partial \F}$.
It is assumed that $\partial \Omega = \overline{\partial \Omega}_\dir \bigcup \overline{\partial \Omega}_\neu$, 
$\partial \Omega_\dir \bigcap \partial \Omega_\neu = \emptyset$ and both 
$\partial \Omega_\dir$ and $\partial \Omega_\neu$ have non-zero $(d-1)$-dimensional Hausdorff measure.

The description of isochoric deformations is obtained satisfying, 
simultaneously, (\ref{eq:linear_momentum}) and the incompressibility constraint $J = 1$,
which directly follows from \eqref{eq:jacMap}.
For an \textit{incompressible} hyperelastic material the strain energy function 
$\eW:\mathbb{R}^{d \times d} \times \mathbb{R} \rightarrow \mathbb{R}$, reads
$$
\eW(\F(\ub),q)\coloneqq \widehat{\eW}(\F(\ub)) + (J-1)q
$$
where $q$ is an arbitrary Lagrange multiplier, and 
$\Piola(\F(\ub),q) = \frac{\partial{\eW}}{\partial \F}$.

Let $\V$ be the set of all kinematically admissible displacements
which satisfy the Dirichlet condition \eqref{eq:lmDbc}
and $Q$ be the set of admissible Lagrange multipliers, 
we define the energy functional $\mathcal{W}: \V \times Q \rightarrow \mathbb{R}$ such that
$$
\mathcal{W}(\vb,q) = \int_\Omega w(\F(\vb),q) \, d\Omega - \int_\Omega \rho \bm{f}^{+} \cdot \vb \, d\Omega - \int_{\partial \Omega_\neu} \bm{g}_\neu \cdot \vb \, d\partial\Omega_\neu.
$$
The static equilibrium of problem \eqref{elastProb} constrained by $J=1$ consists in finding $(\ub,p) \in \V \times Q$
which satisfy the following weak form of the Euler-Lagrange equations
\begin{equation}
\begin{aligned}
\label{elastProbWeak}
0 = \frac{d}{d \epsilon} \mathcal{W}(\ub + \epsilon \delta \vVec,p) \bigg|_{\epsilon=0} &= \int_\Omega \Piola(\F(\ub),p) : \Grad (\delta\vb) - \int_\Omega \rho \bm{f}^{+} \cdot \delta\vb - \int_{\partial \Omega_\neu} \bm{g}_\neu \cdot \delta\vb,  \\
0 = \frac{d}{d \epsilon} \mathcal{W}(\ub ,p+ \epsilon \delta q) \bigg|_{\epsilon=0} &= \int_\Omega (J-1)\, \delta q, 
\end{aligned}
\end{equation}
for all virtual displacements $\delta\vb$ that satisfy a homogeneous Dirichlet condition on $\partial \Omega_\dir$
and for all $\delta q$.
According to \eqref{elastProbWeak}, $(\ub,p)$ is a stationary point of $\mathcal{W}(\vb,q)$.

\subsection{First Piola-Kirchhoff stress tensor} \label{sec:firstPK_tensor}
We restrict our investigation on the following strain-energy functions 
\begin{align}
\text {Saint Venant-Kirchhoff (SVK), \cite{Ciarlet1999a}}:    \widehat{\eW}(\F) &= \mu \, \E : \E + \frac{1}{2} \lambda (\tr{\E})^2,   \label{eq:svk_energy}\\
\text {Compressible neo-Hookean (NHK-C), \cite{Pence2015}}:    \widehat{\eW}(\F) &= \frac{\mu}{2} (\tr{\tnsr{C}} - d) -\mu \ln J + \frac{\lambda}{2} \Theta^2(J), \label{eq:neoh_energy} \\
\text {Incompressible neo-Hookean (NHK-I), \cite{Tadmor2011}}:  \widehat{\eW}(\F) &= \frac{\mu}{2} (\tr{\tnsr{C}} - d).\label{eq:neoh_inc_energy} \\
\text {Cavitating neo-Hookean (NHK-CAV), \cite{Kabaria2015}}:  \widehat{\eW}(\F)&= \frac{2 \mu}{3^{5/4}} (\tr{\tnsr{C}})^{\frac{3}{4}} -\mu \ln J + \frac{\lambda}{2} (\ln{J})^2,
\label{eq:neoh_energy_cavitation}
\end{align}
where  $\Theta(J) = \ln(J)$ (for other choices of $\Theta$ see \cite{Brink1996}),
$\tnsr{C} = \F^\intercal \F$ is the right Cauchy-Green tensor 
and $\E = \frac{1}{2} (\tnsr{C} - \I)$ is the Green-Lagrange strain tensor.

In \eqref{eq:svk_energy}-\eqref{eq:neoh_energy_cavitation}, $\mu$ and $\lambda$ are the Lam\'e parameters 
that can be written in terms of the Poisson's coefficient $\nu$ and the Young's modulus $E$ through the following relations
\begin{equation*}
\mu = \dfrac{E}{2 (1 + \nu) }, \qquad \lambda = \frac{\nu E}{ (1 + \nu) (1 - 2 \nu)}.
\end{equation*}

Introducing the second Piola-Kirchhoff stress tensor 
$
\tnsr{S} = \frac{\partial \widehat{\eW}}{\partial \tnsr{E}} = 2 \frac{\partial \widehat{\eW}}{\partial \tnsr{C}},
$
such that
\begin{align}
\text {Saint Venant-Kirchhoff (SVK)}:      \tnsr{S} &= 2 \mu \E + \lambda (\tr{\E}) \I \label{eq:svk}; \\
\text {Compressible neo-Hookean (NHK-C)}:   \tnsr{S} &= \mu \left( \I - \tnsr{C}^{-1}\right)  + \lambda J \Theta(J) \Theta'(J) \tnsr{C}^{-1} \label{eq:neohookean_compr}; \\
\text {Incompressible neo-Hookean (NHK-I)}: \tnsr{S} &= \mu \I \label{eq:neohookean_inc}; \\
\text {Cavitating neo-Hookean (NHK-CAV)}: \tnsr{S} &= \frac{\mu}{3^{1/4}} \tr{\tnsr{C}}^{-\frac{1}{4}} \I - \mu \tnsr{C}^{-1} + \lambda \ln(J) \tnsr{C}^{-1} \label{eq:neohookean_cavitation}.
\end{align}
the first Piola-Kirchhoff stress tensor can be conveniently rewritten as follows
\begin{equation}
\Piola = 
\Piola(\F(\ub), p)  =  \F \tnsr{S} - p \, J \, \F^{-\intercal},
\label{eq:piola}
\end{equation}
where $p = -q$ is the \emph{hydrostatic pressure}.
For future use, we denote by SVK-C and SVK-I the Saint Venant-Kirchhoff (SVK) law used in the compressible and in the fully incompressible ($J=1$) case, respectively. 

\subsection{Fourth-order elasticity tensor} \label{sec:elastic_fourth-order_tensor}
The fourth-order elasticity tensor associated to the elastic strain-energy energy function  
is computed as follows 
\begin{equation*}
\Tnsr{A} = \Tnsr{A}(\F(\uVec), p)=\frac{\partial }{\partial \F}\frac{\partial {\eW}}{\partial \F}=\frac{\partial \Piola}{\partial \F}.
\label{eq:elastic_modulus_piola}
\end{equation*}
According to \eqref{eq:piola}, the tensor coefficients $\Tnsr{A}_{ijkl}$, such that 
$\Tnsr{A} = \Tnsr{A}_{ijkl}  \; \bm{e}_i \otimes \bm{e}_j \otimes \bm{e}_k \otimes \bm{e}_l$,
reads
\begin{equation*}
\Tnsr{A}_{ijkl}=  \frac{\partial  }{\partial F_{kl}} (\F \tnsr{S} )_{ij} - p \frac{\partial}{\partial F_{kl}} \left( J \, F^{-\intercal}_{ij} \right)
= \delta_{ik} S_{lj} + F_{im} \frac{\partial S_{mj}}{\partial F_{kl}} - p \, J \left( F_{kl}^{-\intercal} F_{ij}^{-\intercal} -  F_{il}^{-\intercal} F_{jk}^{-1}\right).
\end{equation*}
Since $\tnsr{S}$ depends on the strain-energy function, usually expressed in terms of the strain tensors $\tnsr{C}$ and $\tnsr{E}$, the following relations are helpful to obtain the final expression of $\Tnsr{A}$:
\begin{equation*}
\frac{\partial S_{mj}}{\partial F_{kl}} = \frac{\partial S_{mj}}{\partial E_{qr}} \frac{\partial E_{qr}}{\partial F_{kl}}, \quad \frac{\partial S_{mj}}{\partial F_{kl}} = \frac{\partial S_{mj}}{\partial C_{qr}} \frac{\partial C_{qr}}{\partial F_{kl}} \quad \text{and} \quad
\frac{\partial E_{qr}}{\partial F_{kl}} = \frac{1}{2} \frac{\partial C_{qr}}{\partial F_{kl}} =\frac{1}{2} \left(  \delta_{ql} F_{kr} +  \delta_{rl} F_{kq} \right).
\end{equation*}
To conclude, according to \eqref{eq:svk}-\eqref{eq:neohookean_compr}-\eqref{eq:neohookean_inc}-\eqref{eq:neohookean_cavitation}, we get
\begin{align*}
\text {SVK}:      \frac{\partial S_{mj}}{\partial E_{qr}} &= 2 \mu \delta_{mq} \delta_{jr} + \lambda \delta_{mj} \delta_{qr}; \\
\text {NHK-C}:   \frac{\partial S_{mj}}{\partial C_{qr}} &= \mu C^{-1}_{mq} C^{-1}_{rj} + \lambda \bigg(-J \Theta \Theta' C^{-1}_{mq} C^{-1}_{rj} +\frac{J}{2}\left( \Theta \Theta' + J \Theta'^2 + J \Theta \Theta''\right) C^{-1}_{qr} C^{-1}_{mj} \bigg); \\
\text {NHK-I}: \frac{\partial S_{mj}}{\partial C_{qr}} &= 0;\\
\text {NHK-CAV}:   \frac{\partial S_{mj}}{\partial C_{qr}} & = -\frac{1}{4}\frac{\mu}{3^{1/4}}  \tr{\tnsr{C}}^{-\frac{5}{4}} \delta_{mj} \delta_{qr} +\mu C^{-1}_{mq} C^{-1}_{rj} + \lambda \bigg(-\ln( J) \, C^{-1}_{mq} C^{-1}_{rj} +\frac{1}{2} C^{-1}_{qr} C^{-1}_{mj} \bigg). 
\end{align*}

\section{The BR2 dG formulation of nonlinear elasticity} \label{sec:dg_formulation}
\subsection{Mesh setting}
\label{subsec:meshset}
We define a spatial meshes $\Th$ as a finite collections of disjoint mesh elements where
$h_\elem$ denotes the diameter of a \textit{mesh element} $\elem \in \Th$ and $h\coloneqq\max_{\elem\in\Th}h_\elem>0$ is the meshstep size.
$\Th$ is such that $\bigcup_{\elem \in \mathcal{T}_h} \overline{T} = \overline{\Omega}_h$, and either one of the following two conditions is satisfied
\begin{equation}
\label{eq:OmegaD}
\begin{array}{l}
  \Omega_h \equiv \Omega, \\
  \mbox{$\Omega_h$ is a suitable approximation of $\Omega$, meaning that $\lim_{h\rightarrow0} \Omega_h = \Omega$}.
\end{array}
\end{equation}
The mesh skeleton $\bigcup_{\elem \in \Th} \partial \elem$ is partitioned into a finite collection of mesh faces $\Fh$
such that, for each $\face \in \Fh$, one of the following two conditions is satisfied:
\begin{enumerate}[label=(\roman*)]
	\item There exist $\elem,\elem' \in \Th$, with $\elem \neq \elem'$, such that $\face = \partial \elem \cap \partial \elem'$, meaning that $\face$ is an \textit{internal face}.
	\item There is $\elem \in \Th$ such that $\face = \partial \elem \cap \partial \Omega_h$, meaning that $\face$ is a \textit{boundary face}.
\end{enumerate}
For each mesh element $\elem \in \Th$, the set $\FT = \{\face \in \Fh : \face \subset \partial \elem\}$ denotes the faces composing the element boundary $\partial \elem$.

We will consider two strategies for imposing Dirichlet boundary conditions, namely Nitsche method and Lagrange multipliers method. 
Let $\partial \Omega_{h,\dir} = \partial \Omega_{h,\dirw} \bigcup \partial \Omega_{h,\dirlm}$, 
where $\partial \Omega_{h,\dirw}$ and $\partial \Omega_{h,\dirlm}$ are the Nitsche and the Lagrange multipliers partitions of the Dirichlet boundary, respectively. 
We define four disjoint subsets of the set $\FT$: 
\begin{enumerate}
	\item $\FT^\dirw = \{\face \in \FT : \face \subset \partial \Omega_{h,\dirw}\}$: the set of Dirichlet faces where boundary conditions are weakly enforced using Nitsche method;
	\item $\FT^\dirlm = \{\face \in \FT : \face \subset \partial \Omega_{h,\dirlm}\}$: the set of Dirichlet faces where boundary conditions are enforced using Lagrange multipliers;
	\item $\FT^\neu = \{\face \in \FT : \face \subset \partial \Omega_{h,\neu}\}$: the set of Neumann boundary faces;
	\item $\FT^\internal \coloneqq \FT\setminus\big(\FT^\dir\cup\FT^\neu\big)$: the set of internal faces.
\end{enumerate}
For future use, we also let $\FT^{\internal,{\dirw}}\coloneqq\FT^\internal\cup\FT^{\dirw}$. 
For all $T\in\Th$ and all $F\in\FT$, $\normal_{\elem \face}$ denotes the normal vector to $\face$ pointing out of $\elem$. 
We remark that is case of boundary faces $\normal_{\elem \face}$ can be equal to $\N$ or an approximation of $\N$, see \eqref{eq:OmegaD}.

\subsubsection{Settings for Dirichlet BCs imposed by means of the Lagrange multipliers method}
\label{subsec:meshsetlm}
The Lagrange multipliers method requires further settings.
$\partial \Omega_{h,\dirlm}$ is partitioned into smooth patches $\Pi_{h,n}$, $n=1,...,N$, 
such that 
\begin{enumerate}
\item $\bigcup_{\Pi_{h,n} \in \partial \Omega_{h,\dirlm}} \overline{\Pi}_{h,n} = \overline{\partial \Omega}_{h,\dirlm}$
\item the normal vector $\N$ varies continuously over $\Pi_n = \lim_{h\rightarrow0} \Pi_{h,n}$,
\end{enumerate}
We define the sharp corners of $\partial \Omega_{h,\dirlm}$ as $\Gamma^i_{h} := \partial \Pi_{h,n} \cap \partial \Pi_{h,l}$, with $n,l=1,...,N, n \neq l$.
The boundary of $\partial \Omega_{h,\dirlm}$ is defined as follows 
$$ \Gamma^b_{h} := \left(\partial \Omega_{h,\dirlm} \cap \partial \Omega_{h,\dirw} \right) \bigcup \left( \partial \Omega_{h,\dirlm} \cap \partial \Omega_{h,\neu} \right).$$
To conclude we let $\Gamma_h = \Gamma_h^i \bigcup \Gamma^b_h$ be the set collecting all the sharp corners and the boundary of $\partial \Omega_{h,\dirlm}$.

Let $\Fh^\dirlm$ be the set collecting all Dirichlet boundary faces where the Lagrange multipliers method is employed. 
For each mesh face $\face \in \Fh^\dirlm$, the edges composing the face boundary are collected in the set $\EF$ 
such that $\bigcup_{\edge \in \EF} \overline{\edge} = \partial \face$.
We define two disjoint subsets of the set $\EF$: 
\begin{enumerate}
	\item $\EF^b = \{\edge \in \EF : \edge \subset \Gamma_h\}$: the set of \textit{boundary edges};
	\item $\EF^\internal \coloneqq \EF \setminus \EF^b$: the set of \textit{internal edges}.
\end{enumerate}

\subsubsection{Numerical integration over reference mesh entities}
In order to be able to numerically integrate over mesh elements, mesh faces and mesh edges, we require that, for any $Y$ element or face: 
\begin{enumerate}
	\item there exists a reference entity 
	$\widehat{Y}$ of standardized shape and a polynomial mapping $\Psi_Y : \widehat{Y} \rightarrow Y$ such that $Y = \Psi_Y(\widehat{Y})$;
	\item quadrature rules of arbitrary order are available on the reference entity $\widehat{Y}$. 
\end{enumerate}
From the geometrical viewpoint reference entities read as follows
\begin{equation*}
\widehat{Y}=
\begin{cases}
\text{point}, & \text{if $Y \in \Eh$,} \\
\text{line segment}, & \text{if $Y \in \Fh$,} \\
\text{polygon}, & \text{if $Y \in \Th$, }
\end{cases} \; \text{if $d=2$}; \qquad
\widehat{Y}=
\begin{cases}
\text{line segment}, & \text{if $Y \in \Eh$,} \\
\text{polygon}, & \text{if $Y \in \Fh$,} \\
\text{polyhedron}, & \text{if $Y \in \Th$, }
\end{cases} \; \text{if $d=3$}. \\
\end{equation*}

\subsection{dG formulation} \label{subsec:dg_formulation}
We denote by $\Polyd{d}{\ell}$ the space of $d$-variate polynomials of total degree $\le\ell$.
For each $\elem \in \Th$, we denote by $\Poly{\ell}(\elem)$ the space spanned by the restriction of $\Polyd{d}{\ell}$ to $T$
and by $\Poly{\ell}(\face)$ the space spanned by the restriction of $\Polyd{d-1}{\ell}$ to $\face$. 
Fix a polynomial degree $k\ge 1$ and let $\elem\in\Th$.
We define the \emph{local discrete gradient} $\boldsymbol{\mathfrak{G}}_\elem^k : H^1(\Th)^d \rightarrow \Poly{k}(\elem)^{d\times d}$ 
such that, for all $\vVec\in H^1(\Th)^d$,
\[
\int_\elem \boldsymbol{\mathfrak{G}}_\elem^k(\vVec) : \boldsymbol{\tau}
\coloneqq \int_\elem \nabla \vVec_{|\elem} : \boldsymbol{\tau}
- \sum_{\face \in \FT^{\internal,\dirw}} \frac12 \int_\face \left( \normal_{\elem \face} \otimes \jump{\vVec}_{\elem \face} \right) : \boldsymbol{\tau}
\qquad
\forall\boldsymbol{\tau} \in \Poly{k}(\elem)^{d \times d},
\]
where, for each $\face\in\FT^{\internal,\dirw}$, the jump of $\vVec$ across $\face$ is defined as
\begin{equation} 
\label{eq:jumpDG}
\jump{\vVec}_{\elem \face} \coloneqq \begin{cases}
\vVec_{|\elem} - \vVec_{|\elem'} & \text{if $\face \in \FT^\internal\cap\mathcal{\face}_{\elem'}^\internal$ with $\elem,\elem'\in\Th$, $\elem \neq \elem'$}, 
\\
2(\vVec_{|\elem} - \boldsymbol{g}_\dir) & \text{if $\face\in\FT^\dirw$}.
\end{cases}
\end{equation}
Introducing, for any $\face\in\FT^{\internal,\dirw}$ and any integer $\ell\geq0$, the jump lifting operator $\boldsymbol{\mathfrak{R}}_{\face \elem}^\ell: L^2(\face)^d \rightarrow \Poly{\ell}(\elem)^{d \times d}$ such that, for all $\boldsymbol{\varphi}\in L^2(\face)^d$ and all $\boldsymbol{\tau}\in\Poly{\ell}(\elem)^{d\times d}$,
\[
\int_\elem \boldsymbol{\mathfrak{R}}_{\face \elem}^\ell (\boldsymbol{\varphi}) : \boldsymbol{\tau}
=
\frac{1}{2} \int_\face \left(\normal_{\elem \face} \otimes \boldsymbol{\varphi} \right) : \boldsymbol{\tau},
\]
it holds, for all $\vVec\in H^1(\Th)^d$,
\begin{equation}
\label{eq:gradDG}
\boldsymbol{\mathfrak{G}}_T^k(\vVec) =
\nabla \vVec_{|T} - \sum_{\face \in\FT^{\internal,\dirw}} \boldsymbol{\mathfrak{R}}_{\face \elem}^k \left(\jump{\vVec}_{\elem \face} \right).
\end{equation}

Based on definition \eqref{eq:defGrad}, we introduce 
two discrete versions of the deformation gradient  
\begin{alignat*}{2}
\mbox{for any} \, \elem \in \Th &: \Fd_\elem^k(\vVec) &\;\coloneqq \;& \boldsymbol{\mathfrak{G}}_\elem^k(\vVec) + \I; \\
\mbox{for any} \, \elem \in \Th, \, \face \in \FT^{\internal,\dirw} &: \Fd_{\elem\face}^{k+1}(\vVec) &\;\coloneqq\;& \nabla \vVec_{|\elem} - \boldsymbol{\mathfrak{R}}_{\face \elem}^{k+1} \left(\jump{\vVec}_{\elem \face} \right) + \I.
\end{alignat*}
We remark that $\Fd_\elem^k$ relies on jump contributions over $\partial \elem$ and jump lifting operator of degree $k$,
while $\Fd_{\elem\face}^{k+1}$ relies solely on the jump lifting operator of degree $k{+}1$ over $\face$.  
The idea to employ $k{+}1$ lifting operators to get rid of stabilization parameters was first proposed by John \ea~\cite{John2016}
in the context of LDG discretizations of the Laplace operator with rigorous analysis covering the case of simplicial meshes.

We introduce the scalar- and vector-valued \emph{broken} polynomial spaces
\begin{equation}
\begin{aligned}
\Poly{k}(\Th) &\coloneqq \left\{q_h = \text{$(q_{\elem})_{\elem \in \Th}: q_T \in \Poly{k}(\elem)$ for all $\elem\in\Th$} \right\}, \\
\Poly{k}(\Th)^d &\coloneqq \left\{\vVec_h =  \text{$(\vVec_{\elem})_{\elem \in \Th} : \vVec_\elem \in \Poly{k}(\elem)^d$ for all $\elem\in\Th$} \right\}. 
\end{aligned}\nonumber
\end{equation}
and the vector-valued polynomial space for Lagrange multipliers over Dirichlet boundaries
\begin{equation}
\Poly{k}(\Fh^\dirlm)^d \coloneqq \left\{\widehat{\sVec}_h  = \text{$(\widehat{\sVec}_{\face})_{\face \in \Fh^\dirlm} : \widehat{\sVec}_{\face} \in\Poly{k}(\face)^d$ for all $\face \in \Fh^\dirlm$} \right\}. \nonumber
\end{equation}

Let the $k\geq1$ denote the polynomial degree, let a mesh element $\elem \in \Th$ and a mesh face $\face \in \Fh^{\dirlm}$ be fixed. 
Given $(\ub_h, p_h, \bm{\widehat{\lambda}}_h) \in \Poly{k}(\Th)^d \times \Poly{k}(\Th) \times \Poly{k}(\Fh^\dirlm)^d$, 
the local residuals 
\begin{itemize}
\item $r^{\text{Lem}}_\elem\left( (\ub_h,p_h,\bm{\widehat{\lambda}}_h); \bullet \right) : \Poly{k}(\elem)^d \to \mathbb{R}$ of the discrete Lagrangian equation of motion, 
\item $r^{\text{ic}}_\elem((\ub_h,p_h); \bullet): \Poly{k}(\elem) \to \mathbb{R}$ of the discrete incompressibility constraint, 
\item $r^{\dirlm}_\face\left((\ub_h,\widehat{\lambda}_h); \bullet \right): \Poly{k}(\face)^d \to \mathbb{R}$  of the constraints on Dirichlet boundary, 
\end{itemize}
are such that: for all ${\vVec}_\elem \in \Poly{k}(\elem)^d$, all $q_\elem \in \Poly{k}(\elem)$ and all $\widehat{\sVec}_{\face} \in \Poly{k}(\face)$
\begin{align}
r^{\text{Lem}}_\elem\left( (\ub_h, p_h, \bm{\widehat{\lambda}}_h); {\vVec}_\elem \right) &\coloneqq \intE {\Piola}(\Fd^{k}_\elem(\ub_h),p_h) : \nabla {\vVec}_\elem 
-\sumOnFEID \intf \left[ \meanBig{\Piola (\Fd_{\elem\face}^{k+1}(\ub_h),p_h)} \, \normal_{\elem\face} \right] \cdot {\vVec}_\elem + \nonumber \\
&+\sumOnFEID \intf  \eta_{\face} \, \left[\meanBig{\boldsymbol{\mathfrak{R}}_{\face \elem}^k(\jump{\uVec}_{\elem \face}) } \normal_{\elem\face} \right]\cdot {\vVec}_\elem 
 + \nonumber \\  
&{-}\intE \rho \bm{f}^{+} \cdot {\vVec}_\elem -\sumOnFEN \intf \bm{g}_\neu \cdot {\vVec}_\elem 
-\sumOnFDlm \intf \bm{\widehat{\lambda}}_h \cdot {\vVec}_\elem, \label{eq:dg_resT_motion} \\
r^{\text{ic}}_\elem((\ub_h, p_h); q_\elem) &\coloneqq  \intE (\det(\Fd^{k}_\elem(\ub_h)) - 1) \, q_\elem 
+ \sum_{F \in \FT^\internal}  \int_\face \eta_{\text{LBB}} \, h_\face \, \jump{p_h}_{\elem \face}  \, q_\elem, \label{eq:dg_resT_j1}\\
r^{\dirlm}_\face \left((\ub_h, \bm{\widehat{\lambda}}_h); {\widehat{\sVec}}_\face \right) &\coloneqq \intf \left( \ub_h - \bm{g}_\dir \right) \cdot \widehat{\sVec}_\face 
+ \sum_{\edge \in \EF^i} \int_\edge \eta_{\bm{\widehat{\lambda}}} \, h_\edge \, \jump{\bm{\widehat{\lambda}}_h}_{\face \edge} \cdot \widehat{\sVec}_\face. \label{eq:dg_resT_dir}
\end{align}
The average operator in \Eq~\eqref{eq:dg_resT_motion} is such that, for all $\varphi\in H^1(\Th)$ and all $F\in\Fh$,
\[
\mean{\varphi}_F \coloneqq \begin{cases}
\frac12\left(\varphi_{|\elem} + \varphi_{|\elem'}\right) & \text{if $\face\in\FT^\internal\cap\mathcal{\face}_{\elem'}^\internal$ with $\elem,\elem'\in\Th$, $\elem\neq \elem'$},
\\
\varphi_{|\face} & \text{otherwise}.
\end{cases}
\]
with the understanding that $\mean{\bullet}_\face$ acts component-wise when applied to vector and tensor functions. 
Furthermore, for any $\edge \in \EF^{\internal}$, the edge jump operator in \Eq~\eqref{eq:dg_resT_dir} reads
\begin{equation*}
	\jump{\vVec}_{FE} \coloneqq \vVec_{|\face} - \vVec_{|\face'}  \quad \text{if $\edge \in \EF \cap \mathcal{\edge}_{\face'}$ with $\face,\face'\in\Fh^\dirlm$, $\face \neq \face'$}.
\end{equation*}
With $h_\face$ and $h_\edge$, we denote the diameter of a face or an edge, respectively. Note that $h_\edge = 1$ if $d=2$. 

The local residuals in \eqref{eq:dg_resT_motion}-\eqref{eq:dg_resT_dir} contain several stabilization terms 
and user-dependent stabilization parameters $\eta_{(\bullet)}>0$ dictated by stability requirements. 
The first term in the second line of \eqref{eq:dg_resT_motion} is a stabilization term inspired by \cite{BrezziManziniEA99}, where it was 
introduced to ensure coercivity of the BR1 formulation proposed by Bassi and Rebay \cite{Bassi1997}, see also \cite{Arnold2002}.
The adaptive stabilization parameter $\eta_F$ is computed as proposed by Eyck and co-workers ~\cite{TenEyck2006,Eyck2008},
who first introduced the idea of adaptive stabilization in the context of BR1 dG discretizations of nonlinear elasticity problems.
A comprehensive description of the procedure involved in the computation of $\eta_F$ will be given in \Sec~\ref{sec:adaptive_stabilization_strategy}.
The last term in \eqref{eq:dg_resT_j1}, endowed with stabilization parameter $\eta_{\text{LBB}}$, ensures LBB stability 
by penalizing the pressure jumps across internal faces, see \cite{DiPietroStokesArtComp07} and \cite{Baroli2013}. 
The second term in \eqref{eq:dg_resT_dir}, endowed with stabilization parameter $\eta_{\bm{\widehat{\lambda}}}$, 
penalizes the Lagrange multipliers jumps across \emph{internal edges} of Dirichlet boundaries.
We remark that, for the sake of consistency of the dG formulation, we do not penalize the jumps over the sharp corners $\Gamma^i_h$.
Indeed, since $\bm{\widehat{\lambda}}_h$ approximates the stress vector $\Piola \N$ and the normal vector $\N$ is discontinuous at the sharp corners $\Gamma^i_h$, 
Lagrange multipliers are discontinuous over each $\edge \in \EF^b \cap \Gamma^i_h$.
Definitions of $\EF^\internal$, $\EF^b$ and $\Gamma^i_h$ are given in Sec. \ref{subsec:meshsetlm}.  

The global residuals $r_h^{\text{Lem}}\left( (\ub_h,p_h,\widehat{\bm{\lambda}}_h);\bullet\right) : \Poly{k}(\Th)^d \to \mathbb{R}$ 
and $r_h^{\text{ic}} \left( (\ub_h,p_h);\bullet\right) : \Poly{k}(\Th) \to \mathbb{R}$ are obtained 
assembling \textit{element-by-element} the local residuals  \eqref{eq:dg_resT_motion} and \eqref{eq:dg_resT_j1}, i.e. 
\begin{equation*}
\label{eq:resh}
{r}^{\text{Lem}}_h\left( (\ub_h, p_h, \widehat{\bm{\lambda}}_h); {\vVec}_h \right) \coloneqq \sum_{\elem \in \Th} {r}^{\text{Lem}}_\elem\left( (\ub_h, p_h, \widehat{\bm{\lambda}}_h); {\vVec}_h{}_{|\elem} \right),\qquad
r^{\text{ic}}_h\left((\ub_h, p_h); q_h \right) \coloneqq \sum_{\elem \in \Th} r^{\text{ic}}_\elem\left((\ub_h, p_h); q_h{}_{|\elem}\right).
\end{equation*}
Similarly, the global residual $r_h^{\dirlm}\left( (\ub_h,\widehat{\bm{\lambda}}_h);\bullet\right) : \Poly{k}(\Fh^\dirlm)^d \to \mathbb{R}$
is obtained assembling \textit{face-by-face} the local residual \eqref{eq:dg_resT_dir}, i.e. 
\begin{equation*}
\label{eq:reshF}
r_h^{\dirlm}\left( (\ub_h,\widehat{\bm{\lambda}}_h); \widehat{\sVec}_h \right) \coloneqq \sum_{\face \in \Fh^\dirlm} r_\face^{\dirlm}\left( (\ub_h,\widehat{\bm{\lambda}}_h); \widehat{\sVec}_h{}_{|\face} \right).
\end{equation*}


Defining, for the sake of brevity
$$ \W^{k}_{\text{ie},h} = \Poly{k}(\Th)^d \times \Poly{k}(\Th) \times \Poly{k}(\Fh^\dirlm)^d, \qquad \W^{k}_{\text{ce},h} = \Poly{k}(\Th)^d \times \Poly{k}(\Fh^\dirlm)^d.$$
The discrete nonlinear elasticity problems reads as follows: 
\begin{itemize}
	\item \textbf{Incompressible material}: find $ (\ub_h, p_h, \bm{\widehat{\lambda}}_h) \in \W^{k}_{ie,h}$ such that 
	\begin{equation}
	r_{\text{ie},h}\left( (\ub_h,p_h,\bm{\widehat{\lambda}}_h); (\vVec_h,q_h,\widehat{\sVec}_h)\right) =  0  \qquad \forall (\vVec_h,q_h,\widehat{\sVec}_h) \in \W^{k}_{ie,h},
	\label{eq:incompressible_elast_prob}
	\end{equation}
	where, given  $(\ub_h, p_h,\widehat{\sVec}_h) \in \W^{k}_{ie,h}$, $r_{\text{ie},h} \left( (\ub_h,p_h,\widehat{\bm{\lambda}}_h); \bullet \right) : \W^{k}_{ie,h} \rightarrow \mathbb{R}$ is such that $\forall (\vVec_h,q_h,\widehat{\sVec}_h) \in \W^{k}_{ie,h}$
	\begin{equation}
	r_{\text{ie},h} \left( (\ub_h,p_h,\bm{\widehat{\lambda}}_h); (\vVec_h,q_h,\widehat{\sVec}_h)\right) = r_h^{\text{Lem}}\left( (\ub_h,p_h,\bm{\widehat{\lambda}}_h); {\vVec_h}\right) + r_h^{\text{ic}}\left( (\ub_h,p_h); {q_h}\right) + r^{\dirlm}_h((\ub_h, \bm{\widehat{\lambda}}_h); \widehat{\sVec}_h);
	\label{eq:global_residual_incompressible}
	\end{equation}
	\item \textbf{Compressible material}: find $(\ub_h, \bm{\widehat{\lambda}}_h) \in \W^{k}_{ce,h}$ such that
	\begin{equation}
	r_{\text{ce},h} \left( (\ub_h, \bm{\widehat{\lambda}}_h); (\vVec_h, \widehat{\sVec}_h) \right) = 0  \qquad \forall (\vVec_h, \widehat{\sVec}_h) \in \W^{k}_{ce,h},
	\label{eq:compressible_elast_prob}
	\end{equation}
	where, given  $(\ub_h,\bm{\widehat{\lambda}}_h) \in \W^{k}_{ce,h}$, $r_{\text{ce},h} \left( (\ub_h, \bm{\widehat{\lambda}}_h); \bullet \right) : \W^{k}_{ce,h} \rightarrow \mathbb{R}$ is such that $\forall (\vVec_h, \widehat{\sVec}_h) \in \W^{k}_{ce,h}$
	\begin{equation}
	r_{\text{ce},h}\left( \ub_h; \vVec_h\right) = r_h^{\text{Lem}}\left( (\ub_h,0,\bm{\widehat{\lambda}}_h); (\vVec_h, \widehat{\sVec}_h)\right) + r^{\dirlm}_h((\ub_h, \bm{\widehat{\lambda}}_h); \widehat{\sVec}_h). 
	\label{eq:global_residual_compressible}
	\end{equation}
\end{itemize}

\subsection{Adaptive stabilization}\label{sec:adaptive_stabilization_strategy}
If, on the one hand, the amount of stabilization required to ensure the coercivity of 
dG formulations can be precisely estimated in the context of linear elasticity problems, see \cite{Lew2004} and \cite{Wihler2004}, 
on the other hand, penalty parameters are not known a priori in the case of finite deformations of hyperelastic materials, see \cite{Eyck2008a} and \cite{Cockburn2019}.
This is an uncomfortable situation as excessive stabilization worsen the condition number 
of system matrices and insufficient stabilization severely affects the robustness of numerical schemes.
As a result, a tedious trial and error approach would often be required in practice. 
The adaptive stabilization strategies proposed and analysed by Eyck and co-workers \cite{TenEyck2006}-\cite{Eyck2008}
in the context of BR1 dG formulations are a crucial tool for mitigating this drawback. 

The adaptive stabilization  
introduced in \eqref{eq:dg_resT_motion},
can be considered the natural extension of the approach proposed in \cite{Eyck2008} to the BR2 dG discretization.
The stabilization parameter $\eta_{\face}$ is defined as follows
$$
\eta_{\face} = \epsilon + \beta \; \lambda_{\face} 
$$
where $\epsilon,\beta \geq 0$ are user-dependent parameters. 
Let a mesh element $\elem \in \Th$ and a mesh face $\face \in \FT^{\internal, \dirw}$ be fixed,
for each point $\X \in \face $
$$
{ \lambda}_{\elem \face}({\X}) = \max \left\{ 0, - \min_{\bm{0}\neq \bm{G} \in \mathbb{R}^{d \times d}} \frac{\bm{G} : \Tnsr{A}(\Fd_{\elem\face}^{k+1}(\ub_h(\X)),p_h) : \bm{G} }  {\bm{G} : \bm{G} } \right\} 
$$ 
and
$$
\lambda_{\face} =
\begin{cases}
\frac{\overline{ \lambda}_{{\elem\face}} + \overline{ \lambda}_{{\elem' \face}}}{2} & \text{if $\face \in \FT^\internal\cap\mathcal{\face}_{\elem'}^\internal$ with $\elem,\elem'\in\Th$, $\elem \neq \elem'$}, \\
\overline{ \lambda}_{{\elem\face}} & \text{if $\face\in\FT^\dirw$}.
\end{cases}
$$
where $\Tnsr{A}$ is the fourth-order elasticity tensor described in \Sec~\ref{sec:elastic_fourth-order_tensor} and $\overline{\lambda}_{{\elem\face}}$ is the mean value of $\lambda_{{\elem\face}}(\X)$ over $\face$.	

As detailed in Itskov~\cite{Itskov2000}, the nine (in 3D, six in 2D) eigenvalues $\lambda_i$ of $\Tnsr{A}$ 
can be obtained solving the characteristic equation
\begin{equation}
\det \left( \Tnsr{A} - \lambda \Tnsr{I} \right) = 0
\label{eq:eigenvalue_problem_matrix_form}
\end{equation}
where
$\Tnsr{I} = \delta_{ik} \delta_{jl} \bm{e}_i \otimes \bm{e}_j \otimes \bm{e}_k \otimes \bm{e}_l$ is the fourth-order identity tensor. 
Problem \eqref{eq:eigenvalue_problem_matrix_form} is solved for each quadrature point of each mesh face at the first step of each Newton iteration involved 
in the globalisation strategy presented in the next section.

\section{Incremental load method}\label{sec:incremental_method}
Problems \eqref{eq:incompressible_elast_prob} and \eqref{eq:compressible_elast_prob} are solved by means of Newton's method.
In order to globalize the convergence of the Newton iteration towards the equilibrium configuration, we adopt the \textit{incremental method}, see \eg{} \cite{Ciarlet1999a}. 
The idea is to define a quasi-static loading path that allows to reach the final configuration by passing through a sequence of intermediate equilibrium states.
To this end, an incremental percentage of the external solicitation is imposed at each intermediate step $i=1,...N{-}1$, where $N$ is the step corresponding to the final configuration.
In particular, for $i=1,2...,N$, we apply Newton's method to solve the following problem: find $\wVec_{h}^i \in \W^k_{\bullet,h}$ such that 
\begin{equation}
\label{eq:elastProbIncr}
\tilde{r}_{\bullet,h}^i(\wVec_h^i; \zVec_h) = 0, \quad \forall \zVec_h \in \W^k_{\bullet,h}
\end{equation}
where $\tilde{r}^i_{\bullet,h}(\ast;\ast)$ is one of the residuals defined in \eqref{eq:global_residual_incompressible}-\eqref{eq:global_residual_compressible},
\ie{} $\bullet \in \{\text{ce,ie}\}$,
but the external solicitations in \eqref{eq:dg_resT_motion}, see also \eqref{eq:jumpDG},
are replaced by
\begin{equation}
\label{eq:incrForcing}
\tilde{\bm{f}}^{+} = \frac{i}{N}\bm{f}^{+}, \qquad \tilde{\bm{g}}_\neu= \frac{i}{N}\bm{g}_\neu, \qquad \tilde{\bm{g}}_\dir = \frac{i}{N}\bm{g}_\dir.
\end{equation}
Clearly, $\wVec_h^{N}$ is the solution of one of the problems in \eqref{eq:incompressible_elast_prob}-\eqref{eq:compressible_elast_prob}.
Nevertheless, since each state of the sequence is incremental with respect to the previous configuration, convergence of Newton's method is guaranteed providing $N$ big enough.  

Newton's method applied to problem \eqref{eq:elastProbIncr} reads: 
\begin{align*}
&\text{set the initial guess $\wVec_h^i = \wVec_h^{i-1}$}, \nonumber \\
&\text{while $\delta \wVec_h$ is too large, find $\delta \wVec_h \in \W^k_{\bullet,h}$ such that} \\
&\quad \left( \tnsr{J}_{\bullet, h}(\wVec_h^i) \delta \wVec_h, \zVec_h \right)_{L^2(\Omega)}= -\tilde{r}_{\bullet,h}^i(\wVec_h^i; \zVec_h),  \quad \forall \zVec_h \in \W^k_{\bullet,h} \nonumber,\\
&\quad \text{set} \; \wVec_h^i \mathrel{+}= \delta \wVec_h, \nonumber 
\end{align*}
where, for each $\wVec_h \in \W^k_{\bullet,h}$, the Jacobian operator $\tnsr{J}_{\bullet, h}: \W^k_{\bullet,h} \rightarrow \W^k_{\bullet,h}$ is defined such that 
$$ \left( \tnsr{J}_{\bullet, h}(\wVec_h) \yVec_h, \zVec_h \right)_{L^2(\Omega)} = \frac{d}{d \epsilon } \tilde{r}_{\bullet, h}(\wVec_h + \epsilon \yVec_h; \zVec_h)\bigg|_{\epsilon=0}, \qquad \forall \yVec_h, \zVec_h \in \W^k_{\bullet,h}. $$

According to the definitions in \Sec~\ref{sec:elastic_fourth-order_tensor} and \Sec~\ref{subsec:dg_formulation}, the Jacobian operators reads as follows:
\begin{itemize}
	\item \textbf{Incompressible material}: $\left( \tnsr{J}_{\text{ie}, h}(\wVec_h) \delta \wVec_h,  \zVec_h \right)_{L^2(\Omega)} =$ 
	\begin{alignat*}{3}
	&=&\sumOnElem &\bigg( \intE  \left[ \boldsymbol{\mathfrak{G}}_T^k(\delta \uVec_h) : \Tnsr{A} (\Fd^k_{\elem}(\uVec_h), p_h) \right] : \nabla {\vVec}_\elem +\nonumber \\
	&&&{-}\sumOnFEID \intf \left[ \meanBig{  \left( \nabla \delta \uVec_{h|\elem} 
		{-} \boldsymbol{\mathfrak{R}}_{\face \elem}^{k+1}(\jump{\delta \uVec_h}_{\elem \face})\right)  :  \Tnsr{A} (\Fd_{\elem\face}^{k+1}(\uVec_h), p_h) } \, \normal_{\elem\face} \right] \cdot {\vVec}_\elem + \nonumber \\
	&&&{+}\sumOnFEID \intf  \eta_{\face} \, \left[ \meanBig{ \boldsymbol{\mathfrak{R}}_{\face \elem}^k(\jump{\delta \uVec_h}_{\elem \face}) } \normal_{\elem\face} \right] \cdot {\vVec}_\elem +\nonumber \\
	&&&{-}\intE  \delta p_h \, \det(\Fd^k_\elem (\uVec_h) ) \, [\Fd^k_\elem(\uVec_h) ]^{-\intercal}: \nabla {\vVec}_\elem +\nonumber \\
	&&&{-}\sumOnFEID \intf \left[ \meanBig{ \delta p_h \, \det(\Fd_{\elem\face}^{k+1}(\uVec_h) ) \, [\Fd_{\elem\face}^{k+1}(\uVec_h) ]^{-\intercal}} \, \normal_{\elem\face} \right] \cdot {\vVec}_\elem + \nonumber\\
	&&&{+} \intE \left( \det(\Fd^k_\elem (\uVec_h )) [\Fd^k_\elem(\uVec_h )]^{-\intercal} : \boldsymbol{\mathfrak{G}}_T^k(\delta \uVec_h) \right) \, q_\elem  \nonumber 
	{+}\sum_{F \in \FT^\internal} \int_\face \eta_{\text{LBB}} \, h_\face \, \jump{\delta p_h  }_{\elem \face}  \, q_\elem \bigg) +\nonumber \\
	&+ &\sum_{\face \in \Fh^\dirlm}& \bigg({-} \intf \delta\bm{\widehat{\lambda}}_h \cdot {\vVec}_\elem + \intf  \delta \ub_h \cdot \widehat{\sVec}_\face + \sum_{\edge \in \EF^i} \int_\edge \eta_{\bm{\widehat{\lambda}}} \, h_\edge \, \jump{\delta \bm{\widehat{\lambda}}_h}_{\face \edge} \cdot \widehat{\sVec}_\face \bigg);  
	\end{alignat*}
	
	\item \textbf{Compressible material}: $\left( \tnsr{J}_{\text{ce}, h}(\wVec_h) \delta \wVec_h , \zVec_h \right)_{L^2(\Omega)} = $ 
	\begin{alignat*}{3}
	&= &\sumOnElem &\bigg( 
	\intE  \left[ \boldsymbol{\mathfrak{G}}_T^k(\delta \uVec_h) : \Tnsr{A} (\Fd^k_{\elem}(\uVec_h), 0) \right] : \nabla {\vVec}_\elem + \nonumber \\
	& & &{-}\sumOnFEID \intf \left[ \meanBig{  \left( \nabla \delta \uVec_{h|\elem} - \boldsymbol{\mathfrak{R}}_{\face \elem}^{k+1}(\jump{\delta \uVec_h}_{\elem \face})\right)  :  \Tnsr{A} (\Fd_{\elem\face}^{k+1}(\uVec_h), 0) } \, \normal_{\elem\face} \right] \cdot {\vVec}_\elem + \nonumber \\
	& & & +\sumOnFEID \intf  \eta_{\face} \; \meanBig{ \boldsymbol{\mathfrak{R}}_{\face \elem}^k(\jump{\delta \uVec_h}_{\elem \face}) \, \normal_{\elem\face} } \cdot {\vVec}_\elem
	\bigg) + \nonumber \\
	&+ &\sum_{\face \in \Fh^\dirlm} &\bigg({-} \intf \delta\bm{\widehat{\lambda}}_h \cdot {\vVec}_\elem + \intf  \delta \ub_h \cdot \widehat{\sVec}_\face + \sum_{\edge \in \EF^i} \int_\edge \eta_{\bm{\widehat{\lambda}}} \, h_\edge \, \jump{\delta \bm{\widehat{\lambda}}_h}_{\face \edge} \cdot \widehat{\sVec}_\face \bigg). 
	\end{alignat*}
\end{itemize}

\section{Numerical results}\label{sec:numerical_results}
In this section, we numerically validate the BR2 dG discretizations of \Sec~\ref{sec:dg_formulation} solving compressible and incompressible 
nonlinear elasticity problems.
As a first point, we verify the numerical convergence rates for each of the constitutive laws in \eqref{eq:svk_energy}-\eqref{eq:neoh_inc_energy} 
based on manufactured 3D solutions.
To this end, the $L^2$ error norms of the displacement and the displacement gradient are tabulated 
varying the mesh size $h$ and the polynomial degree $k$. 
Afterwords, we challenge the stabilization strategy performing three specifically conceived 2D computations:
the \textit{parabolic indentation problem}, see \Sec~\ref{sec:indentation_problem}, the \textit{beam deformation}, 
see \Sec~\ref{sec:beam_deformation} and the \textit{cavitating voids}, see \Sec~\ref{sec:cavitating_voids}. 
Notice that the constitutive law in \Eq~\eqref{eq:neoh_energy_cavitation} is not included in numerical convergence tests cases 
and is employed solely for the cavitating voids test case of \Sec~\ref{sec:cavitating_voids}. 
To conclude, we tackle 3D computations and analyse the robustness of the $h$-multigrid solution strategy with respect to the stabilization parameter.
\Sec~\ref{sec:beam_torsion} and \Sec~\ref{sec:cylinder_rotation} consider the \textit{torsion of a square-section bar} 
and the \textit{deformation of a hollow cylinder} subjected to the rotation of its top surface, respectively.\\

All numerical test cases require the setup of the incremental strategy presented in \Sec~\ref{sec:incremental_method} whose crucial parameter is the number of loading steps.
On the one hand, an underestimated number of increments leads to Newton's method convergence failure and breakdown of the solution strategy, on the other hand, 
an overestimated number of increments causes an excessive computational cost.
Converge failure is often associated with $\det (\Fd_\elem^k(\ub)) \leq 0$ or $\det (\Fd_{\elem\face}^{k+1}(\ub)) \leq 0$, meaning that the fundamental hypothesis stated in \Eq~\eqref{eq:jacMap} is violated. 
While the number of increments can be adaptively chosen by splitting the problematic step until \Eq~\eqref{eq:jacMap} is satisfied (almost) everywhere in the domain,
we rely on equispaced increments in all the numerical test cases presented hereafter. 
The goal is to show how different test cases are handled in terms of number of increments and 
to stress the crucial role of the strategy employed for imposing Dirichlet boundary conditions, see \Sec~\ref{sec:IncrDBCs}.

At each successful loading path step, Newton iteration achieves a relative residual decrease 
of ten orders of magnitude in less than eight iterations (usually between four and six).
The sequence of linearised equation systems can be solved with either a direct or an iterative solver. 
In the latter case, due to poor performance of standard Incomplete Lower Upper (ILU) factorization preconditioners,
we adopt the $h$-multigrid agglomeration based solution strategy proposed in Botti \ea~\cite{Botti2017}. 
As a distinctive feature, $h$-coarsened mesh sequences are generated on the fly by recursive agglomeration of the fine grid
and, accordingly, arbitrarily unstructured grids can be handled as an input of the agglomeration strategy. 
The cost of numerical integration over agglomerated elements is mitigated by using element-by-element $L^2$ projections 
to build coarse grid operators, with projection operators computed and stored once-and-for-all in a preprocessing phase. 
The performance of the $h$-multigrid preconditioned iterative solver will be evaluated in terms of number of iterations 
required to reach a eight orders of magnitude drop of the relative residual norm. 
We remark that the $h$-multigrid solution strategy has not been implemented and tested in combination 
with boundary conditions enforced by means of Lagrange multipliers, accordingly a direct solver is employed instead.

\subsection{Influence of Dirichlet boundary conditions on the incremental load method}
\label{sec:IncrDBCs}
The BR2 formulations of \Sec~\ref{sec:dg_formulation} admit the imposition of Dirichlet Boundary Conditions (BCs) by means of Nitsche method and Lagrange multipliers method. 
Interestingly, the strategy based on Lagrange multipliers is the most effective, leading to increased robustness of the incremental load method. 
The following reasoning provides an intuitive explanation for the aforementioned behavior.
When using Nitsche method for Dirichlet BCs the occurrence of null or negative Jacobian values, namely $\det (\Fd_{\elem\face}^{k+1}(\ub)) \leq 0$, is often 
triggered by the action of lifting operators on Dirichlet boundaries, see \Eq~\eqref{eq:gradDG} and definition \eqref{eq:jumpDG}. 
On internal faces, the jumps magnitude is controlled mainly by the discretization parameters: 
in particular we expect the jumps to shrink while increasing the polynomial degree $k$ and decreasing the mesh step size $h$.
As opposite, on Dirichlet boundary faces, since Newton's method initial guess is the solution of the previous incremental step, 
jumps magnitude is dictated primarily by the number of increments of the incremental load method, see also \eqref{eq:incrForcing}.

\begin{table}[ht]
	\centering
		\begin{tabular}{|l|c|c|c|c|c|}
			\cline{3-6}
			\multicolumn{2}{l}{\multirow{2}{*}{}}         & \multicolumn{4}{|c|}{\textbf{Number of increments in the loading path}} \\ 
                        \cline{3-6}
		    \multicolumn{2}{l}{}                          & \multicolumn{3}{|c|}{\makecell{\textbf{BCs enforced by} \\ \textbf{Nitsche method}}} 
                                                                  & \makecell{\textbf{BCs enforced by} \\ \textbf{Lagrange multipliers}}\\ \hline 
			\textbf{Model}             & \textbf{Mesh}   & $k=1$         & $k=2$         &  $k=3$    & $k=1,2,3 $       \\ \hline
			\multirow{2}{*}{NHK-C}     & coarse          & 100           & 400           & 400       & \multirow{2}{*}{3}  \\ \cline{2-5} 
			                           & fine            & 400           & 800           & 800       &   \\ \hline
			\multirow{2}{*}{NHK-I }    & coarse          & 400           & 400           & 400       & \multirow{2}{*}{-}  \\ \cline{2-5} 
			                           & fine            & 800           & 1000          & 1500      &   \\ \hline
			\multirow{2}{*}{SVK-C}     & coarse          & 400           & 800           & 800       & \multirow{2}{*}{3} \\ \cline{2-5} 
			                           & fine            & 800           & 2000          & 2500      &   \\ \hline
			\multirow{2}{*}{SVK-I}     & coarse          & 400           & 400           & 800       & \multirow{2}{*}{-}  \\ \cline{2-5} 
		                               & fine            & 1500          & 2000          & 2000      &   \\ \hline
		\end{tabular}
	\caption{Number of increments in the loading path for the manufactured solutions of \Sec~\ref{sec:convergence_tests}.
                 We consider any combination of neo-Hookean and Saint Venant-Kirchhoff constitutive models with compressible and incompressible materials. 
                 Results are given considering the coarsest and the finest grid of the Cartesian grids sequence ($4^3$ and $32^3$ hexahedral elements, respectively) 
                 for different polynomial degrees $k=\left\lbrace 1,2,3 \right\rbrace $.
	\label{tab:increments_conv_test}}
\end{table}

\begin{table}[ht]
	\centering
	\begin{tabular}{|l|c|c|c|}
		\cline{3-4}
	     \multicolumn{2}{c}{} & \multicolumn{2}{|c|}{\textbf{Number of increments in the loading path}}\\ \hline
		\textbf{Test case} &   \textbf{Model} & \makecell{\textbf{BCs enforced by} \\ \textbf{Nitsche method}} 
                                   & \makecell{\textbf{BCs enforced by} \\ \textbf{Lagrange multipliers}} \\ \hline
		\multirow{2}{*}{Parabolic indentation}        &   NHK-C          & 60           & 2         \\ \cline{2-4}
	                                                      &   NHK-I / SVK-I  & 40           & -         \\ \hline
		Beam deformation                              &   NHK-C          & 600          & 15        \\ \hline
		Cavitating voids                              &   NHK-CAV        & -            & 100       \\ \hline
		\multirow{2}{*}{Bar torsion}                  &   NHK-C / SVK-C  & 15 -- 60     & 10        \\ \cline{2-4}
		                                              &   NHK-I / SVK-I  & 15 -- 80     & -         \\ \hline
		\multirow{2}{*}{Cylinder top face rotation}   &   NHK-C          & 1000         & 30        \\ \cline{2-4}
		                                              &   NHK-I          & 650          & -         \\ \hline

	\end{tabular}
	\caption{Number of increments in the loading path for all 2D and 3D test cases considered in {\Sec}s~\ref{sec:2D_simulations} and \ref{sec:3D_simulations}.
                 In case of Bar torsion the number of incremental steps was fine tuned according to the polynomial degree, see \Sec~\ref{sec:beam_torsion} for additional details.
	\label{tab:increments_simulations}}
\end{table}

The number of incremental steps for BR2 formulations with Nitsche method and Lagrange multipliers method Dirichlet BCS 
are reported in \Tab~\ref{tab:increments_conv_test}, where manufactured solutions are considered, and \Tab~\ref{tab:increments_simulations},
where we tabulate data for realistic test cases.
In \Tab~\ref{tab:increments_conv_test}, it is possible to appreciate that, in case of Nitsche method, 
increasingly higher step counts are required as the mesh is refined and the polynomial degree increases.
As opposite, in case of Lagrange multipliers method, three loading steps are employed irrespectively of discretization parameters.
The results of \Tab~\ref{tab:increments_simulations} confirm that Lagrange multipliers method leads to an astonishing decrease 
of the number of steps in all test cases. 
We remark that the Lagrange multiplier method has not been tested in the incompressible regime, further investigation will be carried in future works.

Finally, in \Tab~\ref{tab:stabilization_simulations}, we provide an overview of the stabilization parameters settings for each of the test cases 
presented in {\Sec}s~\ref{sec:2D_simulations} and \ref{sec:3D_simulations}.

\begin{table}[ht]
	\centering
	\begin{tabular}{|l|c|c|c|c|}
		\hline
		\textbf{Test case}         & \textbf{$\beta$} & \textbf{$\epsilon$} & \textbf{$\eta_{\text{LBB}}$} & \textbf{$\eta_{\widehat{\bm{\lambda}}}$} \\ \hline
		Convergence tests          & 1                & 0                   & 1                            & 1                                        \\ \hline
		Parabolic indentation      & 0                & 0                   & 1                            & 1                                        \\ \hline
		Beam deformation           & 1                & 1                   & -                            & 1                                        \\ \hline
		Cavitating voids           & 1                & 1                   & -                            & 1                                        \\ \hline
		Bar torsion                & 0--5             & 0                   & 1                            & 1                                        \\ \hline
		Cylinder top face rotation & 4                & 1                   & 1                            & 1                                        \\ \hline
	\end{tabular}
	\caption{Stabilization parameters for all 2D and 3D test cases considered in {\Sec}s~\ref{sec:2D_simulations} and \ref{sec:3D_simulations}.
                 In case of Bar torsion $\beta$ was fine tuned according to the polynomial degree, see \Sec~\ref{sec:beam_torsion} for additional details.
	\label{tab:stabilization_simulations}}
\end{table}

\subsection{Evaluation of convergence rates}\label{sec:convergence_tests}
Convergence tests consider the neo-Hookean (NHK) and Saint Venant-Kirchhoff (SVK) constitutive models in both the compressible (-C) and the incompressible (-I) regime.
Numerical solutions are obtained over a four grids $h$-refined mesh sequence of the unit cube $\Omega: [0,1]^3$.
The uniform hexahedral elements have diameter $h$ ranging from 0.25 (coarse mesh) to 0.03125 (fine mesh), halving $h$ at each refinement step. 
We apply first, second and third degree BR2 dG discretizations and enforce boundary conditions based on smooth analytical displacement fields, 
see \Sec~\ref{sec:SolComprMat} and \Sec~\ref{sec:SolIncomprMat}. 
Dirichlet boundary conditions based on the exact displacement are imposed with Nitsche method or Lagrange multipliers method on five of the six surfaces composing $\partial \Omega$.
A Neumann boundary condition based on the exact deformation gradient is imposed on the unaccounted surface.
Convergence is evaluated based on the $L^2$-norm of the error on the displacement, the displacement gradient and, eventually, the pressure. 
Forcing terms are computed based on analytical solutions by means of SageMath~\cite{sagemath}, an open-source library featuring symbolic calculus.

\subsubsection{Compressible materials}
\label{sec:SolComprMat}
Let's denote by $\X = (X,Y,Z)$ the Cartesian coordinates in the reference configuration and by $u,v,w$ the three components of the displacement vector $\ub$.\\
In the compressible regime, we consider the following displacement field proposed by Abbas \ea~\cite{Abbas2018}
\begin{equation*}
\begin{dcases}
u(\X) = \left( \frac{1}{\lambda }  + \alpha \right) X + \psi(Y)\\
v(\X) = -\left( \frac{1}{\lambda} + \dfrac{\alpha + \gamma + \alpha \gamma}{1+\alpha + \gamma + \alpha \gamma}\right) Y \\
w(\X) = \left( \frac{1}{\lambda} + \gamma \right) Z + \omega(X) + \xi (Y) 
\end{dcases}
\end{equation*}
where $\alpha = \gamma = 0.1$, $\psi(Y) = \alpha \sin(\pi Y)$, $\omega(X) = \gamma \sin(\pi X)$ and $\xi(Y) = 0$. 
Relevant parameters of the NHK-C and SVK-C constitutive laws are defined setting $\mu = 1$ and $\lambda = 10$, which corresponds to a Poisson's ration of $\nu \simeq 0.455$.
The adaptive stabilization parameters are taken as $\beta = 1 $ and $\epsilon = 0$, respectively.
In case of Dirichlet boundary conditions enforced by means of Lagrange multipliers, we set $\eta_{\bm{\widehat{\lambda}}} = 1$.
Asymptotic convergence rates of order $k+1$ and $k$ for the displacement and the displacement gradient 
can be appreciated in \Tab~\ref{tab:conv_result_neo_comp_3D} and \Tab~\ref{tab:conv_result_svk_comp_3D} for the NHK-C and the SVK-C model, respectively.

\begin{table}[ht]
        \centering
       \small
       \begin{tabular}{|c|cccc|cccc|}\hline
       	$\mathrm{card}(\Th)$ 
      & $\| \ub {-} \ub_h \|_{L^2(\Omega)}$ & \textbf{rate}
      & $\| \nabla(\ub {-} \ub_h) \|_{L^2(\Omega)}$ & \textbf{rate}
      & $\| \ub {-} \ub_h \|_{L^2(\Omega)}$ & \textbf{rate}
      & $\| \nabla(\ub {-} \ub_h) \|_{L^2(\Omega)}$ &\textbf{rate}\\ 
       	\hline                                                                                            
       	&\multicolumn{4}{c|}{$k=1$, BCs by Nitsche method}                     &  \multicolumn{4}{c|}{$k=1$, BCs by Lagrange multipliers method}                                                    \\ 
       	\hline           
       	64    & 2.457e-03 &  -   & 7.116e-02 &  -   & 2.622e-03 &  -   & 7.250e-02 &  -   \\
       	512   & 6.104e-04 & 2.00 & 3.561e-02 & 0.99 & 6.126e-04 & 2.10 & 3.565e-02 & 1.02 \\ 
       	4096  & 1.536e-04 & 1.99 & 1.781e-02 & 0.99 & 1.538e-04 & 1.99 & 1.780e-02 & 1.00 \\ 
       	32768 & 3.853e-05 & 1.99 & 8.903e-03 & 1.00 & 3.858e-05 & 1.99 & 8.904e-03 & 1.00 \\ 
       	\hline                                   
       	&\multicolumn{4}{c|}{$k=2$, BCs by Nitsche method}                     & \multicolumn{4}{c|}{$k=2$, BCs by Lagrange multipliers method}         \\ 
       	\hline                                  
       	64    & 2.241e-04 &  -   & 7.374e-03 &  -   & 6.589e-04 &  -   & 1.295e-02 &  -   \\ 
       	512   & 2.859e-05 & 2.97 & 1.838e-03 & 2.00 & 4.982e-05 & 3.73 & 2.323e-03 & 2.48 \\ 
       	4096  & 3.568e-06 & 3.00 & 4.580e-04 & 2.00 & 4.484e-06 & 3.47 & 4.944e-04 & 2.23 \\ 
       	32768 & 4.421e-07 & 3.01 & 1.141e-04 & 2.00 & 4.765e-07 & 3.23 & 1.165e-04 & 2.09 \\ 
       	\hline                                 
       	&\multicolumn{4}{c|}{$k=3$, BCs by Nitsche method}                     &\multicolumn{4}{c|}{$k=3$, BCs by Lagrange multipliers method}          \\ 
       	\hline                                
       	64    & 1.130e-05 &  -   & 4.926e-04 &  -   & 1.257e-05 &  -   & 5.400e-04 &  -   \\ 
       	512   & 7.462e-07 & 3.92 & 6.110e-05 & 3.01 & 7.480e-07 & 4.07 & 6.183e-05 & 3.13 \\ 
       	4096  & 4.796e-08 & 3.96 & 7.579e-06 & 3.00 & 4.794e-08 & 3.96 & 7.598e-06 & 3.02 \\ 
       	32768 & 3.043e-09 & 3.98 & 9.420e-07 & 3.00 & $\star$   &  -   & $\star$   &  -   \\ 
       	\hline
       \end{tabular}
	\caption{Errors and convergence rates for BR2 dG discretizations of degree $k=\left\lbrace 1,2,3\right\rbrace $ over a $h$-refined mesh sequence of the unit cube, NHK-C constitutive model.
         $\star$ indicates unavailable data due to excessive memory consumption of the LU solver. 
	\label{tab:conv_result_neo_comp_3D}}
\end{table}

\begin{table}[ht]
        \centering
       \small
       \begin{tabular}{|c|cccc|cccc|}\hline
       	$\mathrm{card}(\Th)$ 
      & $\| \ub {-} \ub_h \|_{L^2(\Omega)}$ & \textbf{rate}
      & $\| \nabla(\ub {-} \ub_h) \|_{L^2(\Omega)}$ & \textbf{rate}
      & $\| \ub {-} \ub_h \|_{L^2(\Omega)}$ & \textbf{rate}
      & $\| \nabla(\ub {-} \ub_h) \|_{L^2(\Omega)}$ &\textbf{rate}\\ 
       	\hline                                                                                            
       	&\multicolumn{4}{c|}{$k=1$, BCs by Nitsche method}                     &  \multicolumn{4}{c|}{$k=1$, BCs by Lagrange multipliers method}                                                    \\ 
       	\hline           
       	64    & 2.486e-03 &  -   & 7.117e-02 &  -   & 2.648e-03 &  -   & 7.248e-02 &  -    \\
       	512   & 6.211e-04 & 2.00 & 3.562e-02 & 0.99 & 6.259e-04 & 2.08 & 3.567e-02 & 1.02  \\ 
       	4096  & 1.558e-04 & 1.99 & 1.781e-02 & 1.00 & 1.562e-04 & 1.99 & 1.781e-02 & 1.00  \\ 
       	32768 & 3.901e-05 & 1.99 & 8.904e-03 & 1.00 & 3.906e-05 & 1.99 & 8.904e-03 & 1.00  \\ 
       	\hline                                   
       	&\multicolumn{4}{c|}{$k=2$, BCs by Nitsche method}                     & \multicolumn{4}{c|}{$k=2$, BCs by Lagrange multipliers method}         \\ 
       	\hline                                  
       	64    &  2.170e-04 &   -  & 7.345e-03 &  -   & 6.217e-04 &  -   & 1.262e-02 &  -   \\ 
       	512   &  2.789e-05 & 2.99 & 1.831e-03 & 2.00 & 4.810e-05 & 3.69 & 2.283e-03 & 2.47 \\ 
       	4096  &  3.504e-06 & 2.99 & 4.566e-04 & 2.00 & 4.401e-06 & 3.45 & 4.909e-04 & 2.22 \\ 
       	32768 &  4.383e-07 & 2.99 & 1.140e-04 & 2.00 & 4.731e-07 & 3.22 & 1.162e-04 & 2.08 \\ 
       	\hline                                 
       	&\multicolumn{4}{c|}{$k=3$, BCs by Nitsche method}                     &\multicolumn{4}{c|}{$k=3$, BCs by Lagrange multipliers method}          \\ 
       	\hline                                
       	64    & 1.138e-05 &  -   & 4.894e-04 &  -   & 1.261e-05 &  -   & 5.385e-04 &  -   \\ 
       	512   & 7.513e-07 & 3.92 & 6.069e-05 & 3.01 & 7.516e-07 & 4.07 & 6.186e-05 & 3.12 \\ 
       	4096  & 4.818e-08 & 3.96 & 7.542e-06 & 3.00 & 4.804e-08 & 3.97 & 7.600e-06 & 3.03 \\ 
       	32768 & 3.049e-09 & 3.98 & 9.396e-07 & 3.00 & $\star$   & -    & $\star$   &  -   \\ 
       	\hline
       \end{tabular}
	\caption{Errors and convergence rates for BR2 dG discretizations of degree $k=\left\lbrace 1,2,3\right\rbrace $ over a $h$-refined mesh sequence of the unit cube, SVK-C constitutive model.
         $\star$ indicates unavailable data due to excessive memory consumption of the LU solver. 
	\label{tab:conv_result_svk_comp_3D}}
\end{table}

\subsubsection{Incompressible nonlinear elasticity}
\label{sec:SolIncomprMat}
The fully incompressible nonlinear elasticity problem is defined
according to the following isochoric displacement field
\begin{equation*}
\begin{dcases}
u(\X)= (a^2 - 1) X + \frac{b}{2} \sin^2(Y) + \frac{c}{2} \sin^2(Z)\\
v(\X) = \left(\frac{1}{a} - 1 \right) Y \\
w(\X) = \left(\frac{1}{a} - 1 \right) Z \\
\end{dcases}
\label{eq:exact_displ_inc}
\end{equation*}
where $a=1.1$, $b=1$ and $c=1$.
The pressure field reads $p = \dfrac{1}{3} \tr{\tnsr{\sigma}}$, where $\tnsr{\sigma} = \dfrac{1}{J} \F \tnsr{S} \F^\intercal$ is the Cauchy stress tensor. 
Based on the NHK-I model, the exact pressure reads
\begin{equation*}
p^{\text{NHK-I}} = \frac{{{c}^{2}}\, {{\cos^{2}{(Z)}}}\, {{\sin^{2}{(Z)}}} \mu +{{b}^{2}}\, {{\cos^{2}{(Y)}}}\, {{\sin^{2}{(Y)}}} \mu +{{a}^{4}} \mu +\frac{2 \mu }{{{a}^{2}}}}{3}.
\end{equation*}
Based on the SVK-I model, the exact pressure reads
\begin{align*}
p^{\text{SVK}} = \frac{1}{6} \mu \bigg(& 2 {{\cos^{4}{(Z)}}}\, {{\sin^{4}{(Z)}}} + \Big( 4 {{\cos^{2}{(Y)}}}\, {{\sin^{2}{(Y)}}}+6\Big) \, {{\cos^{2}{(Z)}}}\, {{\sin^{2}{(Z)}}}+ \nonumber \\
& +2 {{\cos^{4}{(Y)}}}\, {{\sin^{4}{(Y)}}}+6 {{\cos^{2}{(Y)}}}\, {{\sin^{2}{(Y)}}} \bigg) + \nonumber\\
+   \frac{1}{6} \lambda \bigg(& {{\cos^{4}{(Z)}}}\, {{\sin^{4}{(Z)}}}+\Big( 2 {{\cos^{2}{(Y)}}}\, {{\sin^{2}{(Y)}}}+3\Big) \, {{\cos^{2}{(Z)}}}\, {{\sin^{2}{(Z)}}}+ \nonumber \\
&+{{\cos^{4}{(Y)}}}\, {{\sin^{4}{(Y)}}} +3 {{\cos^{2}{(Y)}}}\, {{\sin^{2}{(Y)}}}\bigg).
\label{eq:exact_pressure_inc} 
\end{align*}
As in the previous section, the adaptive stabilization parameters are set as $\beta = 1$ and $\epsilon = 0$,
while pressure jumps stabilization coefficient is taken as $\eta_{\text{LBB}} = 1$. 
Similarly to the compressible regime, asymptotic convergence rates of order $k+1$ and $k$ 
are observed for the displacement and the displacement gradient over $h$-refined meshes. 
The pressure error in $L^2$-norm exhibits a rate of convergence between $k$ and $k+1$.
Convergence results are reported in \Tab~\ref{tab:conv_result_neo_inc_3D} and \Tab~\ref{tab:conv_result_svk_inc_3D} for the NHK-I and SVK-I models, respectively.\\

\begin{table}[H]
        \centering
       \small
       \begin{tabular}{|c|cccccc|}\hline
       	$\mathrm{card}(\Th)$ 
      & $\| \ub {-} \ub_h \|_{L^2(\Omega)}$ & \textbf{rate}
      & $\| p {-} p_h \|_{L^2(\Omega)}$ & \textbf{rate}
      & $\| \nabla(\ub {-} \ub_h) \|_{L^2(\Omega)}$ & \textbf{rate} \\
       	\hline                                                                                            
       	&\multicolumn{6}{c|}{$k=1$, BCs by Nitsche method}  \\ 
       	\hline           
       	64    & 2.665e-03 &  -   & 3.137e-02 & -    & 6.589e-02 &  -   \\ 
       	512   & 6.531e-04 & 2.02 & 1.002e-02 & 1.64 & 3.270e-02 & 1.01 \\ 
       	4096  & 1.638e-04 & 1.99 & 2.764e-03 & 1.85 & 1.628e-02 & 1.01 \\ 
       	32768 & 4.122e-05 & 1.99 & 7.464e-04 & 1.88 & 8.129e-03 & 1.00 \\
       	\hline                                   
       	&\multicolumn{6}{c|}{$k=2$, BCs by Nitsche method}  \\ 
       	\hline                                  
       	64    & 1.801e-04 &  -  & 2.612e-03 & -     & 5.165e-03 &  -   \\ 
       	512   & 2.246e-05 & 3.00 & 2.882e-04 & 3.18 & 1.279e-03 & 2.01 \\ 
       	4096  & 2.808e-06 & 2.99 & 3.338e-05 & 3.11 & 3.187e-04 & 2.01 \\ 
       	32768 & 3.514e-07 & 2.99 & 4.405e-06 & 2.92 & 7.958e-05 & 2.00 \\ 
       	\hline                                 
       	&\multicolumn{6}{c|}{$k=3$, BCs by Nitsche method} \\ 
       	\hline                                
       	64    & 4.257e-06 &  -   & 5.952e-05 & -    & 1.817e-04 &  -   \\ 
       	512   & 2.782e-07 & 3.94 & 4.538e-06 & 3.73 & 2.252e-05 & 3.01 \\ 
       	4096  & 1.780e-08 & 3.97 & 3.392e-07 & 3.74 & 2.794e-06 & 3.01 \\ 
       	32768 & 1.126e-09 & 3.98 & 2.709e-08 & 3.65 & 3.478e-07 & 3.00 \\ 
       	\hline
       \end{tabular}
	\caption{Errors and convergence rates for BR2 dG discretizations of degree $k=\left\lbrace 1,2,3\right\rbrace $ over a $h$-refined mesh sequence of the unit cube, NHK-I constitutive model.	\label{tab:conv_result_neo_inc_3D}}
\end{table}

\begin{table}[H]
        \centering
       \small
       \begin{tabular}{|c|cccccc|}\hline
       	$\mathrm{card}(\Th)$ 
      & $\| \ub {-} \ub_h \|_{L^2(\Omega)}$ & \textbf{rate}
      & $\| p {-} p_h \|_{L^2(\Omega)}$ & \textbf{rate}
      & $\| \nabla(\ub {-} \ub_h) \|_{L^2(\Omega)}$ & \textbf{rate} \\
       	\hline                                                                                            
       	&\multicolumn{6}{c|}{$k=1$, BCs by Nitsche method}  \\ 
       	\hline           
       	64    & 2.596e-03 &  -   & 3.691e-02 & -    & 6.555e-02 &  -   \\
       	512   & 6.522e-04 & 1.99 & 1.437e-02 & 1.36 & 3.265e-02 & 1.01 \\
       	4096  & 1.653e-04 & 1.98 & 4.348e-03 & 1.72 & 1.628e-02 & 1.00 \\
       	32768 & 4.181e-05 & 1.98 & 1.180e-03 & 1.88 & 8.128e-03 & 1.00 \\
       	\hline                                   
       	&\multicolumn{6}{c|}{$k=2$, BCs by Nitsche method}  \\ 
       	\hline                                  
       	64    & 1.707e-04 &  -   & 4.901e-03 & -    & 5.175e-03 &  -   \\
       	512   & 2.127e-05 & 3.01 & 6.606e-04 & 2.89 & 1.282e-03 & 2.01 \\
       	4096  & 2.660e-06 & 2.99 & 1.339e-04 & 2.30 & 3.195e-04 & 2.00 \\
       	32768 & 3.328e-07 & 2.99 & 3.227e-05 & 2.05 & 7.977e-05 & 2.00 \\
       	\hline                                 
       	&\multicolumn{6}{c|}{$k=3$, BCs by Nitsche method} \\ 
       	\hline                                
       	64    & 4.248e-06 &  -   & 1.020e-04 & -    & 1.814e-04 &  -   \\
       	512   & 2.779e-07 & 3.93 & 6.755e-06 & 3.92 & 2.248e-05 & 3.01 \\
       	4096  & 1.779e-08 & 3.97 & 4.701e-07 & 3.85 & 2.792e-06 & 3.01 \\
       	32768 & 1.126e-09 & 3.98 & 3.463e-08 & 3.76 & 3.477e-07 & 3.00 \\
       	\hline
       \end{tabular}
	\caption{Errors and convergence rates for BR2 dG discretizations of degree $k=\left\lbrace 1,2,3\right\rbrace $ over a $h$-refined mesh sequence of the unit cube, SVK-I constitutive model. \label{tab:conv_result_svk_inc_3D}}
\end{table}

\subsection{2D simulations} \label{sec:2D_simulations}
\subsubsection{Parabolic indentation problem} \label{sec:indentation_problem}
The 2D \textit{parabolic indentation} problem imposes a severe deformation of parabolic shape to the top surface of a unit-length square computational domain $\Omega : [0,1]^2 $.
As proposed by Eyck and coworkers \cite{Eyck2008a}, the parabolic profile reads $v(\X) = 3 (X - 0.5)^2$ and the bottom surface is clamped. 
The computational mesh consist of 512 triangular elements and we consider a first degree BR2 dG discretization.
Dirichlet boundary conditions are imposed on the top and bottom surfaces while homogeneous Neumann boundary conditions are enforced on the rest of the boundary.
We set $\beta = 0$ and $\epsilon = 0$, meaning that the adaptive stabilization strategy is switched-off, and, in the incompressible regime, we also set $\eta_{\text{LBB}} = 1$.

\begin{figure}[H]
	\vspace{-1.5cm}
	\centering
	\subfloat[NHK-C]
	{
		\hspace{-0.5cm}
		\includegraphics[width=0.24\textwidth]{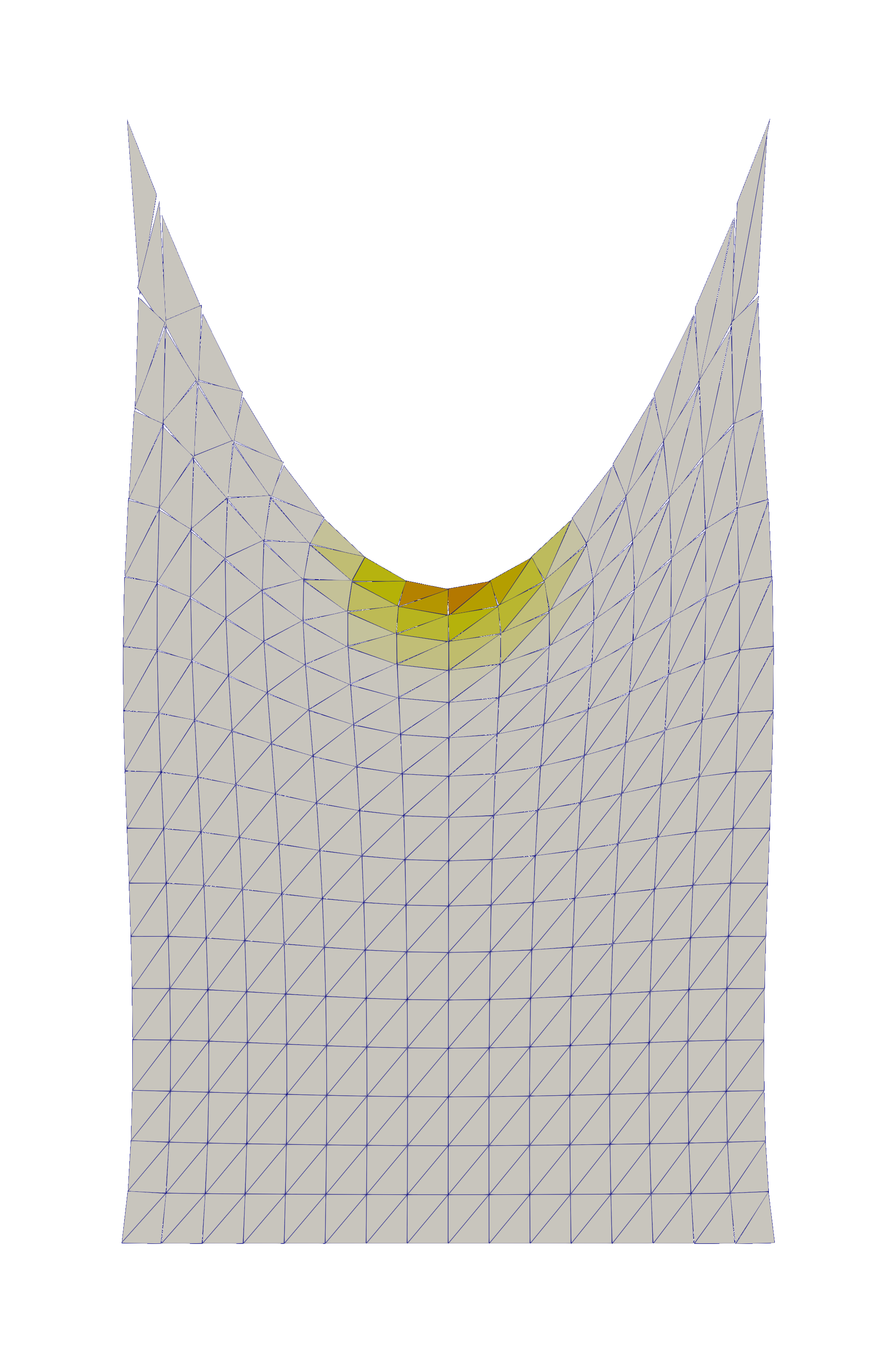}}
	\subfloat[NHK-I]
	{
		\hspace{-0.5cm}
		\includegraphics[width=0.24\textwidth]{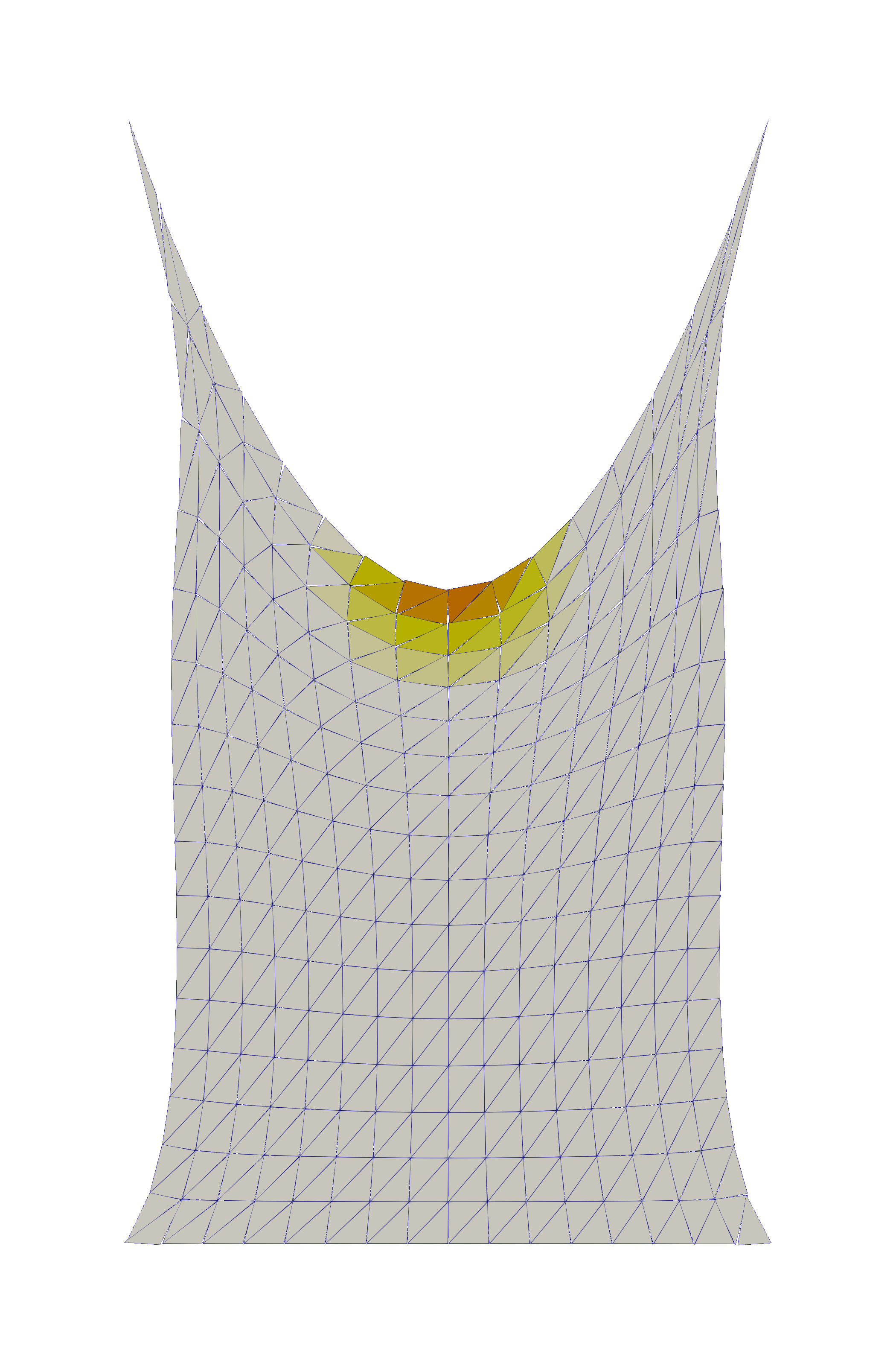}}
	\subfloat[SVK-C]
	{
		\hspace{-0.5cm}
		\includegraphics[width=0.24\textwidth]{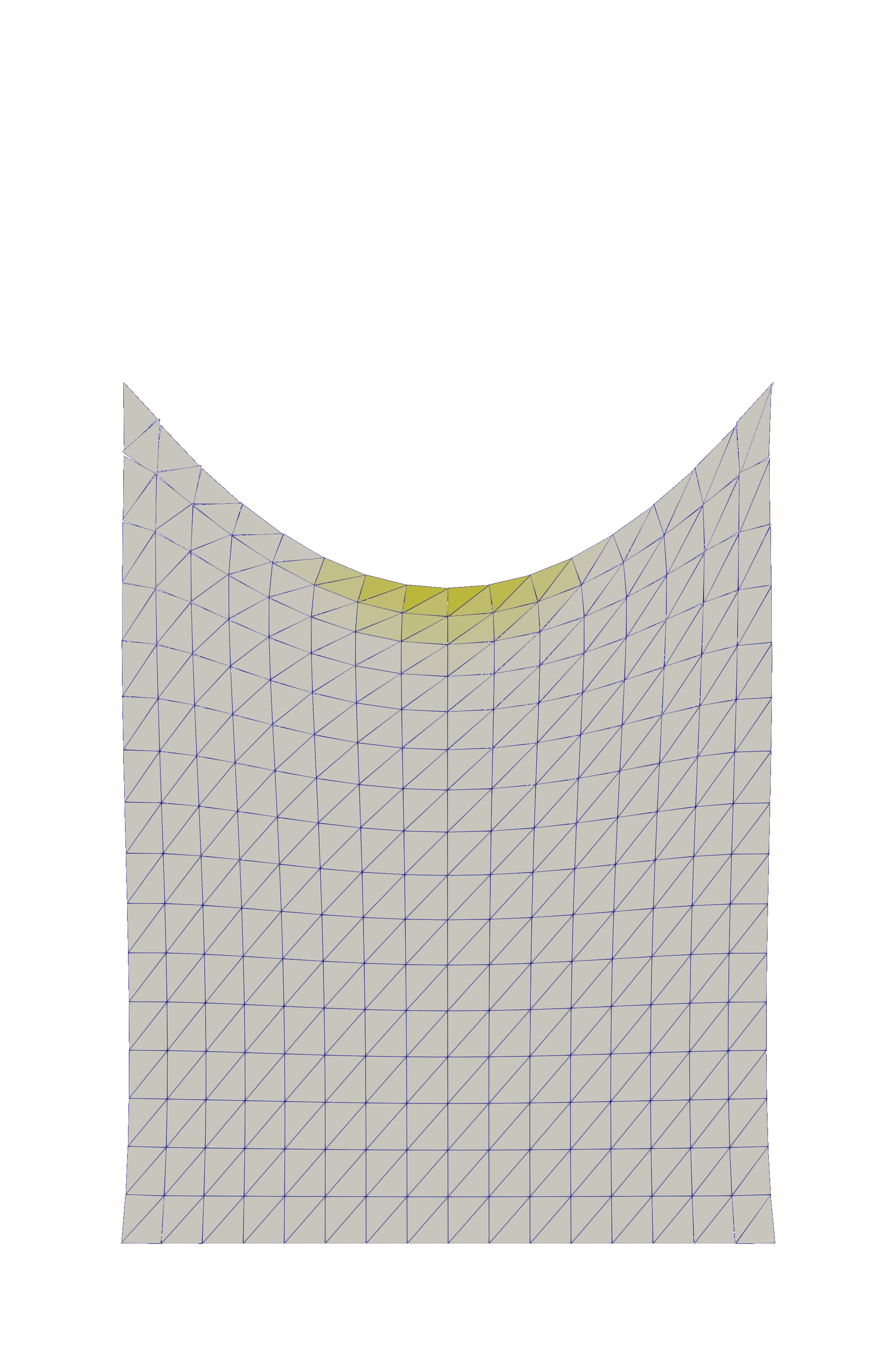}}
	\subfloat[SVK-I]
	{
		\hspace{-0.5cm}
		\includegraphics[width=0.24\textwidth]{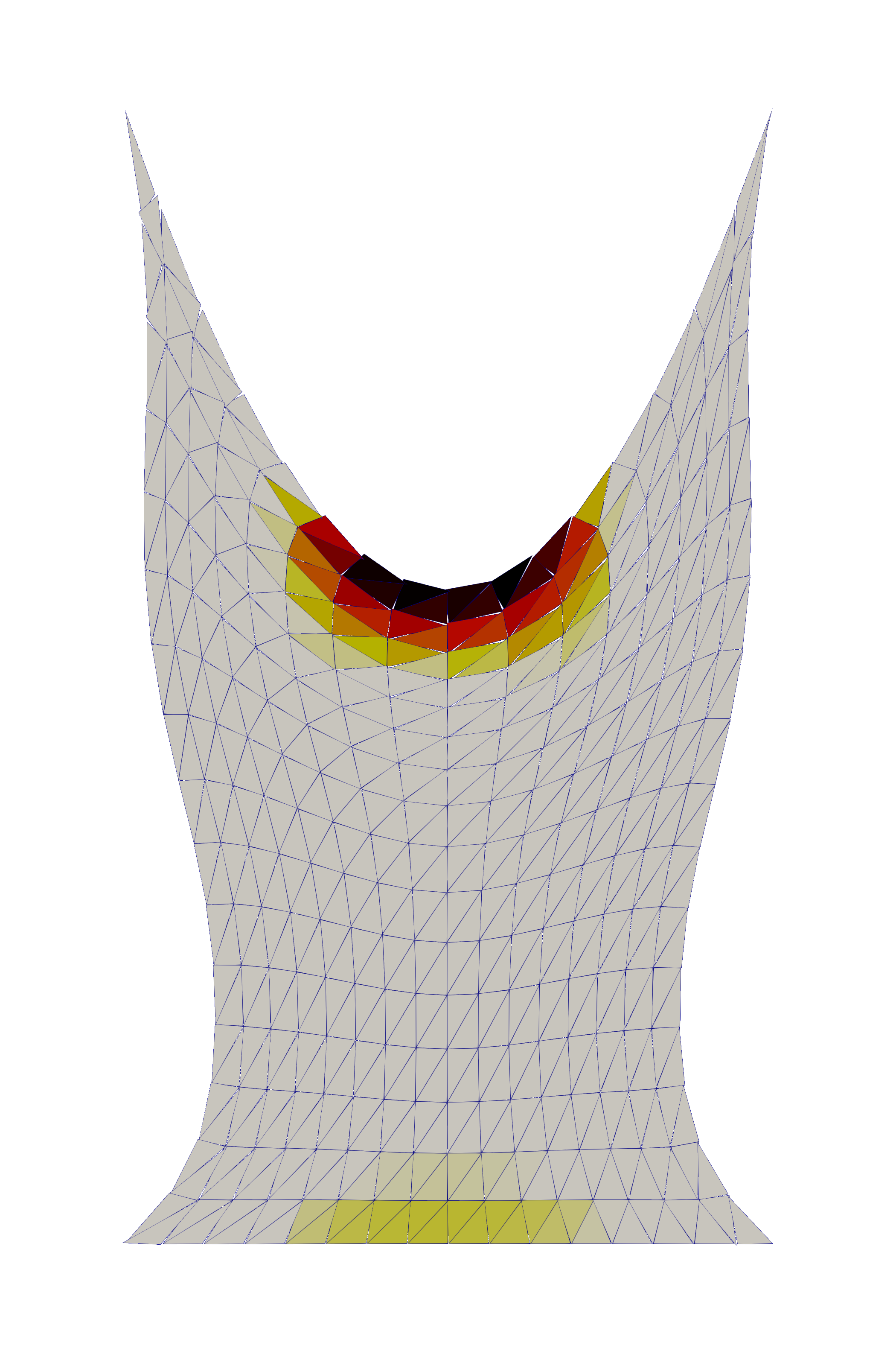}}
	\subfloat
	{
		\includegraphics[width=0.1\textwidth]{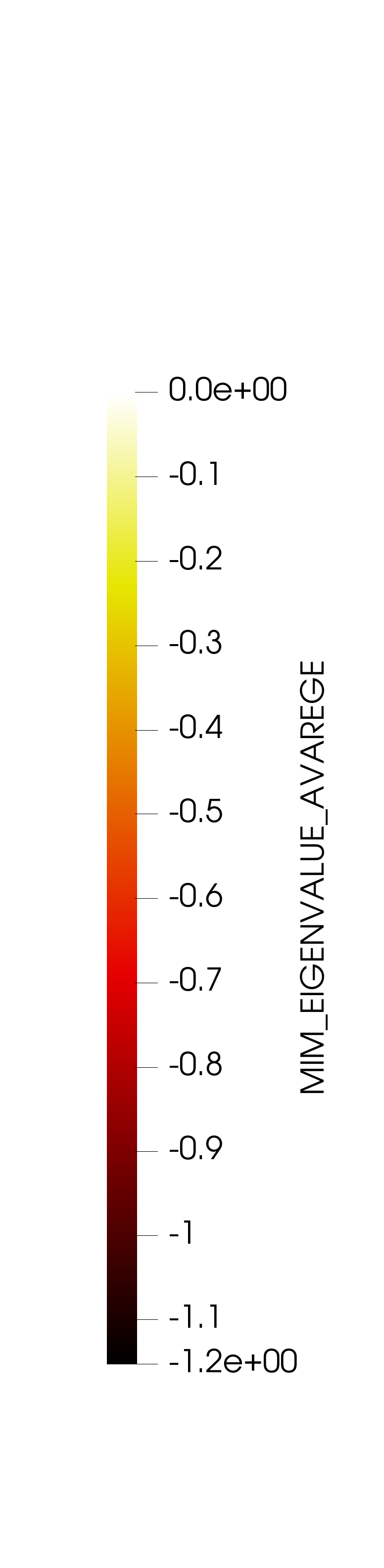}}
	\caption{Deformed configurations of the parabolic indentation problem using NHK and the SVK constitutive models. 
                 Images are colour coded based on the minimum negative eigenvalue of the fourth order elasticity tensor $\Tnsr{A}$.
	\label{fig:parabolic_indentation_eigenvalues}}.
\end{figure}

The deformed states obtained with all constitutive models relevant for this configuration are depicted in \Fig~\ref{fig:parabolic_indentation_eigenvalues},
material parameters reads $\mu = \lambda = 0.4$.
When using the SVK-C model, Newton's method fails to converge when reaching 50\% of the loading path, irrespectively of the amount of stabilization introduced. 
Accordingly, the final configuration is not attained. 
This behavior can be explained by noticing that, as opposite to the NHK-C model, 
the SVK strain-energy function \eqref{eq:svk_energy} lacks of any term preventing the onset of negative Jacobian values, see also \cite{Ciarlet1999a}. 
We remark that the SVK-I model is successful, because $\det(\F){=}1$ is weakly enforced in accordance with the incompressibility constraint. 

\subsubsection{2D beam deformation} \label{sec:beam_deformation}

As proposed by Eyck and co-workers \cite{Eyck2008}, we challenge the adaptive stabilisation strategy considering the deformation of a 2D beam: 
the bottom surface is clamped while the upper surface of the beam is first rotated by $\pi/2$ and, then, translated in the direction orthogonal to the beam axis. 
We consider a NHK-C constitutive law and we set $\nu =0.3$ and $E = 1$.
The computational mesh consists of 110 triangular elements and we employ a $k=1$ BR2 dG discretization.  
Dirichlet boundary conditions are imposed on the top and bottom surfaces while homogeneous Neumann boundary conditions are enforced on the rest of the boundary.

\begin{figure}[H]
	\subfloat{
		\hspace{-4.cm}
		\includegraphics[width=0.43\textwidth]{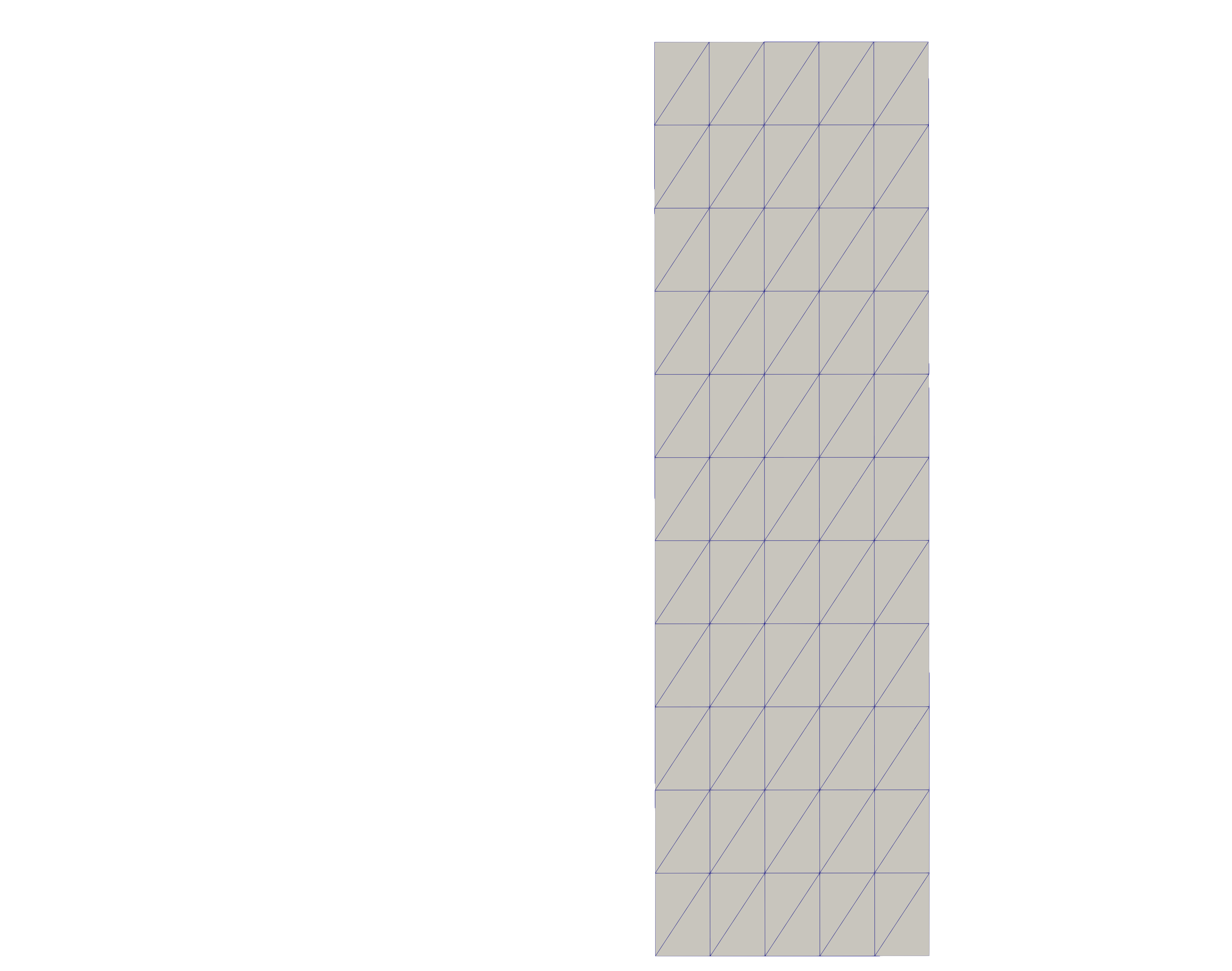}}
		\subfloat{
		\hspace{-5.3cm}
		\includegraphics[width=0.43\textwidth]{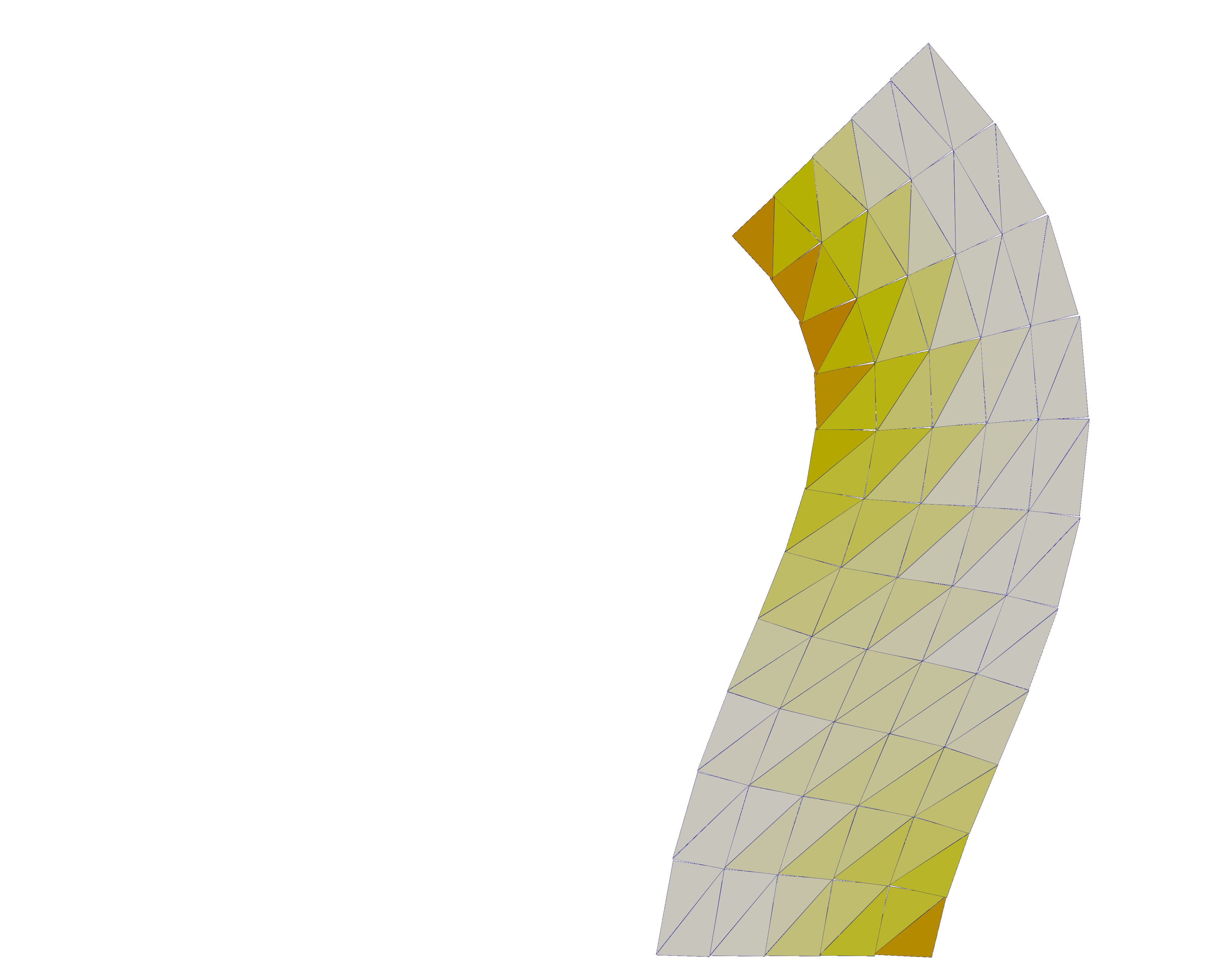}}
		\subfloat{
		\hspace{-5.3cm}
		\includegraphics[width=0.43\textwidth]{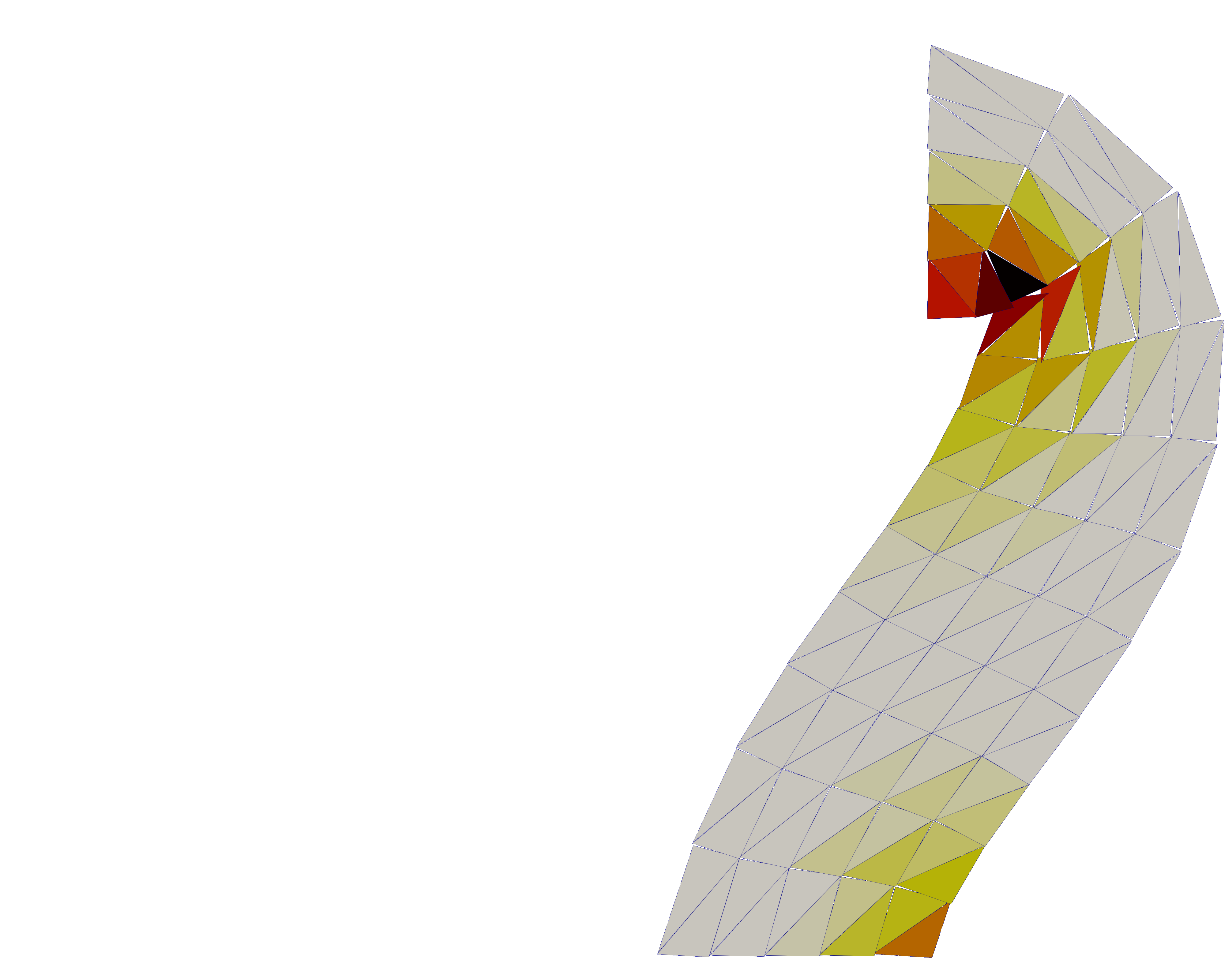}}
		\subfloat{
		\hspace{-2.8cm}
		\includegraphics[width=0.43\textwidth]{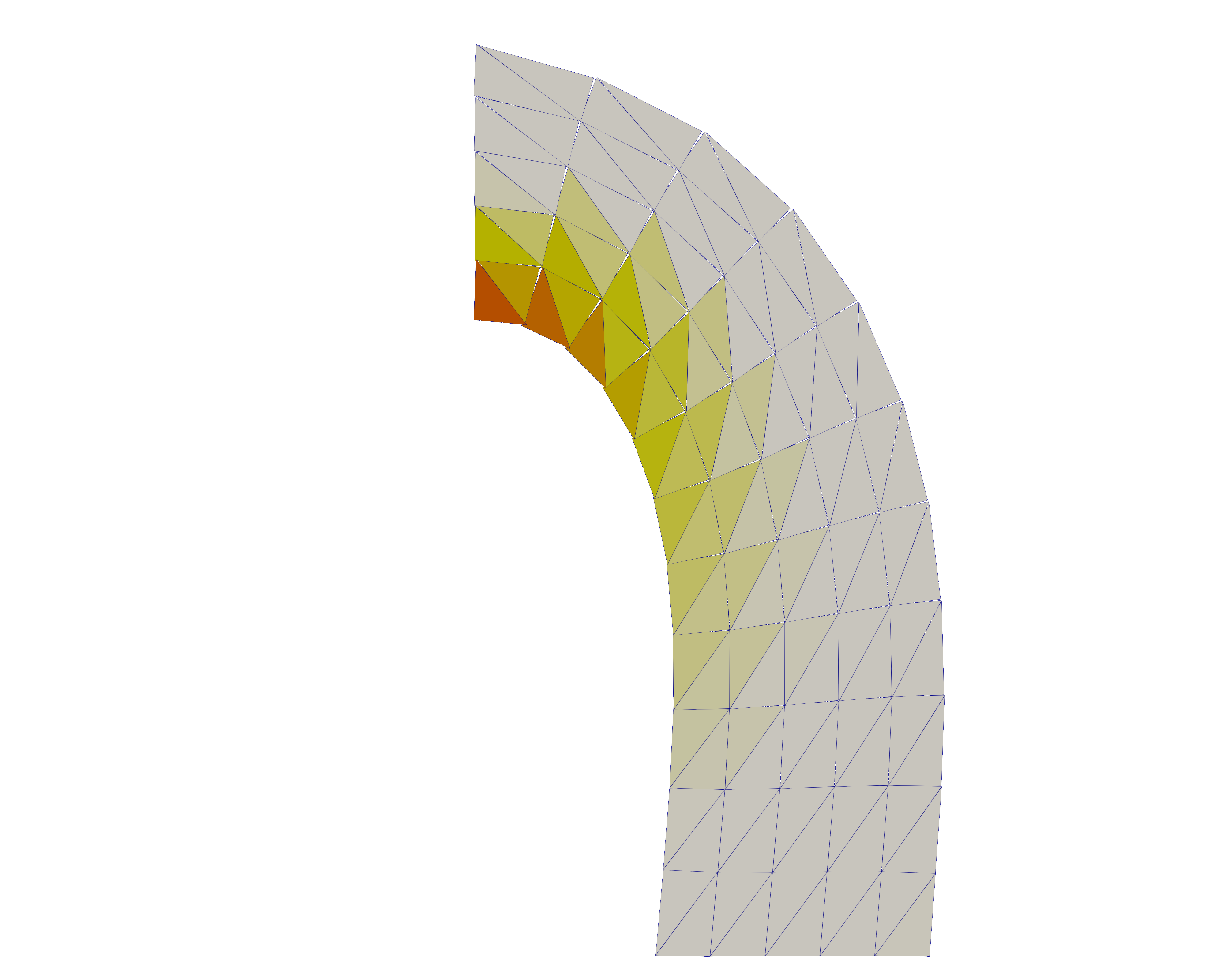}}
			\subfloat{
		\hspace{-2.0cm}
		\includegraphics[width=0.43\textwidth]{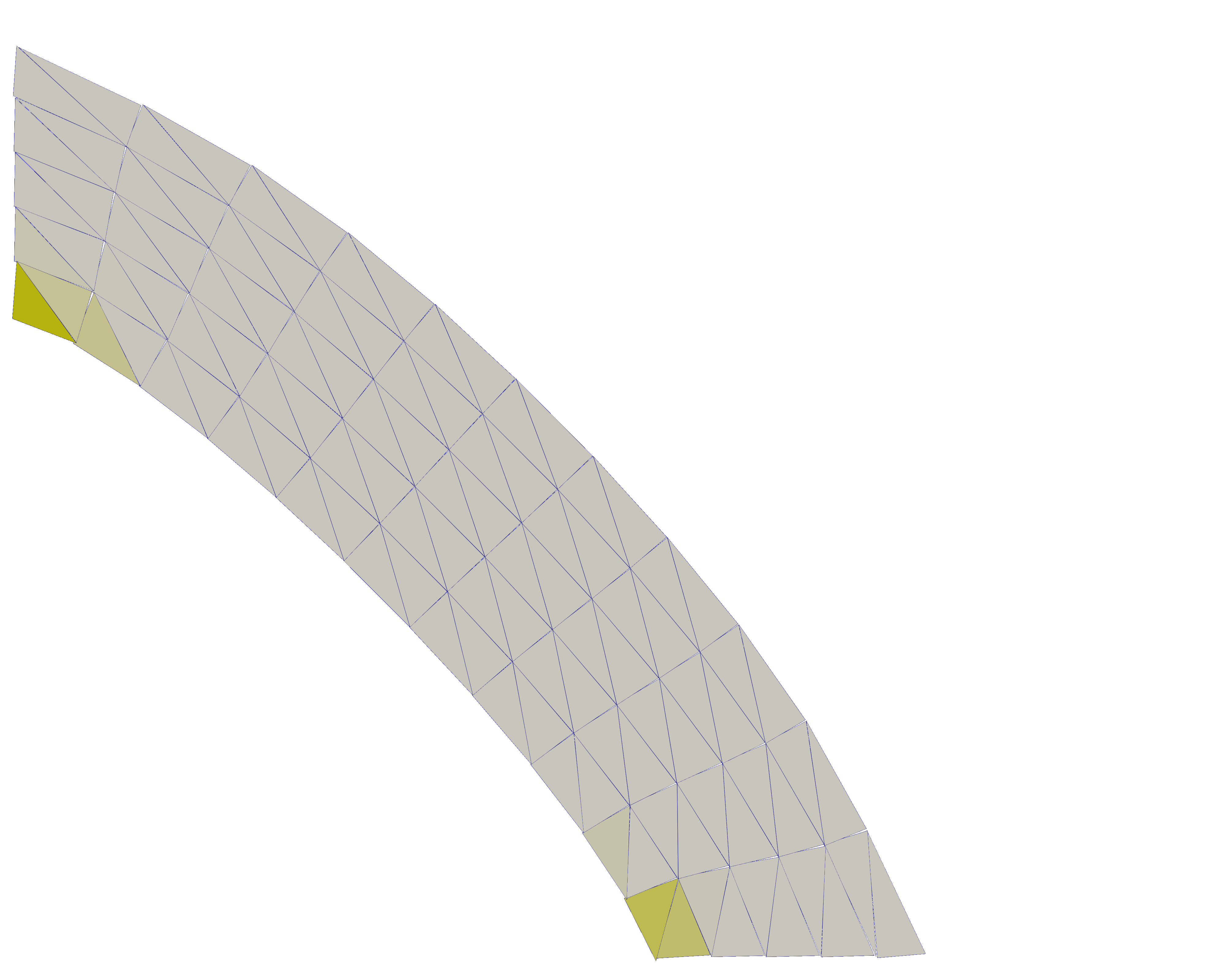}}
		\subfloat{
		\hspace{-2.3cm}
		\includegraphics[width=0.08\textwidth]{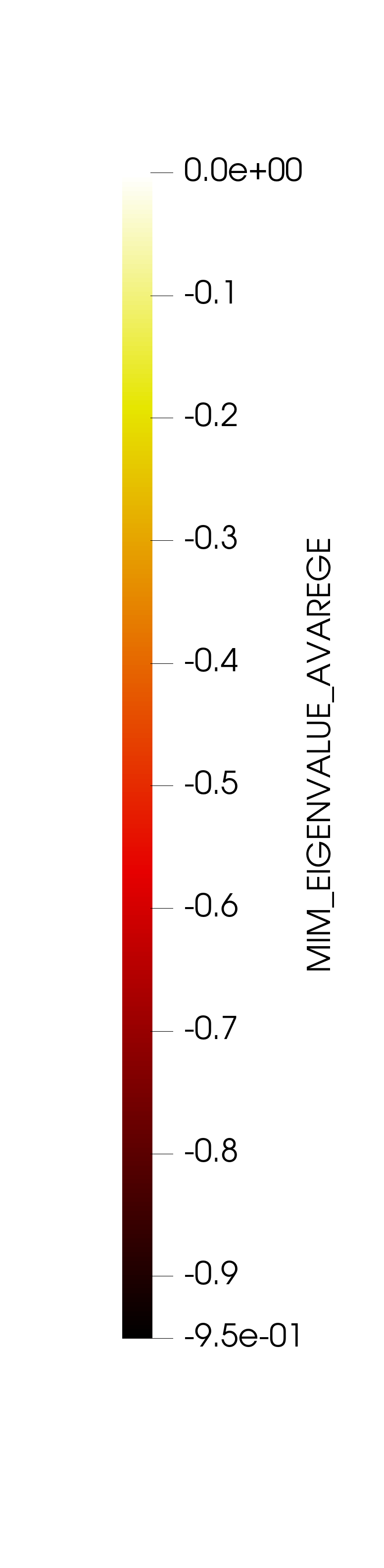}}
	\caption{Deformation of a 2D NHK-C beam. Images are colour-coded based on the minimum negative eigenvalue of the fourth-order elasticity tensor $\Tnsr{A}$.
	\label{fig:beam_2D_def_eigenvalues}}
\end{figure}

\Fig~\ref{fig:beam_2D_def_eigenvalues} depicts the beam deformation by showing a sequence of deformed states consistent with the loading path.
Deformed states are colour-coded with minimum negative eigenvalues of the elasticity tensor 
allowing to appreciate that compression of the beam material triggers the adaptive stabilization strategy.  

\begin{figure}[H]
	\subfloat[$\epsilon = 0$]{ \includegraphics[width=0.5\linewidth]{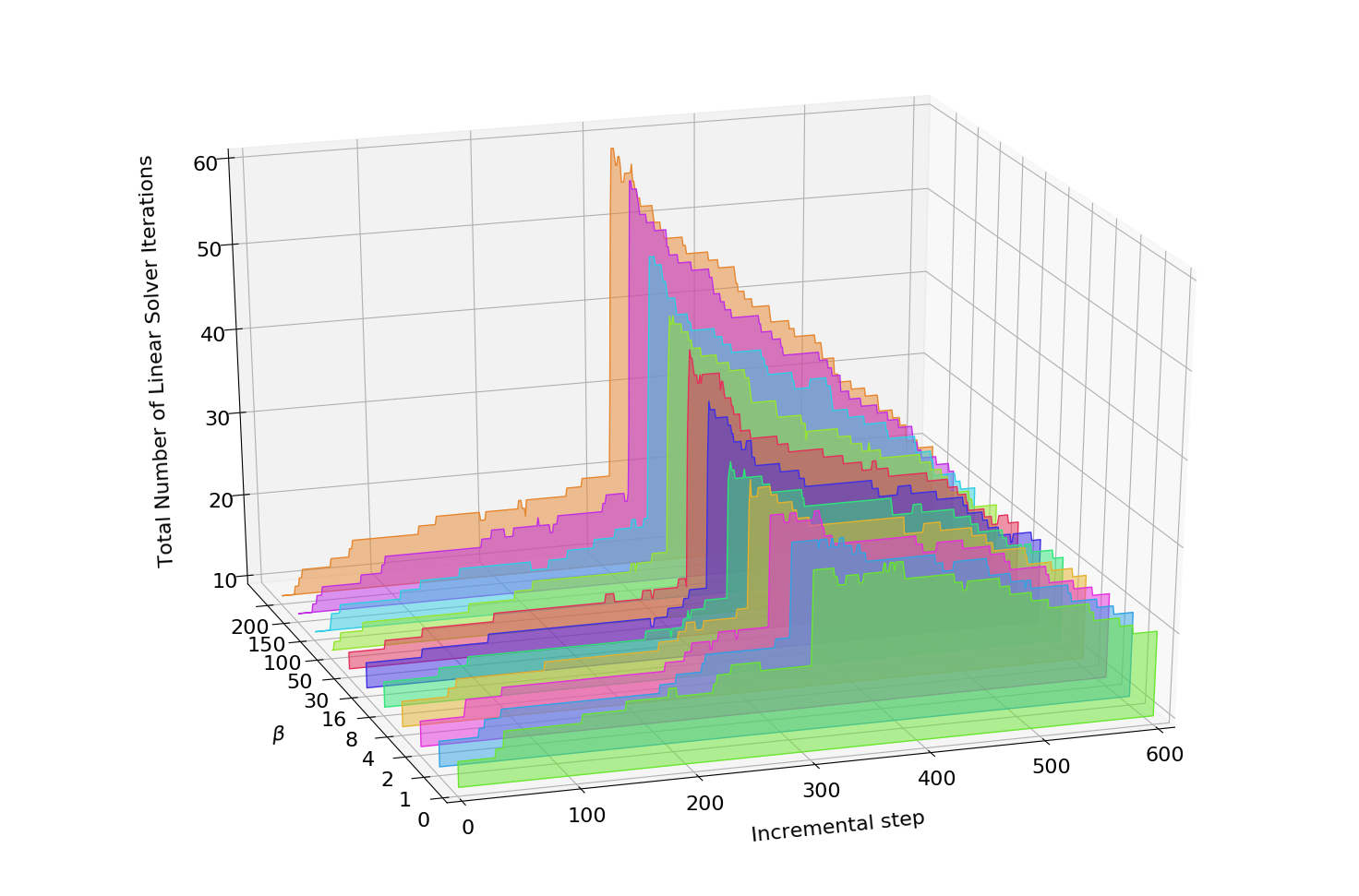} }
	\subfloat[$\epsilon = 1$]{ \includegraphics[width=0.5\linewidth]{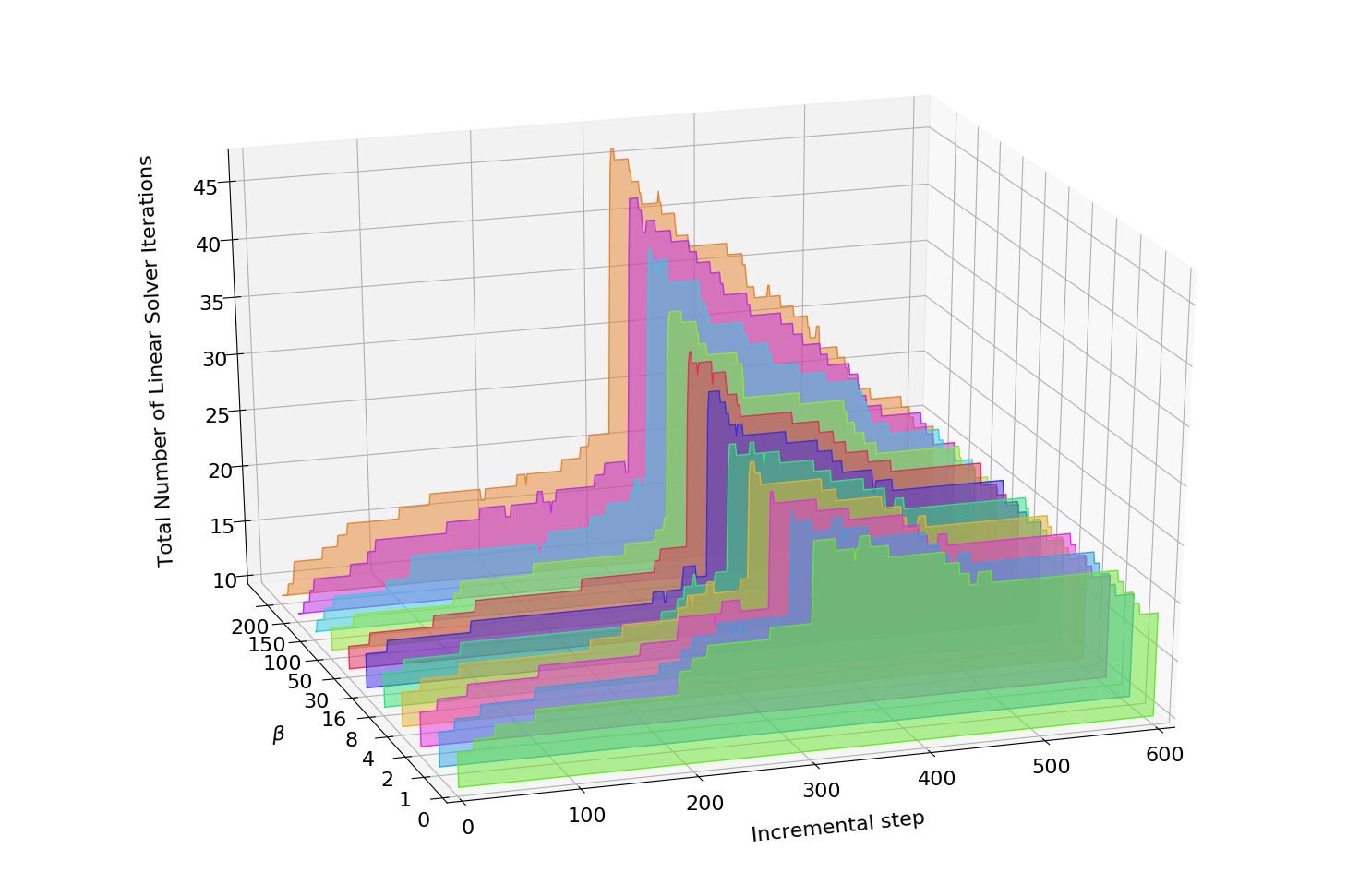} }
	
	\hspace{-0.8cm}
	\subfloat[$\epsilon = 10$]{ \includegraphics[width=0.5\linewidth]{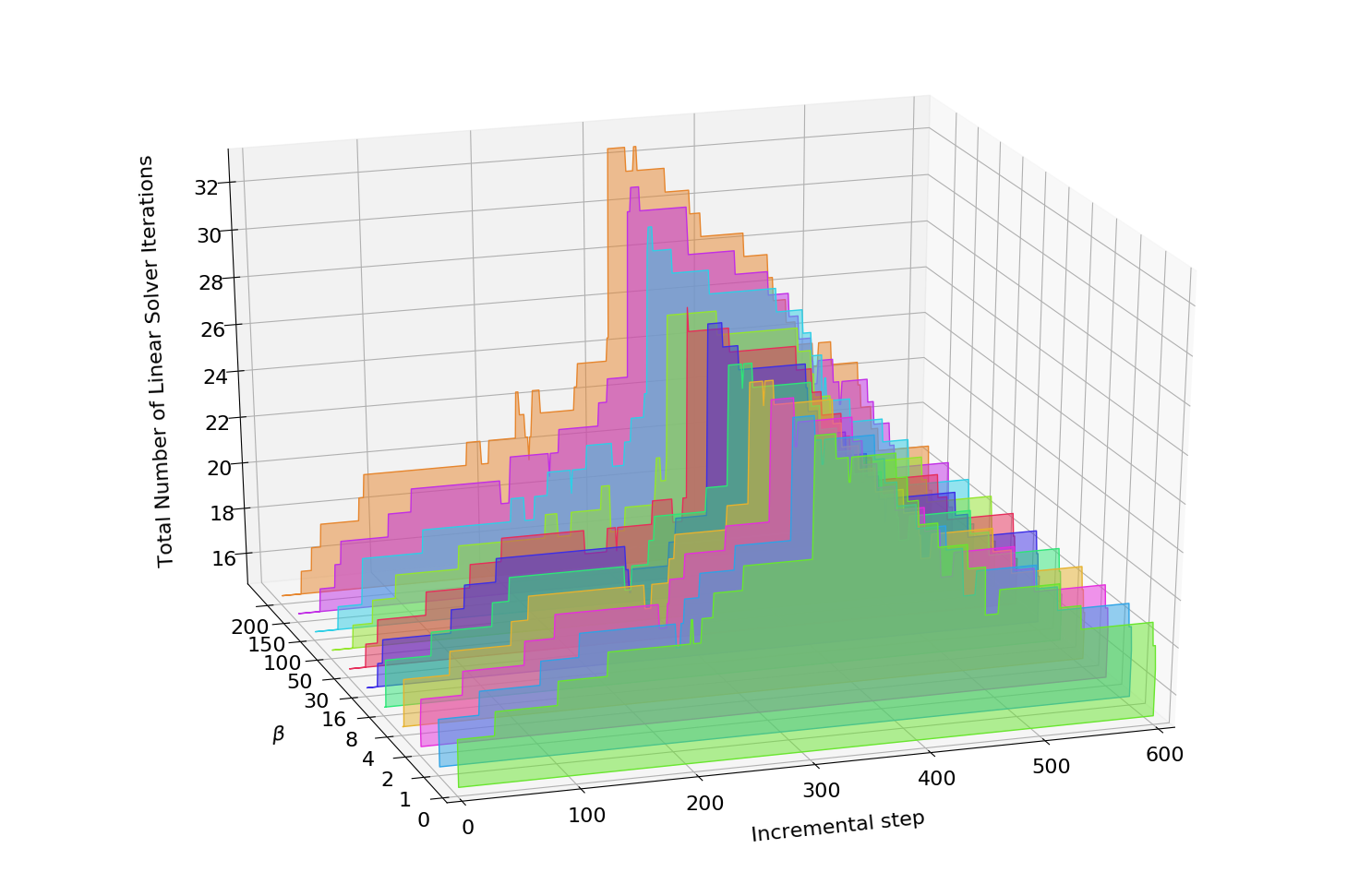} }
	\subfloat[$\epsilon = 20$]{ \includegraphics[width=0.5\linewidth]{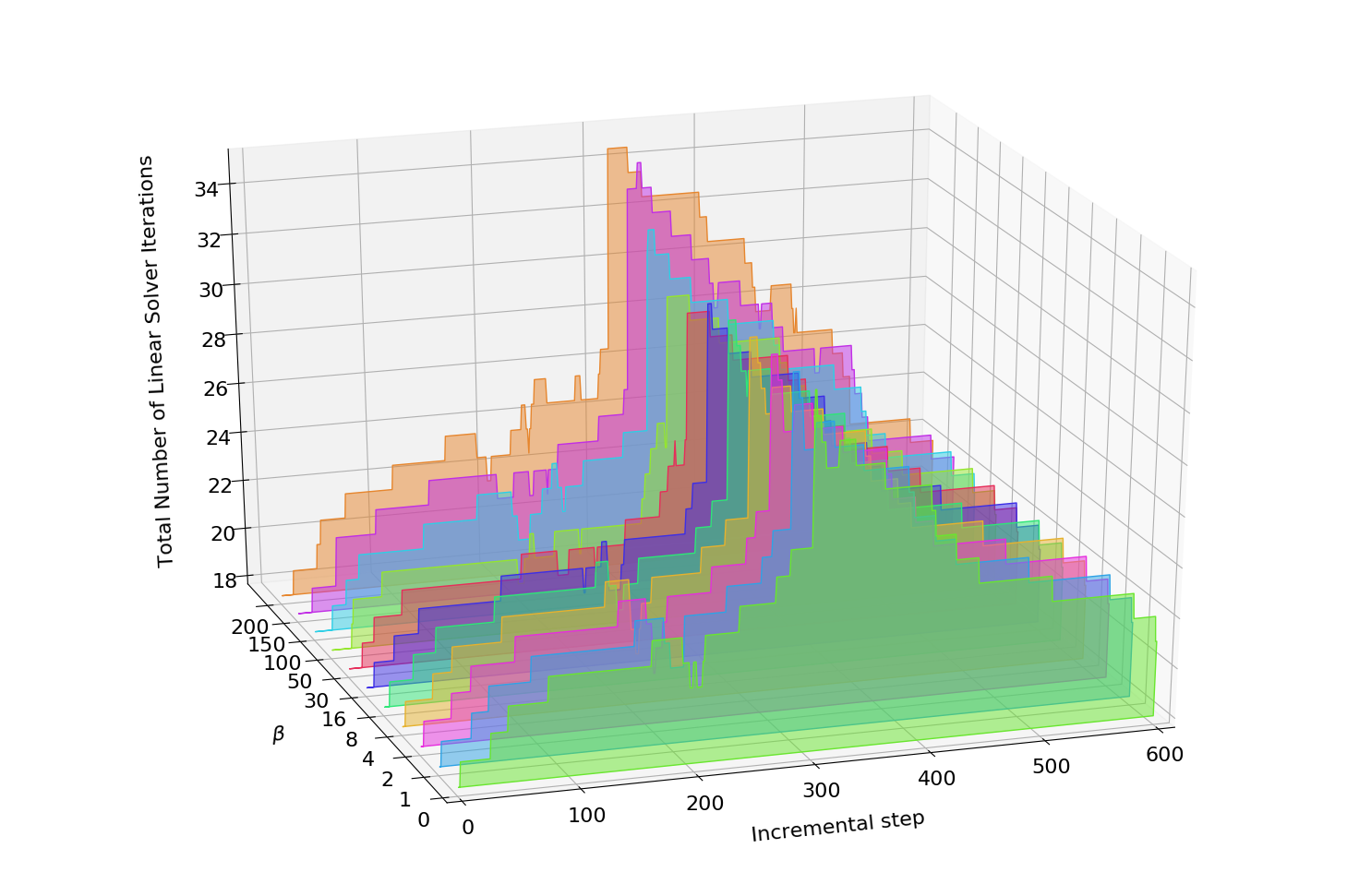} }
    \caption{2D NHK-C beam: total number of linear solver iterations recorded along the loading path (600 increments). 
             Results are obtained varying the stabilization parameters $\beta$ and $\epsilon$ in order to show their influence on the performance of the solution strategy.
	\label{fig:beam_rot_trasl_beta_test_total_iterations}}
\end{figure}

In order to study the influence of the stabilization parameter on the performance of the $h$-multigrid solution strategy, 
we run a series of test varying $\beta$ in the interval $[0,200]$ for each $\epsilon$ in $\left\lbrace 0, 1, 10, 20 \right\rbrace $.
For each combination of $\beta$ and $\epsilon$, the history of total linear solver iterations 
recorded along the loading path, counting of 600 incremental steps, is depicted in \Fig~\ref{fig:beam_rot_trasl_beta_test_total_iterations}. 
Furthermore, the average and maximum number of Newton iteration as well as the average and maximum number of linear solver iterations
are tabulated in \Tab~\ref{tab:2D_beam_beta_tests}.

\begin{table}[H]
	\centering
	\begin{tabular}{|c|l|c|c|c|c|c|c|c|c|c|c|c|c|}
		\hline 
		
		\multirow{2}{*}{\begin{tabular}[l]{@{}l@{}}$\epsilon$ \end{tabular}}  & \multicolumn{2}{|l|}{\multirow{2}{*}{}} & \multicolumn{11}{c|}{$\beta$}\\ \cline{4-14}
		&  \multicolumn{2}{|l|}{}  &  \multicolumn{1}{l|}{\textbf{0}}  &  \multicolumn{1}{l|}{\textbf{1}}  &  \multicolumn{1}{l|}{\textbf{2}}  &  \multicolumn{1}{l|}{\textbf{4}}  &  \multicolumn{1}{l|}{\textbf{8}}  &  \multicolumn{1}{l|}{\textbf{16}}  &  \multicolumn{1}{l|}{\textbf{30}}  &  \multicolumn{1}{l|}{\textbf{50}}  &  \multicolumn{1}{l|}{\textbf{100}}  &  \multicolumn{1}{l|}{\textbf{150}}  &  \multicolumn{1}{l|}{\textbf{200}}\\ 
		\hline
		\hline
		
		\multirow{4}{*}{\begin{tabular}[l]{@{}l@{}} 0 \end{tabular}} & \multirow{2}{*}{\begin{tabular}[l]{@{}l@{}}\textbf{Newton} \\  \textbf{iterations}\end{tabular}}   & \textbf{mean}   &  6  &  6  &  6  &  6  &  6  &  6  &  6  &  6  &  6  &  6  &  6\\ \cline{3-14}
		&  & \textbf{max}   &  7  &  7  &  7  &  8  &  8  &  8  &  9  &  9  &  10  &  10  &  10\\ \cline{2-14}
		&  \multirow{2}{*}{\begin{tabular}[l]{@{}l@{}}\textbf{Linear Solver} \\ \textbf{iterations} \end{tabular}}   & \textbf{mean}  &  4  &  4  &  4  &  4  &  4  &  4  &  4  &  4  &  5  &  5  &  5\\ \cline{3-14}
		&  & \textbf{max}    &  \multicolumn{1}{c|}{6}    &  \multicolumn{1}{c|}{6}    &  \multicolumn{1}{c|}{7}    &  \multicolumn{1}{c|}{7}    &  \multicolumn{1}{c|}{7}    &  \multicolumn{1}{c|}{7}    &  \multicolumn{1}{c|}{7}    &  \multicolumn{1}{c|}{7}    &  \multicolumn{1}{c|}{8}    &  \multicolumn{1}{c|}{8}    &  \multicolumn{1}{c|}{9}  \\ 
		\hline
		\hline
		\multirow{4}{*}{\begin{tabular}[l]{@{}l@{}} 1 \end{tabular}} & \multirow{2}{*}{\begin{tabular}[l]{@{}l@{}}\textbf{Newton} \\  \textbf{iterations}\end{tabular}}   & \textbf{mean}   &  5  &  5  &  6  &  6  &  6  &  5  &  5  &  5  &  5  &  5  &  5\\ \cline{3-14}
		&  & \textbf{max}   &  6  &  7  &  7  &  7  &  7  &  7  &  7  &  7  &  8  &  8  &  8\\ \cline{2-14}
		&  \multirow{2}{*}{\begin{tabular}[l]{@{}l@{}}\textbf{Linear Solver} \\ \textbf{iterations} \end{tabular}}   & \textbf{mean}  &  4  &  4  &  4  &  4  &  4  &  4  &  4  &  4  &  5  &  5  &  5\\ \cline{3-14}
		&  & \textbf{max}    &  \multicolumn{1}{c|}{6}    &  \multicolumn{1}{c|}{7}    &  \multicolumn{1}{c|}{7}    &  \multicolumn{1}{c|}{7}    &  \multicolumn{1}{c|}{7}    &  \multicolumn{1}{c|}{7}    &  \multicolumn{1}{c|}{7}    &  \multicolumn{1}{c|}{7}    &  \multicolumn{1}{c|}{8}    &  \multicolumn{1}{c|}{8}    &  \multicolumn{1}{c|}{8}  \\ 
		\hline
		\hline
		\multirow{4}{*}{\begin{tabular}[l]{@{}l@{}} 10 \end{tabular}} & \multirow{2}{*}{\begin{tabular}[l]{@{}l@{}}\textbf{Newton} \\  \textbf{iterations}\end{tabular}}   & \textbf{mean}   &  5  &  5  &  5  &  5  &  5  &  5  &  5  &  5  &  5  &  5  &  5\\ \cline{3-14}
		&  & \textbf{max}   &  5  &  5  &  5  &  5  &  5  &  5  &  5  &  6  &  6  &  6  &  5\\ \cline{2-14}
		&  \multirow{2}{*}{\begin{tabular}[l]{@{}l@{}}\textbf{Linear Solver} \\ \textbf{iterations} \end{tabular}}   & \textbf{mean}  &  5  &  5  &  5  &  5  &  5  &  5  &  5  &  5  &  6  &  6  &  6\\ \cline{3-14}
		&  & \textbf{max}    &  \multicolumn{1}{c|}{8}    &  \multicolumn{1}{c|}{8}    &  \multicolumn{1}{c|}{8}    &  \multicolumn{1}{c|}{8}    &  \multicolumn{1}{c|}{8}    &  \multicolumn{1}{c|}{8}    &  \multicolumn{1}{c|}{8}    &  \multicolumn{1}{c|}{8}    &  \multicolumn{1}{c|}{8}    &  \multicolumn{1}{c|}{9}    &  \multicolumn{1}{c|}{9}  \\ 
		\hline
		\hline
		\multirow{4}{*}{\begin{tabular}[l]{@{}l@{}} 20 \end{tabular}} & \multirow{2}{*}{\begin{tabular}[l]{@{}l@{}}\textbf{Newton} \\  \textbf{iterations}\end{tabular}}   & \textbf{mean}   &  5  &  5  &  5  &  5  &  5  &  5  &  4  &  4  &  4  &  4  &  4\\ \cline{3-14}
		&  & \textbf{max}   &  5  &  5  &  5  &  5  &  5  &  5  &  5  &  5  &  5  &  5  &  5\\ \cline{2-14}
		&  \multirow{2}{*}{\begin{tabular}[l]{@{}l@{}}\textbf{Linear Solver} \\ \textbf{iterations} \end{tabular}}   & \textbf{mean}  &  6  &  6  &  6  &  6  &  6  &  6  &  6  &  6  &  7  &  7  &  7\\ \cline{3-14}
		&  & \textbf{max}    &  \multicolumn{1}{c|}{9}    &  \multicolumn{1}{c|}{9}    &  \multicolumn{1}{c|}{9}    &  \multicolumn{1}{c|}{9}    &  \multicolumn{1}{c|}{9}    &  \multicolumn{1}{c|}{9}    &  \multicolumn{1}{c|}{9}    &  \multicolumn{1}{c|}{9}    &  \multicolumn{1}{c|}{9}    &  \multicolumn{1}{c|}{10}    &  \multicolumn{1}{c|}{10}  \\ 
		\hline
	\end{tabular}
    \caption{2D NHK-C beam: average and maximum number of Newton and linear solver iterations recorded along the loading path. 
             Results are obtained varying the stabilization parameters $\beta$ and $\epsilon$ in order to show their influence on the performance of the solution strategy.
	\label{tab:2D_beam_beta_tests}}
\end{table}

We remark that, in the range $ 0 \leq \beta \leq 200$, the number of linear solver iterations is pretty stable, 
and setting $\epsilon = 1$ reduces the number of Newton iteration resulting in a decrease of the total number of the linear solver iterations per incremental step.

\subsubsection{Cavitating voids} \label{sec:cavitating_voids}
In solid mechanics, the term \textit{cavitation} refers to the formation and rapid expansion 
of voids that occurs when a solid is subjected to sufficiently large tensile stresses. 
Some experiments on the cavitation are reported by Gent and Lindley~\eg~\cite{Gent1959} 
where unusual internal flaws appear in vulcanized rubber cylinders under a well-defined relatively small tensile load.
Since, during the growth of voids, significant deformation occurs near the cavities, 
the numerical simulation of cavitation requires numerical methods that are robust respect to mesh distortion. 
A Crouzeix-Raviart nonconforming finite element method was presented in Xu and Henao~\cite{Xu2011} while, more recently, 
the cavitation problem has been studied using HDG \cite{Kabaria2015} and HHO \cite{Abbas2018} discretizations. 
The interested reader may refer to the review by Xu \ea~\cite{Xu2011}.

We consider a unit radius disk centered at the origin with two holes: 
the first centered at $C_1=(-0.3, 0, 0)$ with radius 0.25 and the second centered at $C_2=(0.3, 0, 0)$ with radius 0.2. 
The disc is expanded by imposing Dirichlet boundary conditions $\bm{g}_\dir = (\alpha - 1) \X$, with $\alpha \geq 1$,
on the outer surface ($|\X| = 1$ in reference configuration), while imposing traction-free Neumann boundary conditions on the inner walls of the holes.
For cavitation to occur, we rely on the strain energy function in \eqref{eq:neoh_energy_cavitation}.
Note that, with respect to the standard NHK-C  law \eqref{eq:neoh_energy}, NHK-CAV has been modified
to reduce the rate of the strain-energy growth with respect to the deformation gradient. 
In order to enable direct comparison with \cite{Kabaria2015}, we use the same material configuration setting $\mu = 0.1$ and $\lambda = 1$ and $\alpha = 4.7$. 

\begin{figure}[H]
	\hspace{-1cm}
	\centering
	\subfloat[$k = 1$]{
		\includegraphics[width=0.3\textwidth]{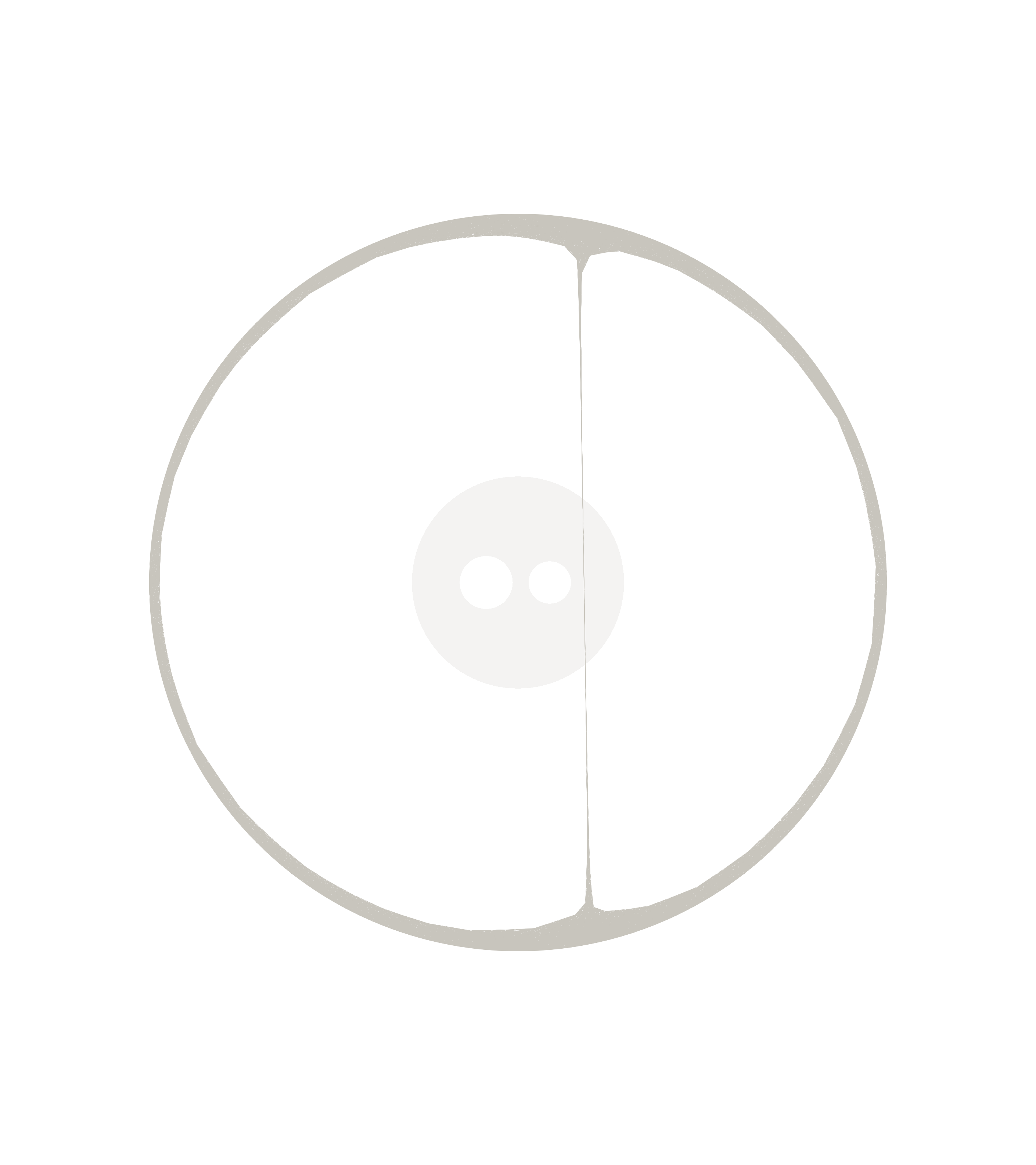} \label{fig:circle_voids_p_order_1}}
	\quad
	\subfloat[$k = 2$]{
		\includegraphics[width=0.3\textwidth]{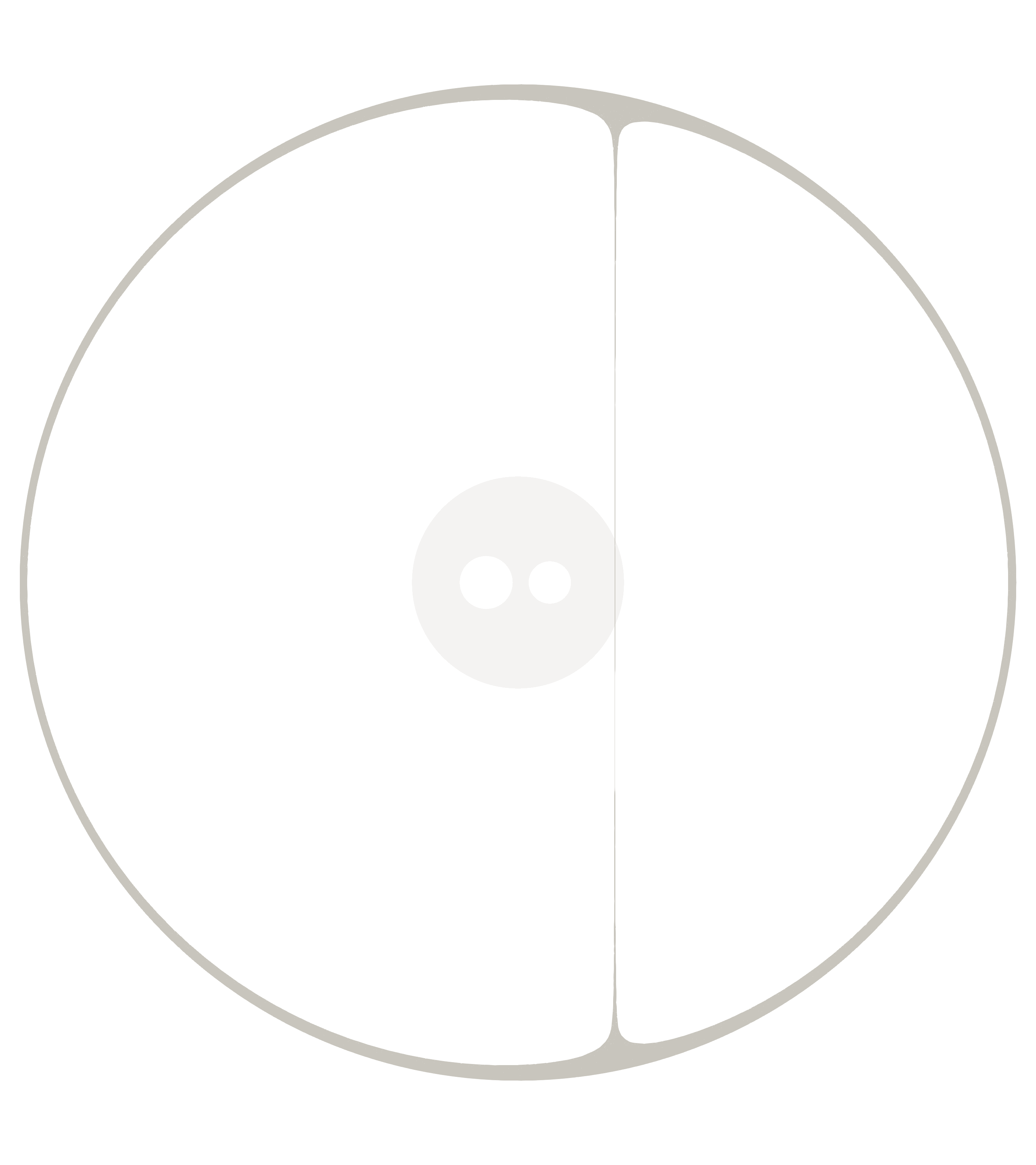}\label{fig:circle_voids_p_order_2}}
	\quad
	\subfloat[$k = 3$]{
		\includegraphics[width=0.3\textwidth]{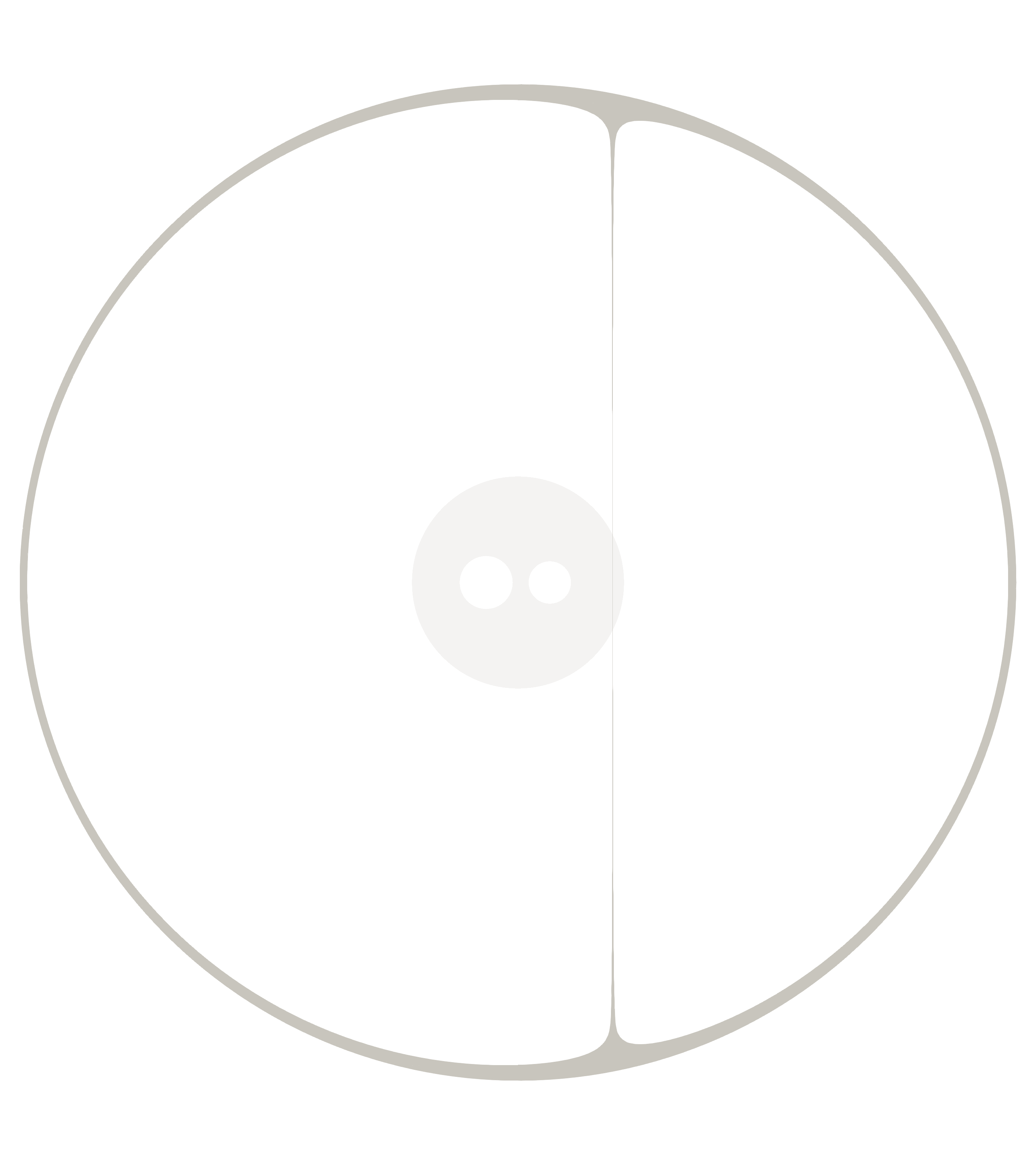}\label{fig:circle_voids_p_order_3}}
    \quad
	\caption{NHK-CAV disk with two holes subjected to tensile stresses. For each polynomial degree $k=\{ 1,2,3 \}$ the reference configuration (small disk with circular holes) 
         and deformed configurations (big disk with stretched holes) are shown.
	\label{fig:circle_voids_p_order}}
\end{figure}

\Fig~\ref{fig:circle_voids_p_order} reports the results obtained choosing $k = \{ 1,2,3 \}$. 
It is interesting to remark that only higher-order $k=2,3$ dG discretizations are able to reach the final configuration
while the first degree $k=1$ discretization fails at 64\% of the loading path due to the onset of negative Jacobian values, see \Fig~\ref{fig:circle_voids_p_order_1}. 

\begin{figure}[H]
	\hspace{-1cm}
	\centering
	\subfloat[$k = 1$]{
		\includegraphics[width=0.25\textwidth]{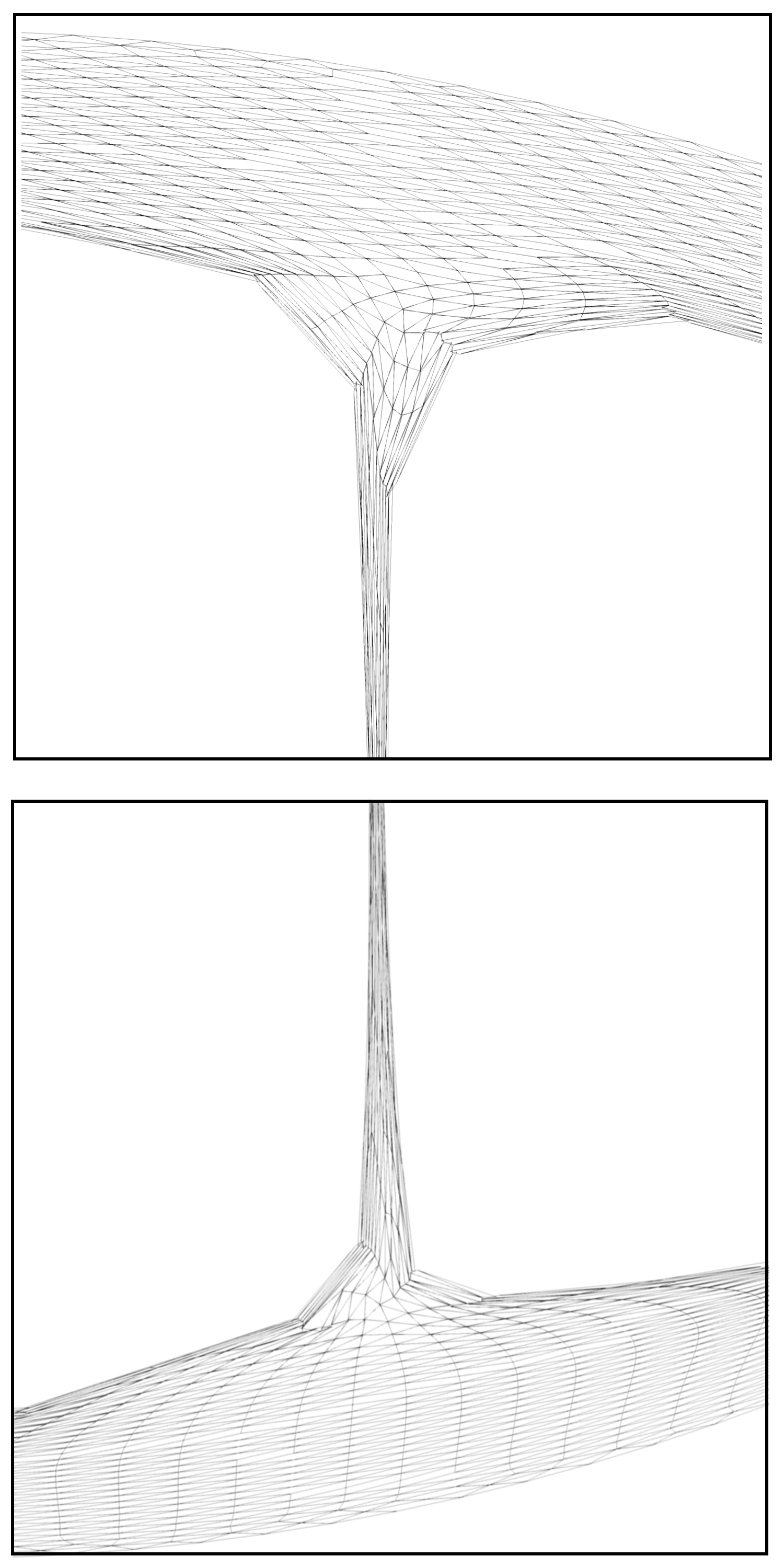}}
	\quad
	\subfloat[$k = 2$]{
		\includegraphics[width=0.25\textwidth]{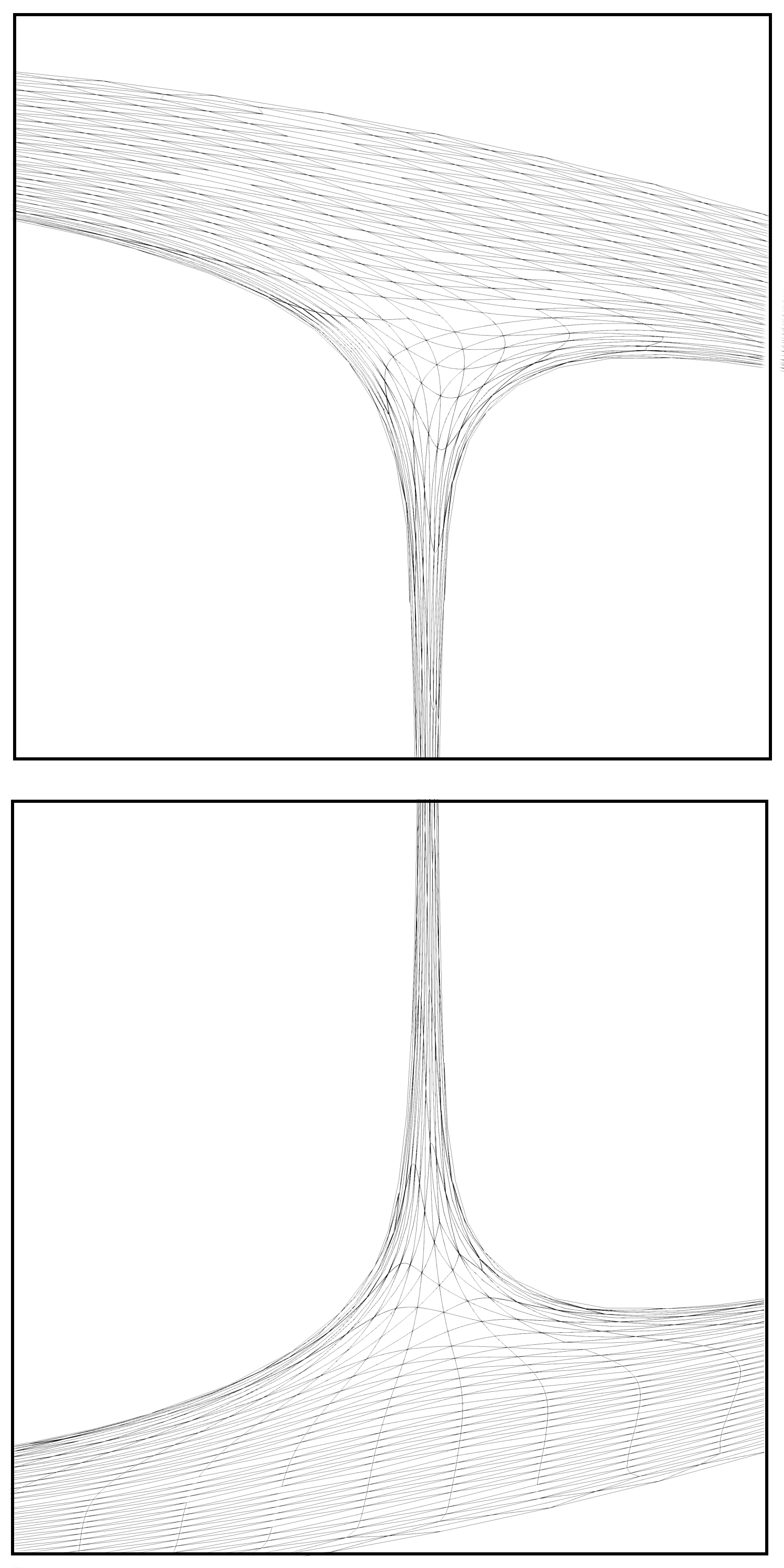}}
	\quad
	\subfloat[$k = 3$]{
		\includegraphics[width=0.25\textwidth]{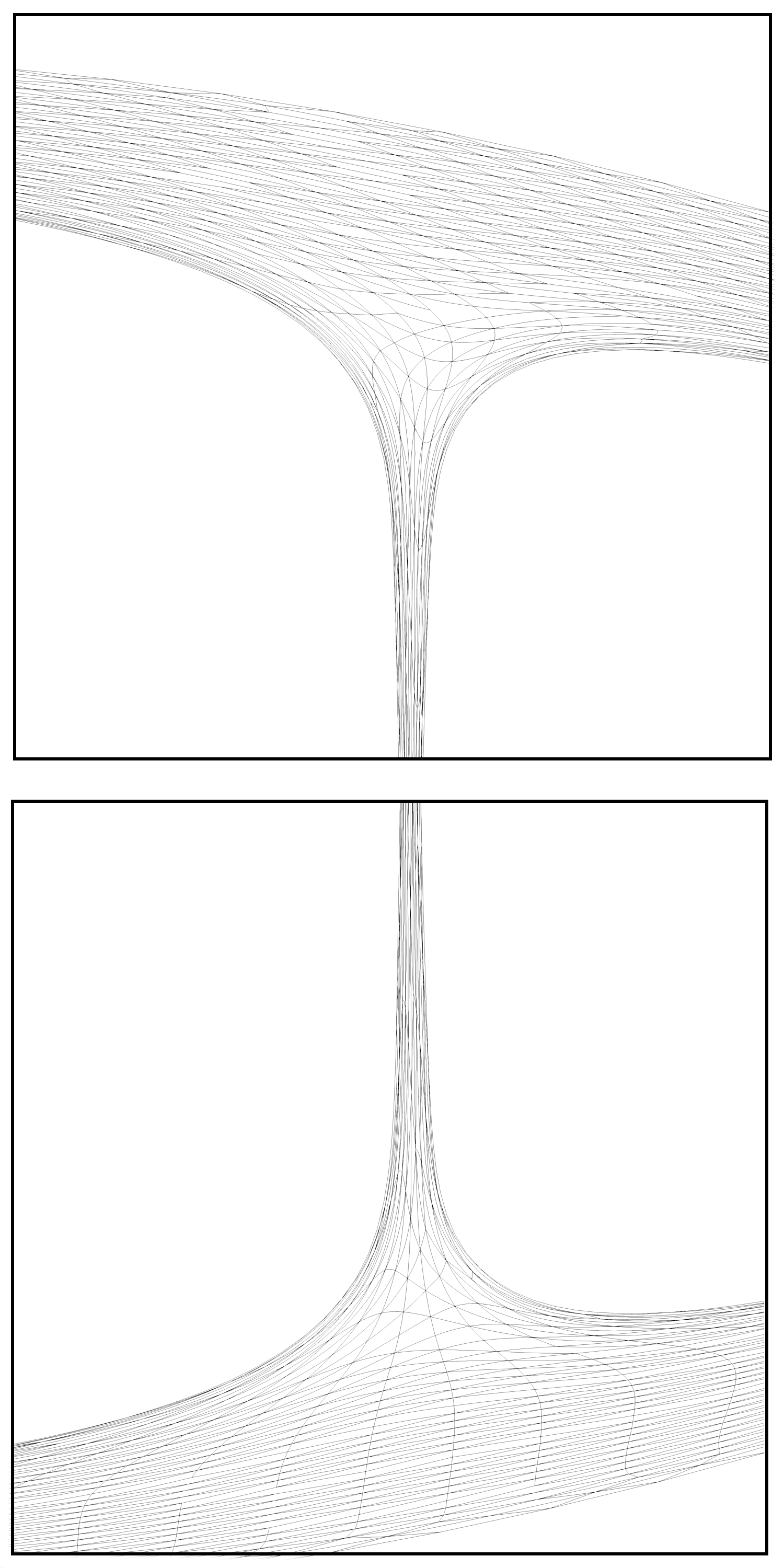}}
	\quad
	\caption{ NHK-CAV disk with two holes subjected to tensile stresses.
                  Details of the strip separating the two holes at different polynomial degrees $k = \{ 1,2,3 \}$. 
	\label{fig:circle_voids_details}}
\end{figure}

The adaptive stabilization parameters are set as $\beta = 1$ and $\epsilon = 1$, the latter helping Newton's method convergence rates. 
This test case requires 100 incremental steps with Dirichlet boundary conditions imposed by means of the Lagrange multipliers method, 
as opposite, we were unable to succeed with Dirichlet BCs imposed by means of Nitsche method. 
As reported in \Tab~\ref{tab:increments_simulations}, 100 steps is the highest number of increments required among all test cases employing Lagrange multipliers.
This confirms that the test case challenges the robustness of the numerical strategy. 
The final highly distorted computational mesh, consisting of 8982 triangular elements, is depicted in \Fig~\ref{fig:circle_voids_details}.

\subsection{3D simulations} \label{sec:3D_simulations}
\subsubsection{Torsion of a square section bar} \label{sec:beam_torsion}
\begin{figure}[H]
	\hspace{-0.4cm}
	\subfloat[3D Mesh]{
	\includegraphics[width=0.15\textwidth]{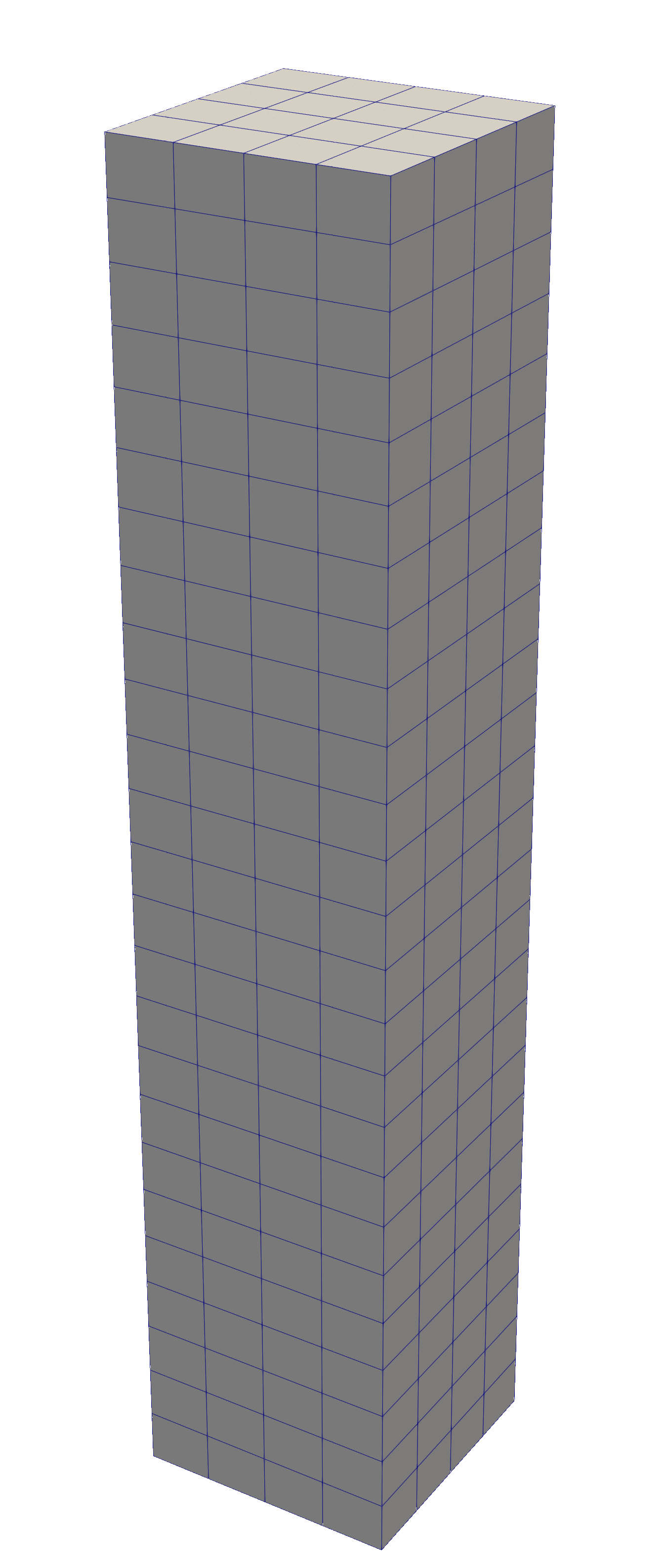}}
	\subfloat[$\alpha = \frac{\pi}{2}$\newline $k=1$]{
		\includegraphics[width=0.15\textwidth]{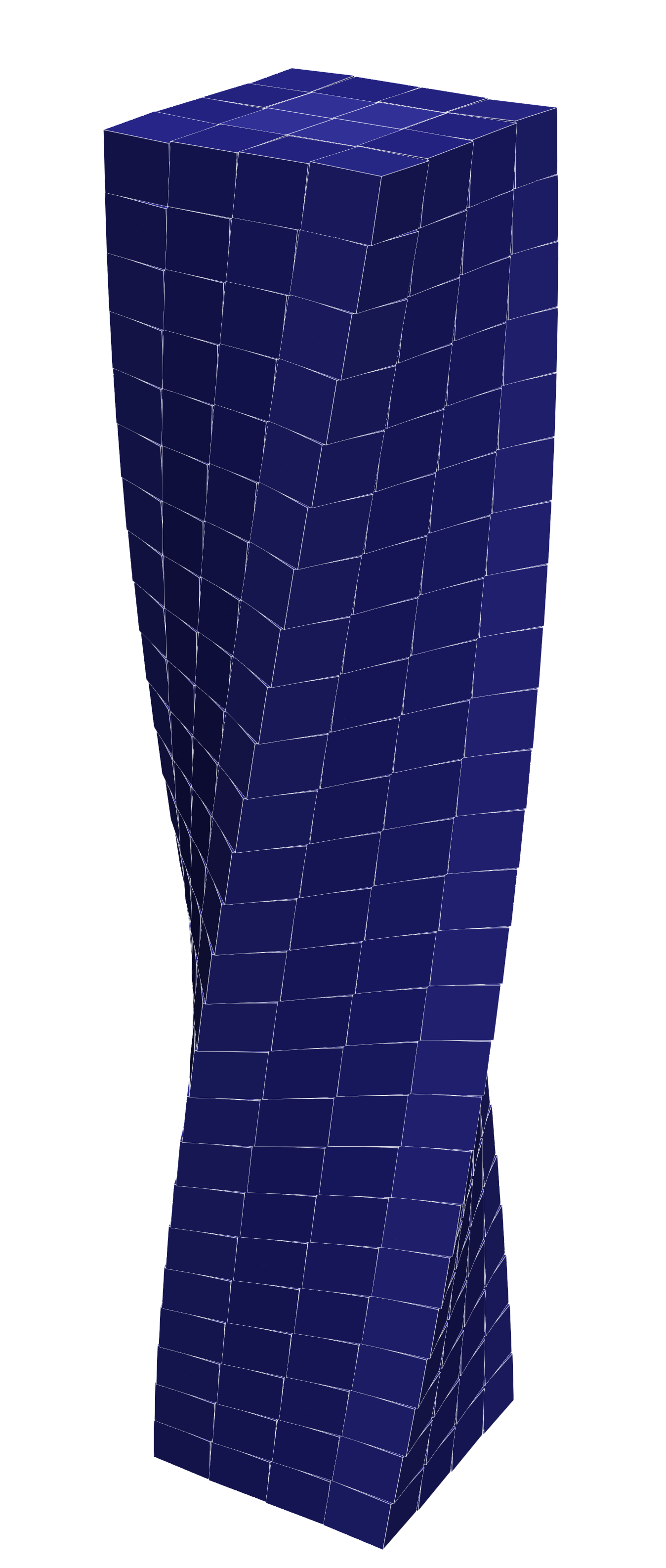}}
	\subfloat[$\alpha = \frac{3 \pi}{2}$\newline $k=1$]{
	\includegraphics[width=0.15\textwidth]{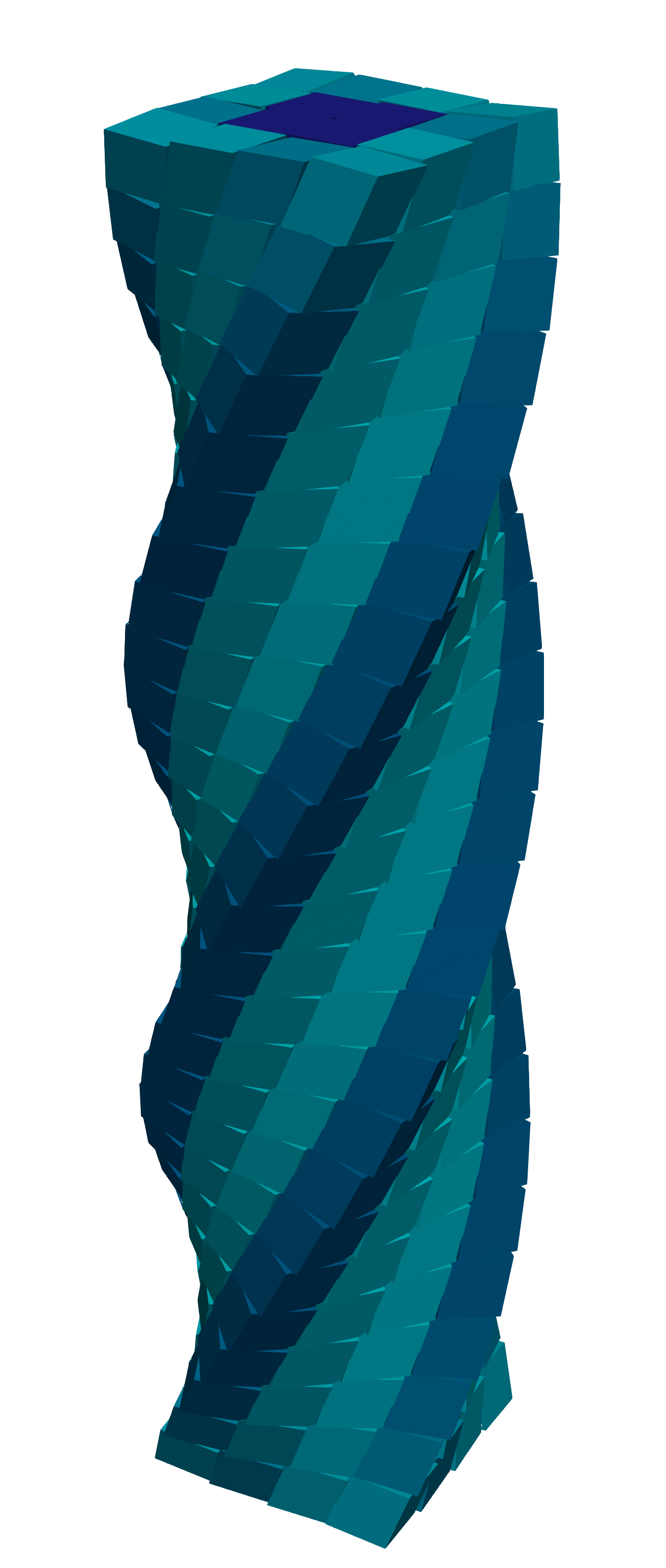}}
	\subfloat[$\alpha = 2 \, \pi$\newline $k=1$]{
 	\includegraphics[width=0.15\textwidth]{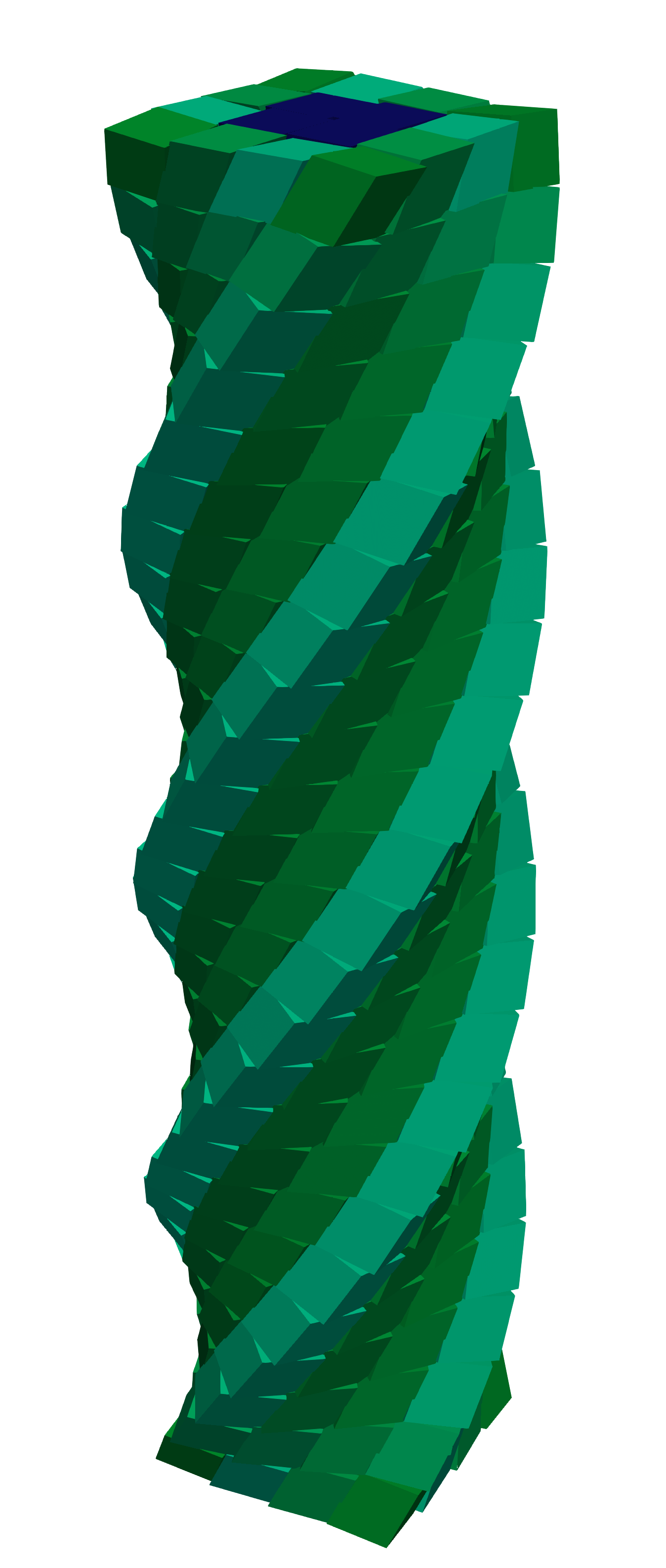}}
    \subfloat[$\alpha = 2 \, \pi$\newline $k=2$]{
	\includegraphics[width=0.15\textwidth]{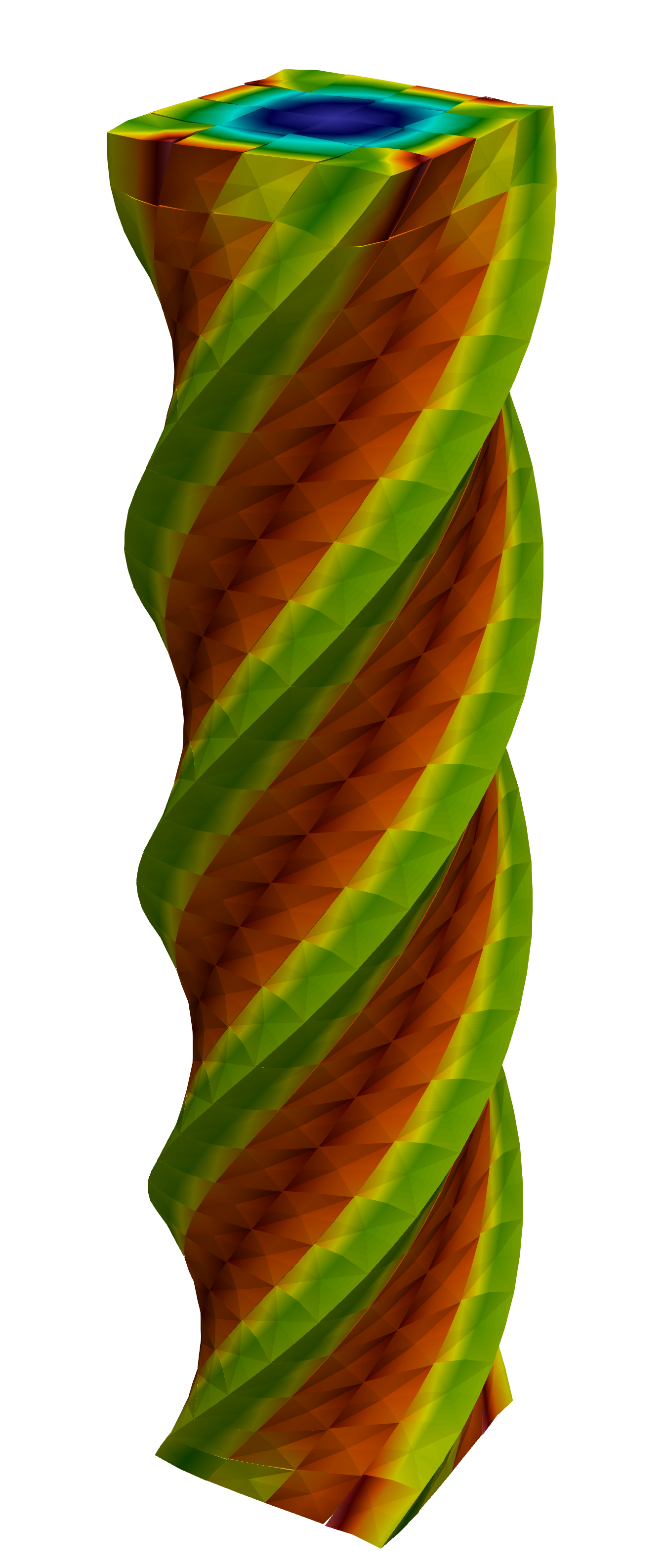}}
    \subfloat[$\alpha = 2 \, \pi$\newline $k=3$]{
	\includegraphics[width=0.15\textwidth]{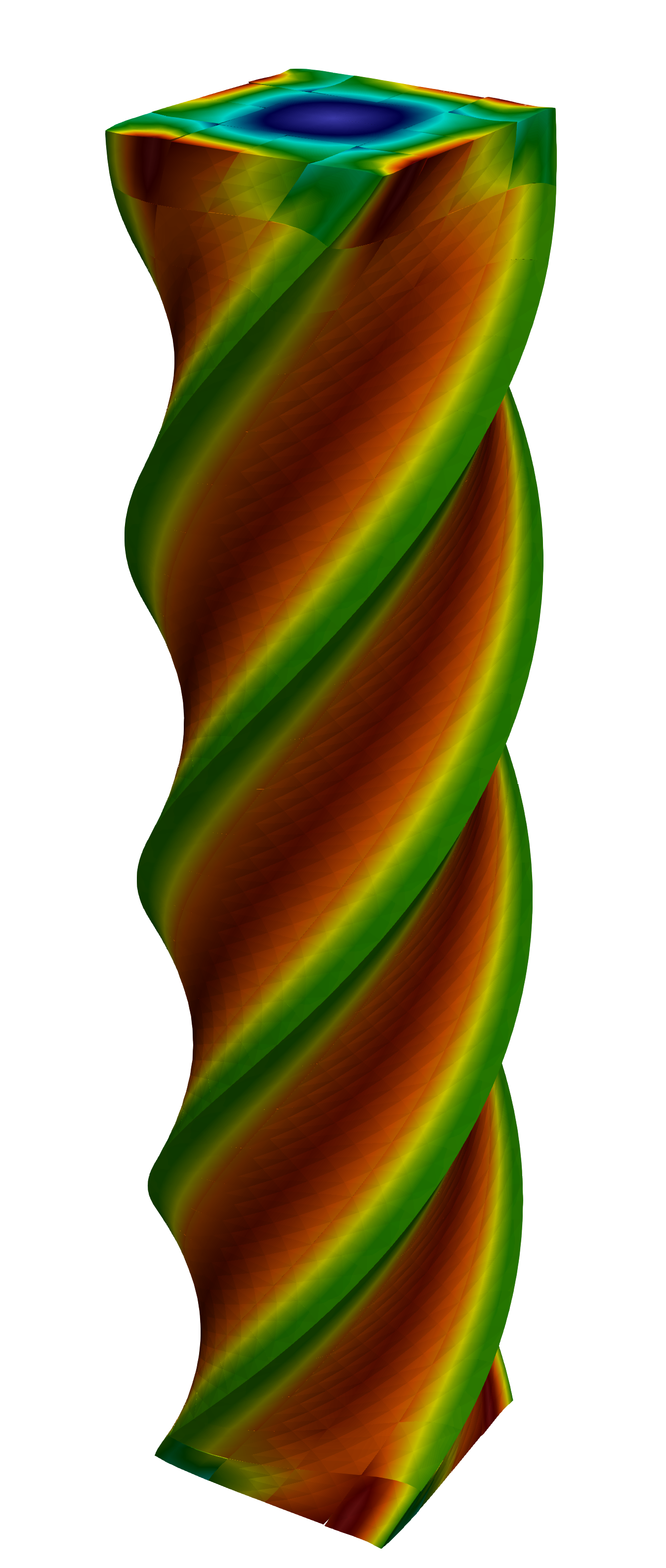}}
	\subfloat{
	\includegraphics[width=0.1\textwidth]{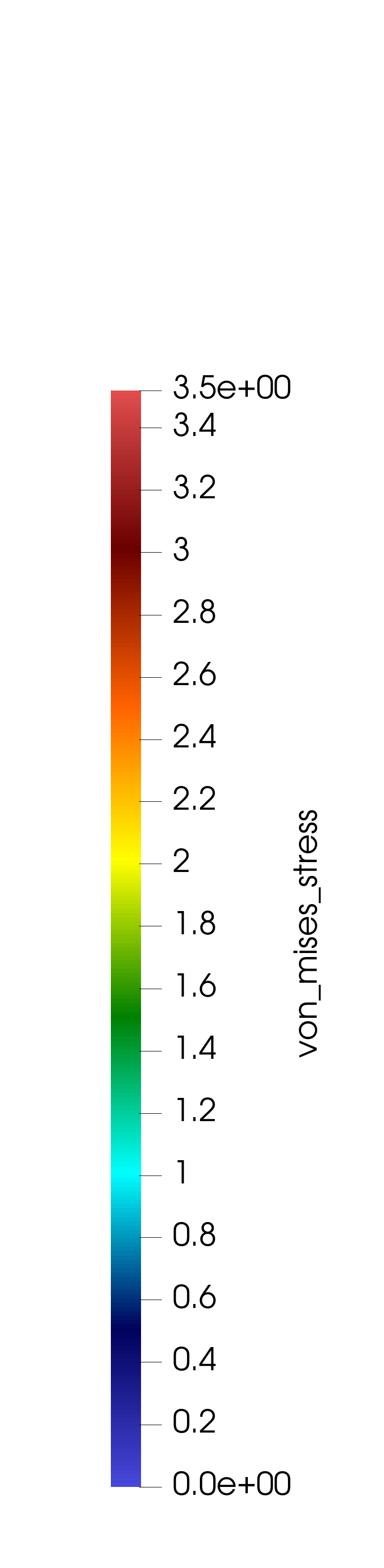}}
	
	\caption{Torsion of a SVK-I square section bar. (a) Computational mesh; (b-f) Von Mises stress distribution for $k=1,2,3$ when the top surface is rotated by an angle $\alpha$. 
	\label{fig:bar_torsion_inc}}
\end{figure}

We consider a square section bar such that $H/L = 5$, where $L$ is the edge length of the square cross-section and H is the extension of the bar in the axial direction. 
The mesh, consisting of 400 uniform hexahedral elements, is shown in \Fig~\ref{fig:bar_torsion_inc}. 
The bottom surface is clamped while the top surface is subjected to a 360 degrees plane rotation around its centroid. 
We employ the fully incompressible SVK-I model with material parameters $\mu=1$ and $\lambda=1$. 
We remark that this test case was successfully completed also based on the SVK-C, NHK-C and NHK-I constitutive laws but the results are not presented for the sake of conciseness.
Dirichlet boundary conditions are imposed on the top and bottom surfaces while homogeneous Neumann boundary conditions are enforced on the rest of the boundary.
Regarding the stabilization parameters, adaptive stabilization is mandatory only in the NHK-C case: we set $\beta=1$ for $k=1$, $\beta=2$ for $k=2$ and $\beta=5$ for $k=3$. 
Other relevant parameters are $\epsilon = 0$ and $\eta_{\text{LBB}} = 1$ in the incompressible regime. 

The results for the SVK-I model are displayed in \Fig~\ref{fig:bar_torsion_inc} considering first, second and third degree BR2 dG discretizations. 
Increasing the polynomial degree reduces the amplitude of discontinuities in the displacement field resulting in a more precise representation of the geometry of the deformed bar. 
Furthermore, the regions where the stresses intensify are more accurately captured: note that, for $k=1$, the stress is constant inside each mesh element. 

\subsubsection{Cylinder deformation} \label{sec:cylinder_rotation}
We consider a hollow cylinder such that $H/R = 4$ and $r = 0.7R$, where $R$ and $r$ are, 
respectively, the external and internal radius of the annulus cross-section and $H$ is the extension of the cylinder in the axial direction.
As proposed in \cite{Eyck2008}, the top surface of the cylinder is rotated while keeping the bottom surface clamped.
We consider both NHK-C and NHK-I constitutive laws with $\nu =0.25$ and $E = 1$. 
Dirichlet boundary conditions are imposed on the top and bottom surfaces while homogeneous Neumann boundary conditions are enforced on the rest of the boundary.

\begin{figure}[H]
	\subfloat[$\alpha = 0$]{
		\hspace{-1.cm}
		\includegraphics[width=0.3\textwidth]{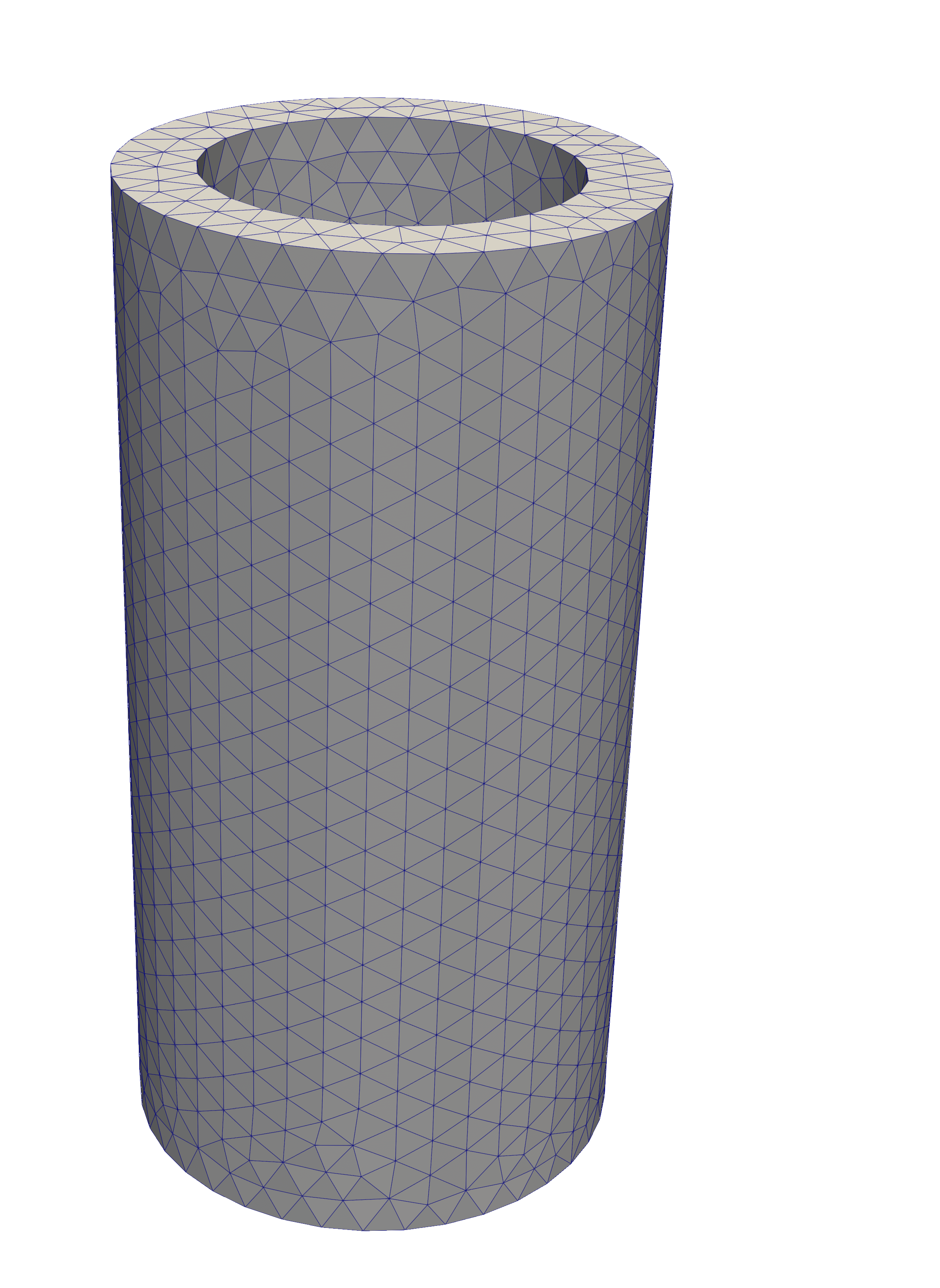}}
	\subfloat[$\alpha = 0.15 \pi$]{
		\hspace{-1.5cm}
		\includegraphics[width=0.3\textwidth]{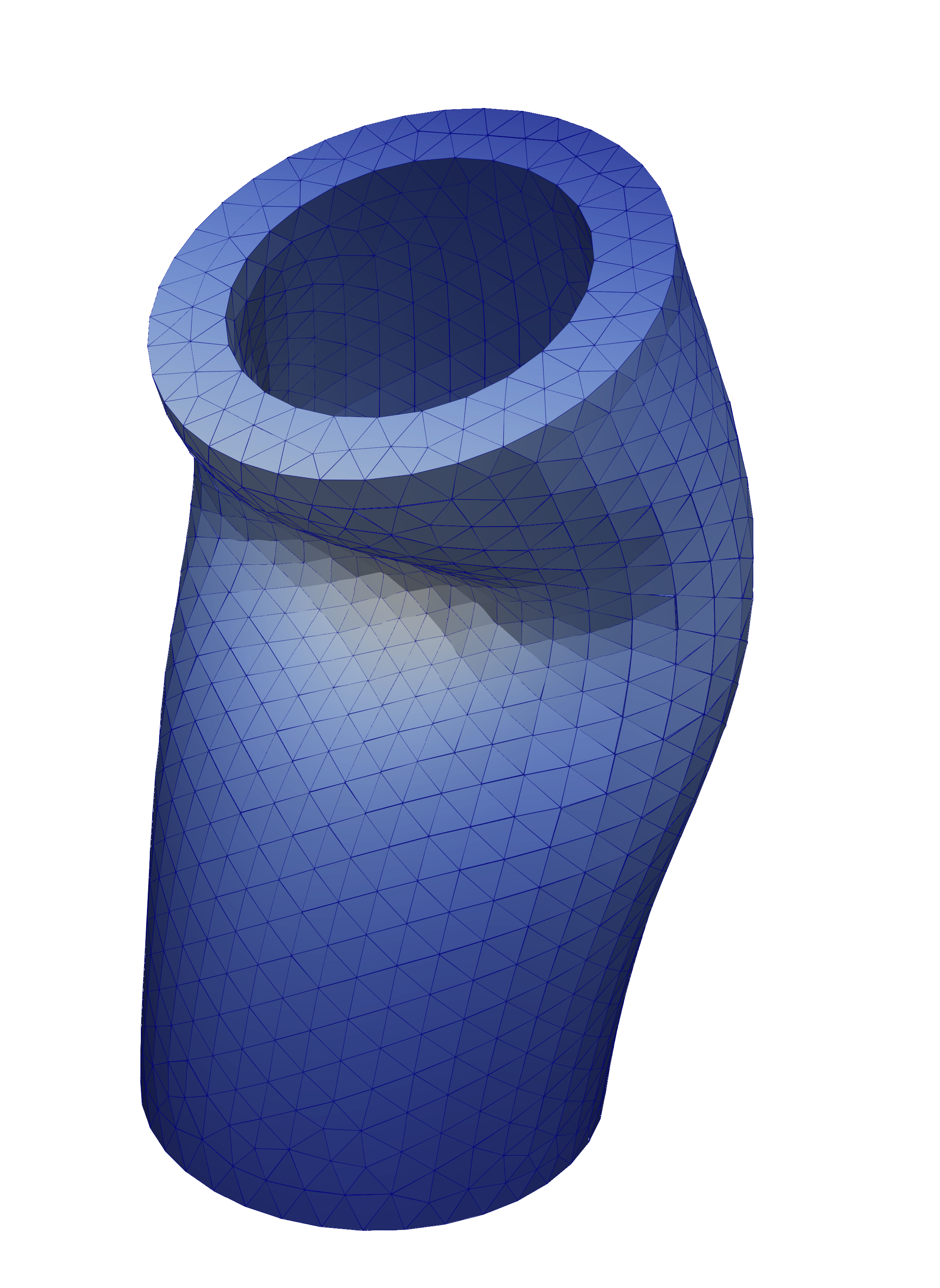}}
	\subfloat[$\alpha = 0.30 \pi$]{
		\hspace{-1.5cm}
		\includegraphics[width=0.3\textwidth]{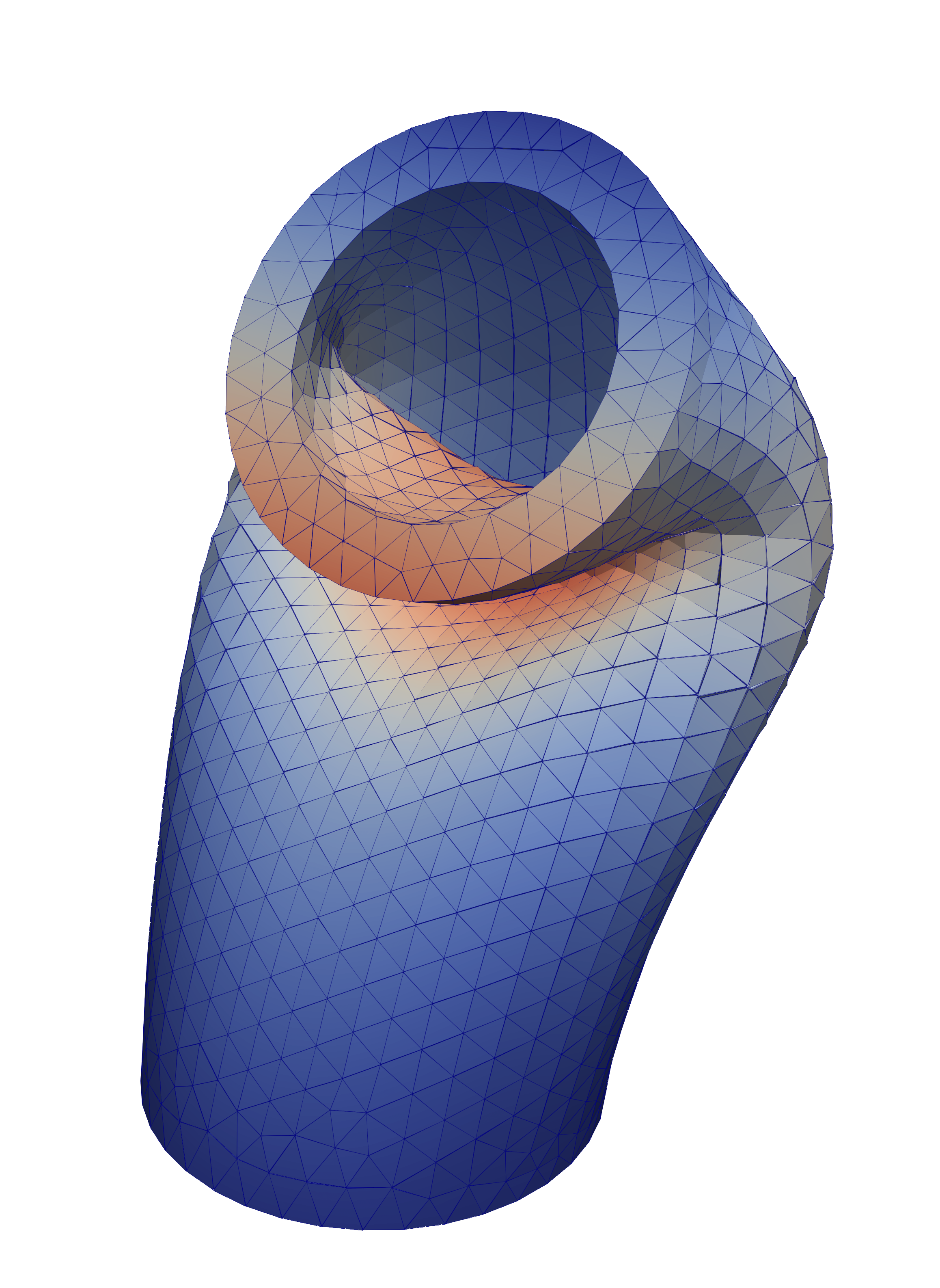}}
	\subfloat[$\alpha = 0.5 \pi$]{
		\hspace{-1.5cm}
		\includegraphics[width=0.3\textwidth]{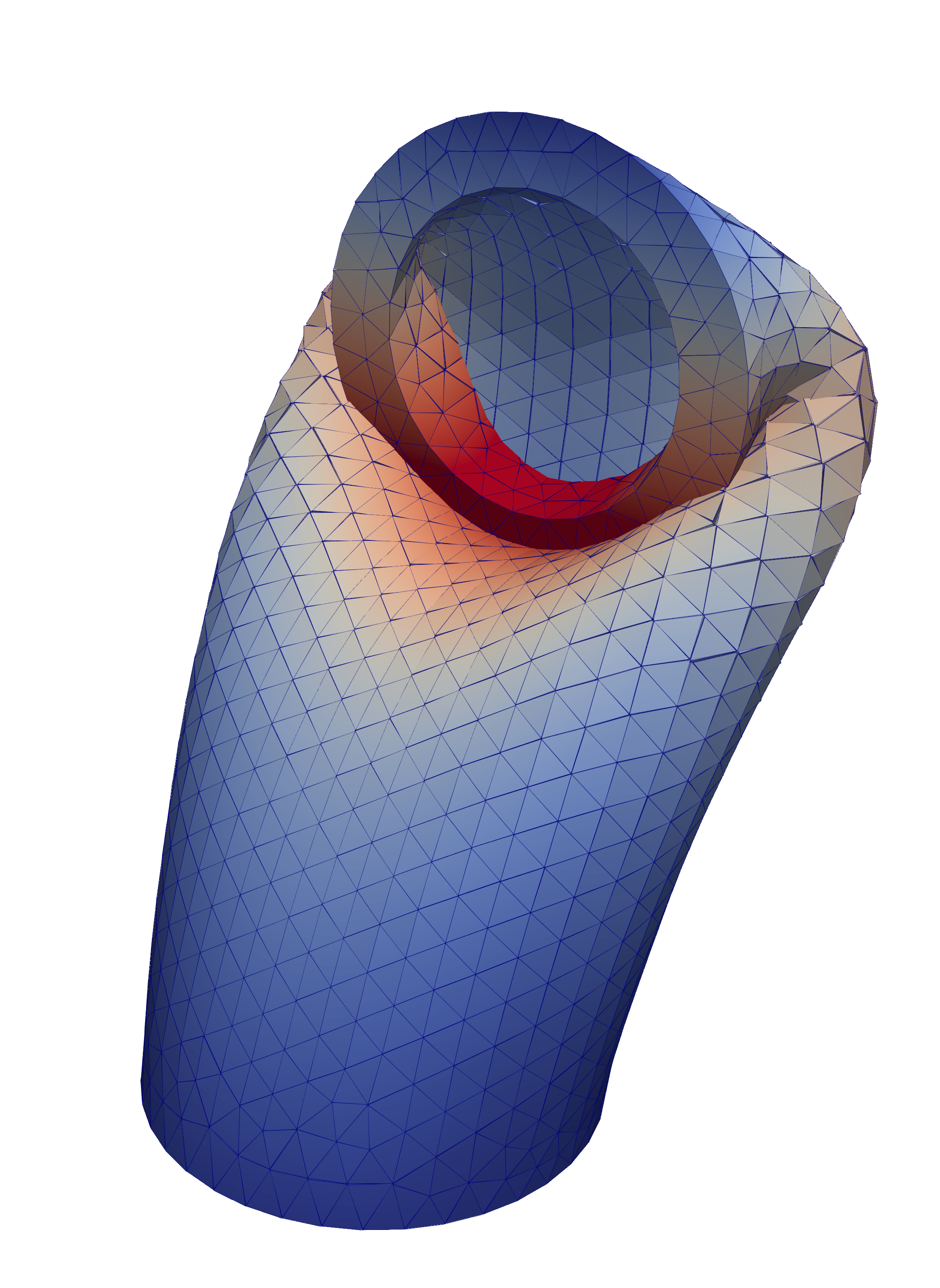}}
	\subfloat{
		\hspace{-0.9cm}
		\includegraphics[width=0.15\textwidth]{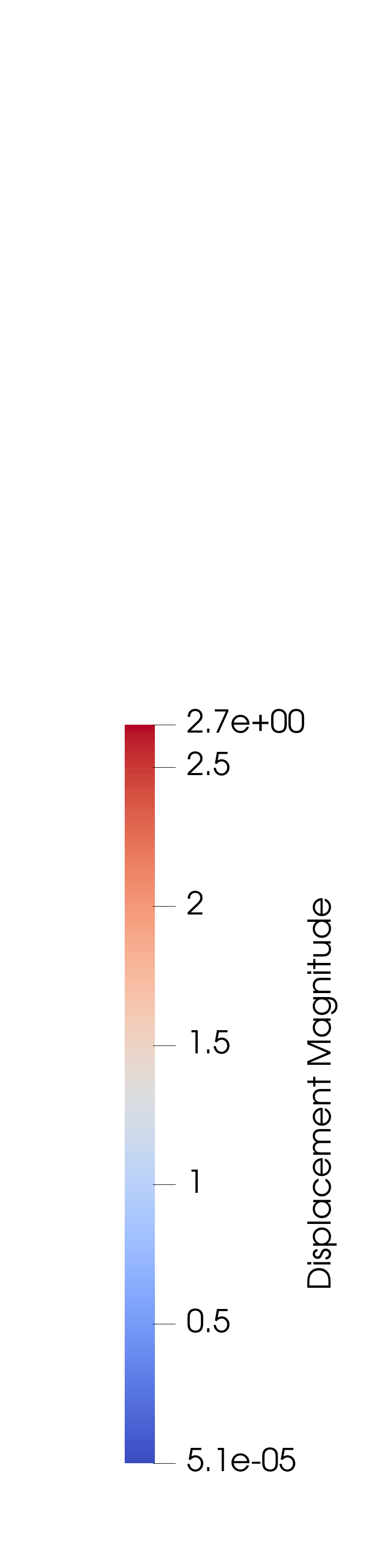}}
	
	\caption{Rotation of the top surface of a NHK-C hollow cylinder: 
                 sequence of equilibrium states obtained by the incremental load method while increasing the rotation angle $\alpha$, $k=1$.
	\label{fig:cylinder_rot_comp}}
\end{figure}

\Fig~\ref{fig:cylinder_rot_comp} reports the computational mesh, consisting of 8906 tetrahedral elements, and
the deformed configurations at different rotation angles $\alpha$, with $0 \leq \alpha \leq \frac{\pi}{2}$,
obtained for the compressible case with a first degree BR2 dG formulation. 
It is worth mentioning that for $\alpha \geq \frac{\pi}{4}$ the cylinder penetrates itself due to the lack of contact boundary conditions. 
Despite the lack of meaningfulness from the physical viewpoint, this result emphasises the capability of dealing with large deformations. 

\begin{figure}[H]
	\centering
	\subfloat[NHK-C]{
		\includegraphics[width=0.3\textwidth]{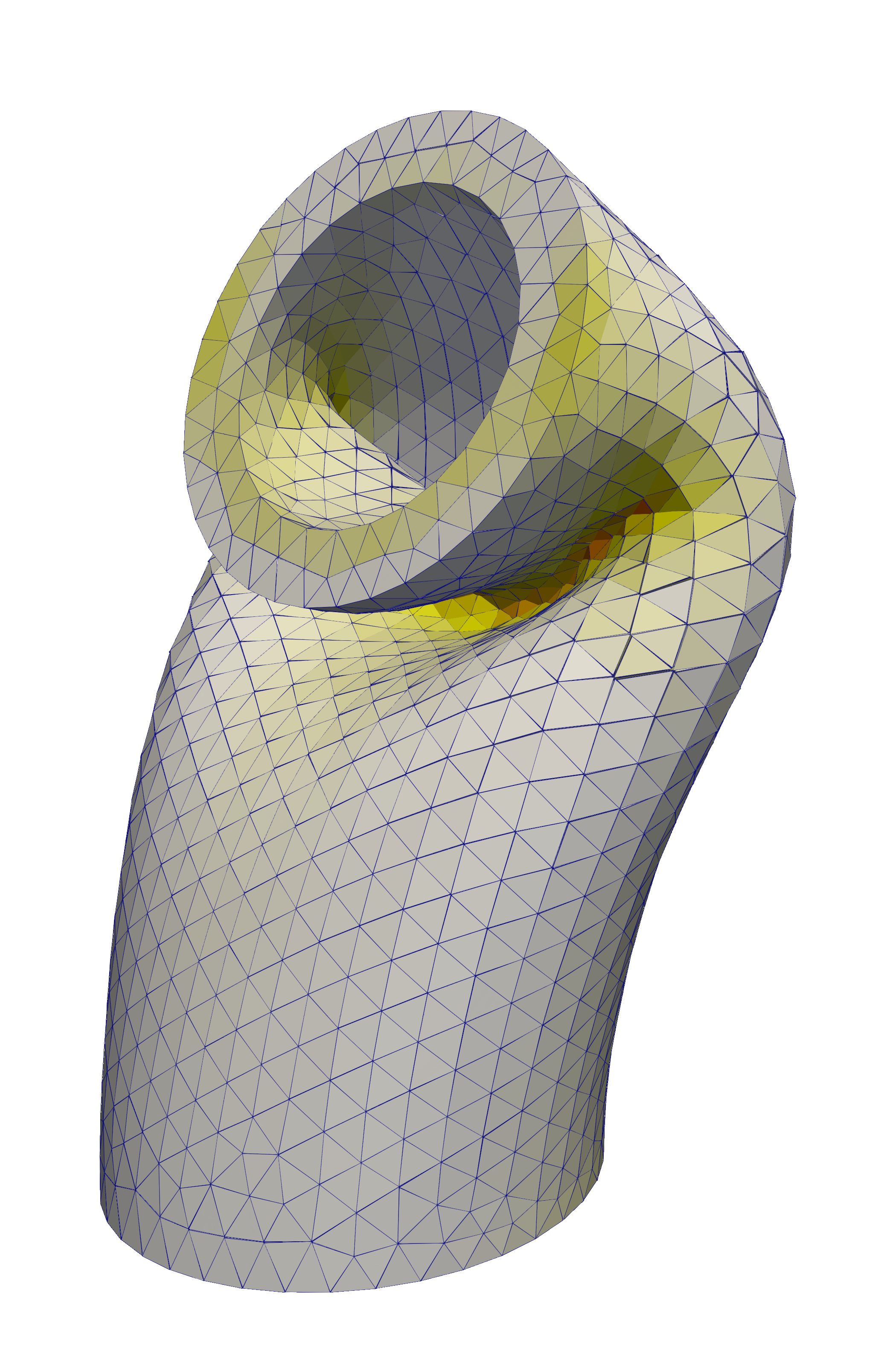}\label{fig:cylinder_rot_comp_eigenvalues}}
	\subfloat[NHK-I]{
		\includegraphics[width=0.3\textwidth]{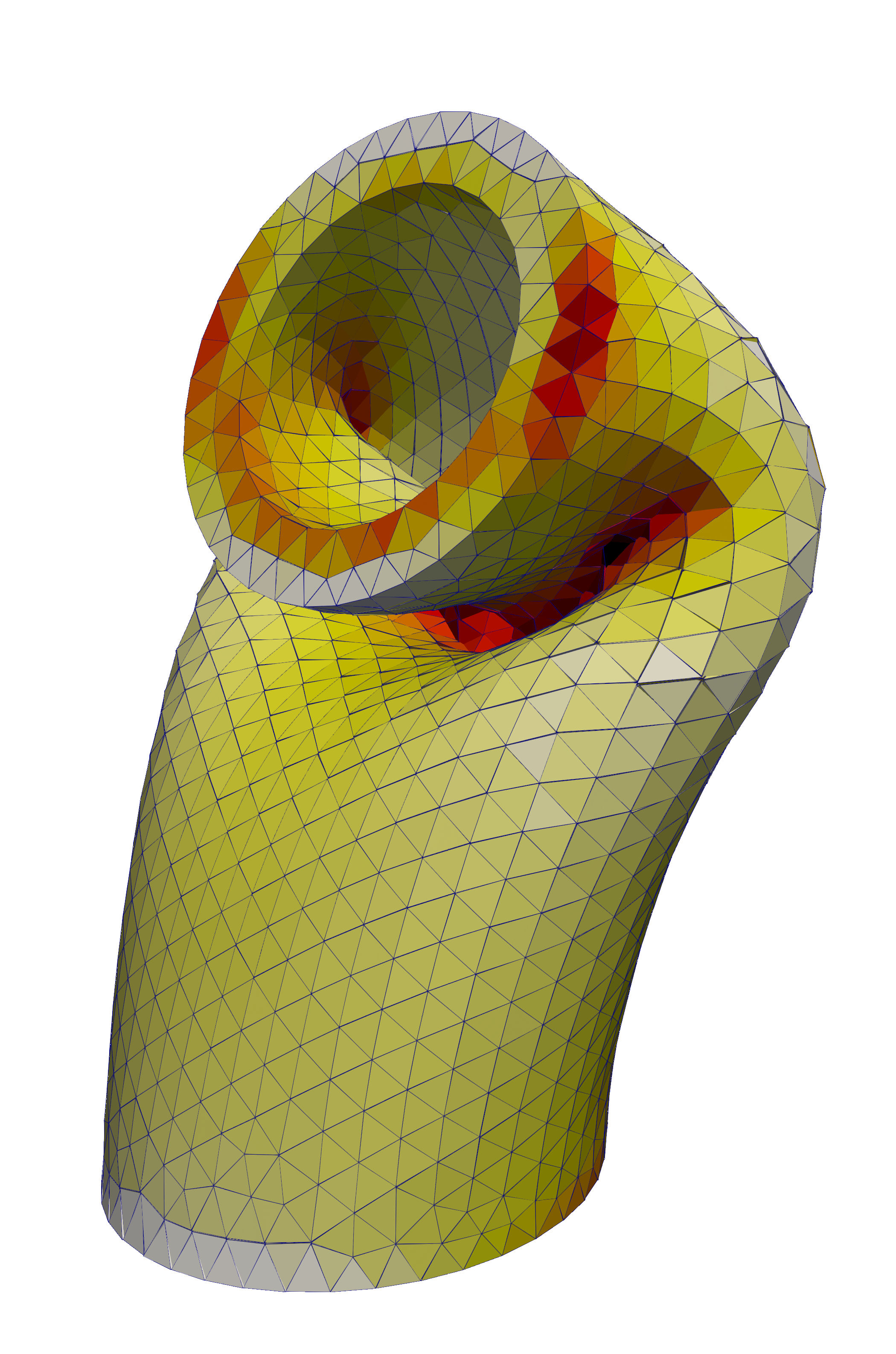} \label{fig:cylinder_rot_inc_eigenvalues}}
	\subfloat{
		\includegraphics[width=0.12\textwidth]{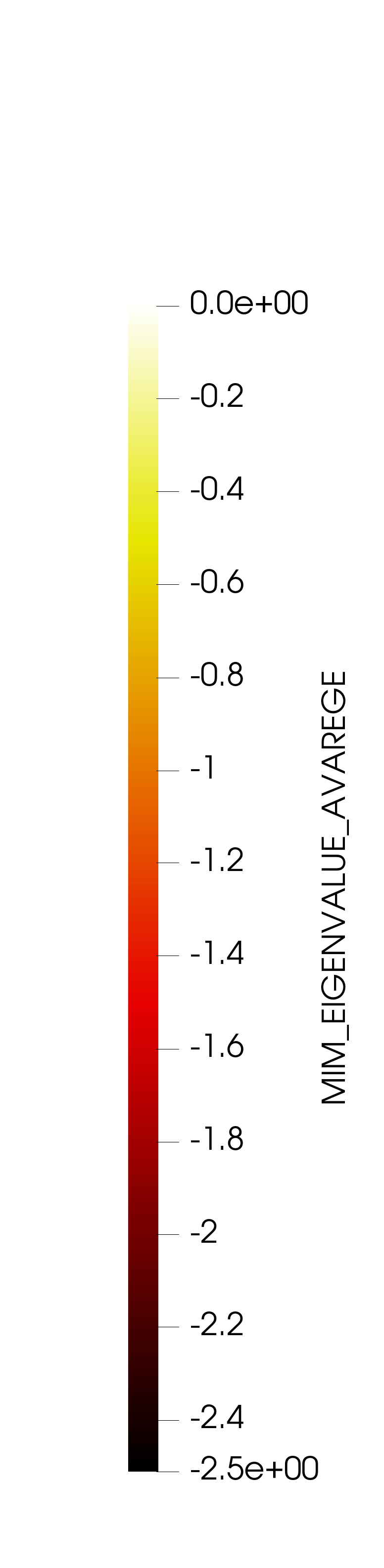}}
	\caption{Deformation of a hollow cylinder. \emph{Left} and \emph{right}: NHK-C and NHK-I constitutive laws, respectively. 
                 Images are colour-coded based on the minimum negative eigenvalue of the fourth-order elasticity tensor $\Tnsr{A}$.
	\label{fig:cylinder_rot_comp_vs_inc_eigenvalues}}
\end{figure}

The stabilization parameters reads $\beta = 4$ and $\epsilon = 1$ and $\eta_{\text{LBB}} = 1$ for the incompressible model. 
In \Fig~\ref{fig:cylinder_rot_comp_vs_inc_eigenvalues}, deformed states are colour-coded with minimum negative eigenvalues of the fourth-order elasticity tensor 
allowing to appreciate that compression of the beam material triggers the adaptive stabilization strategy. 
\Fig~\ref{fig:cylinder_rot_inc_eigenvalues} also depicts the deformed state reached 
at 65\% of the entire rotation by the NHK-I cylinder. 
After approaching this configuration Newton's method struggles to converge irrespectively of the amount of stabilization introduced.

\begin{figure}[H]
\centering
	\subfloat[$\epsilon = 0$]{ \includegraphics[width=0.5\linewidth]{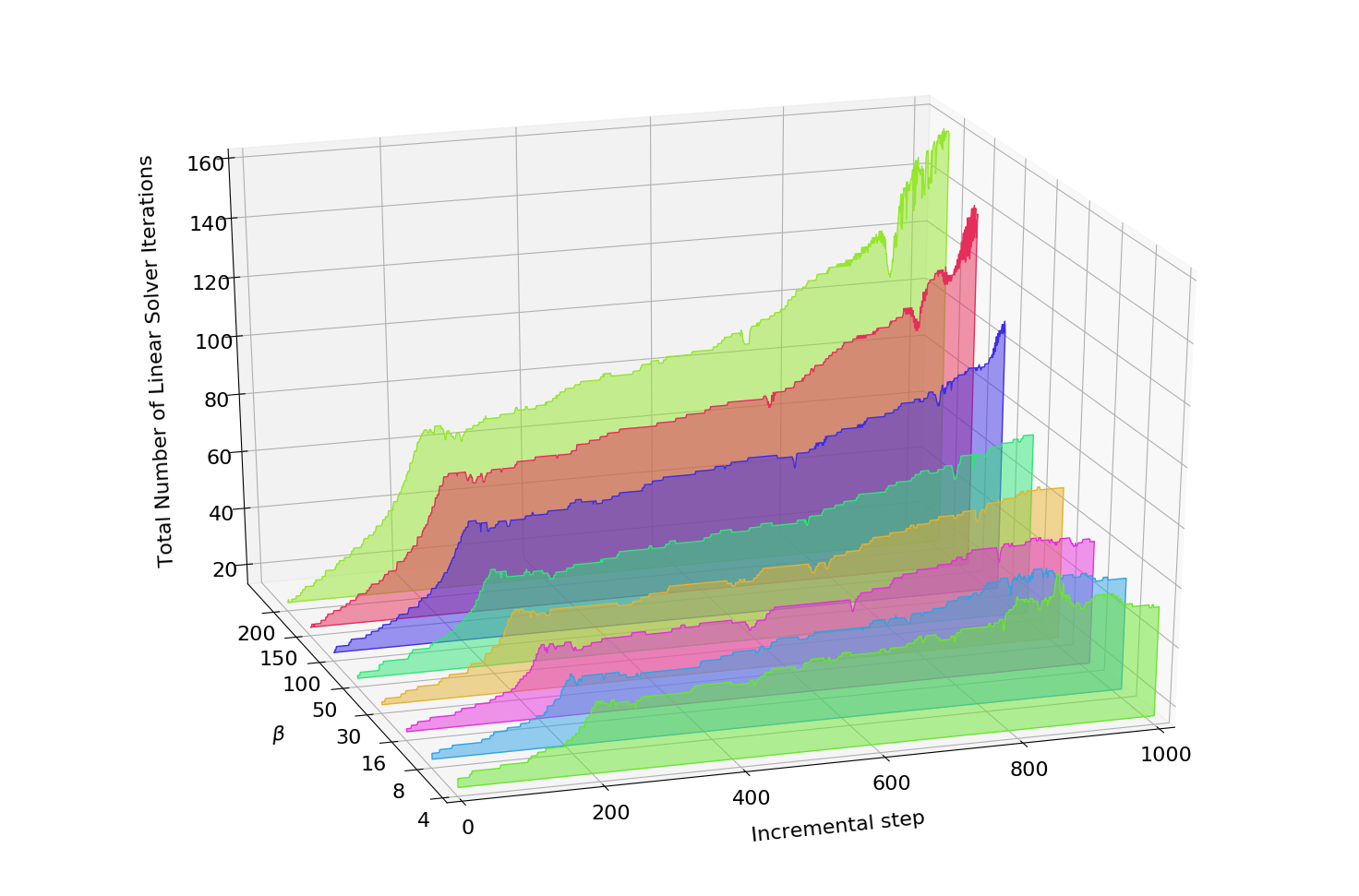} }
	\subfloat[$\epsilon = 1$]{ \includegraphics[width=0.5\linewidth]{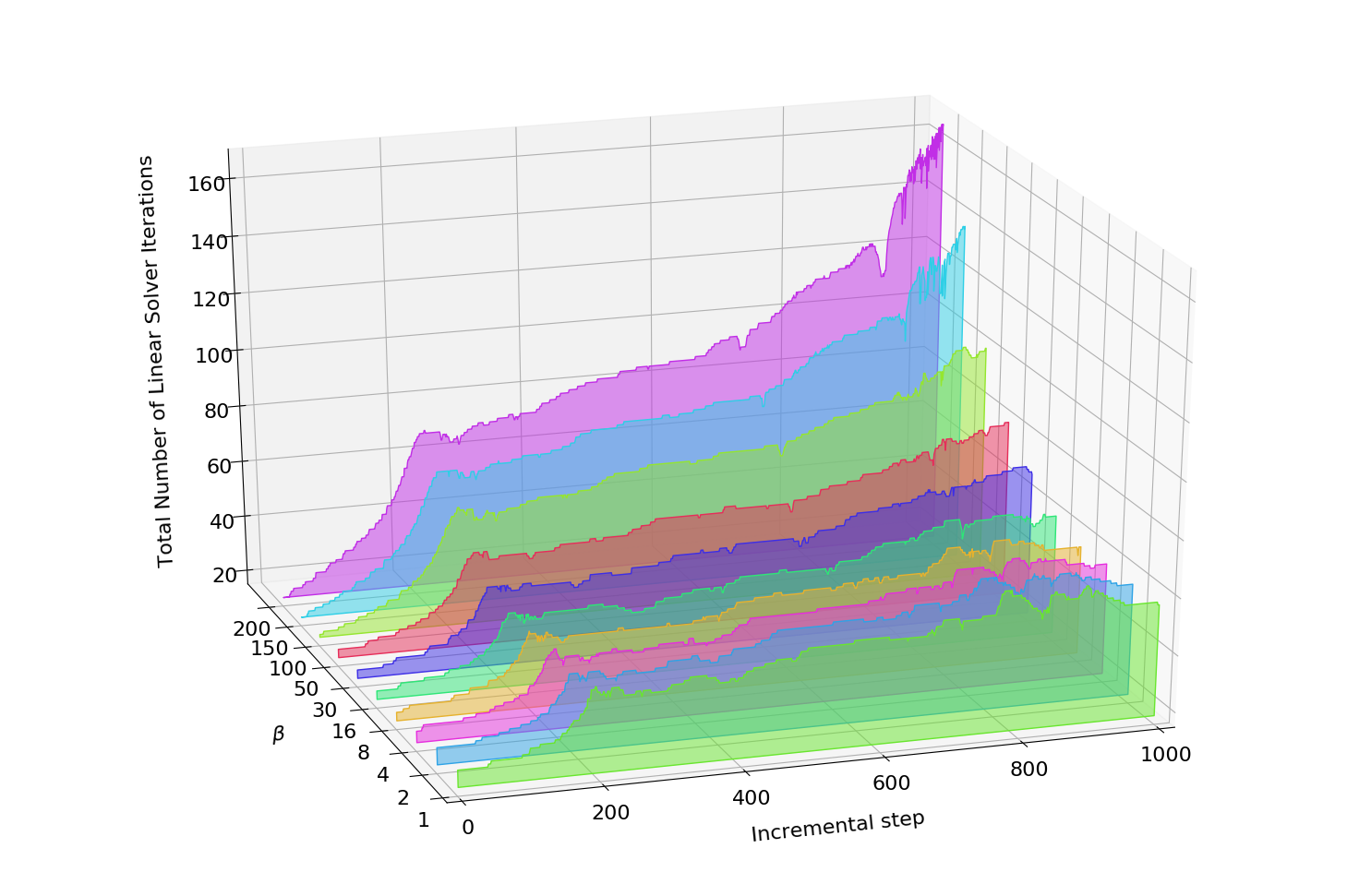} }
	\vspace{-1cm}	
	\subfloat[$\epsilon = 10$]{ \includegraphics[width=0.5\linewidth]{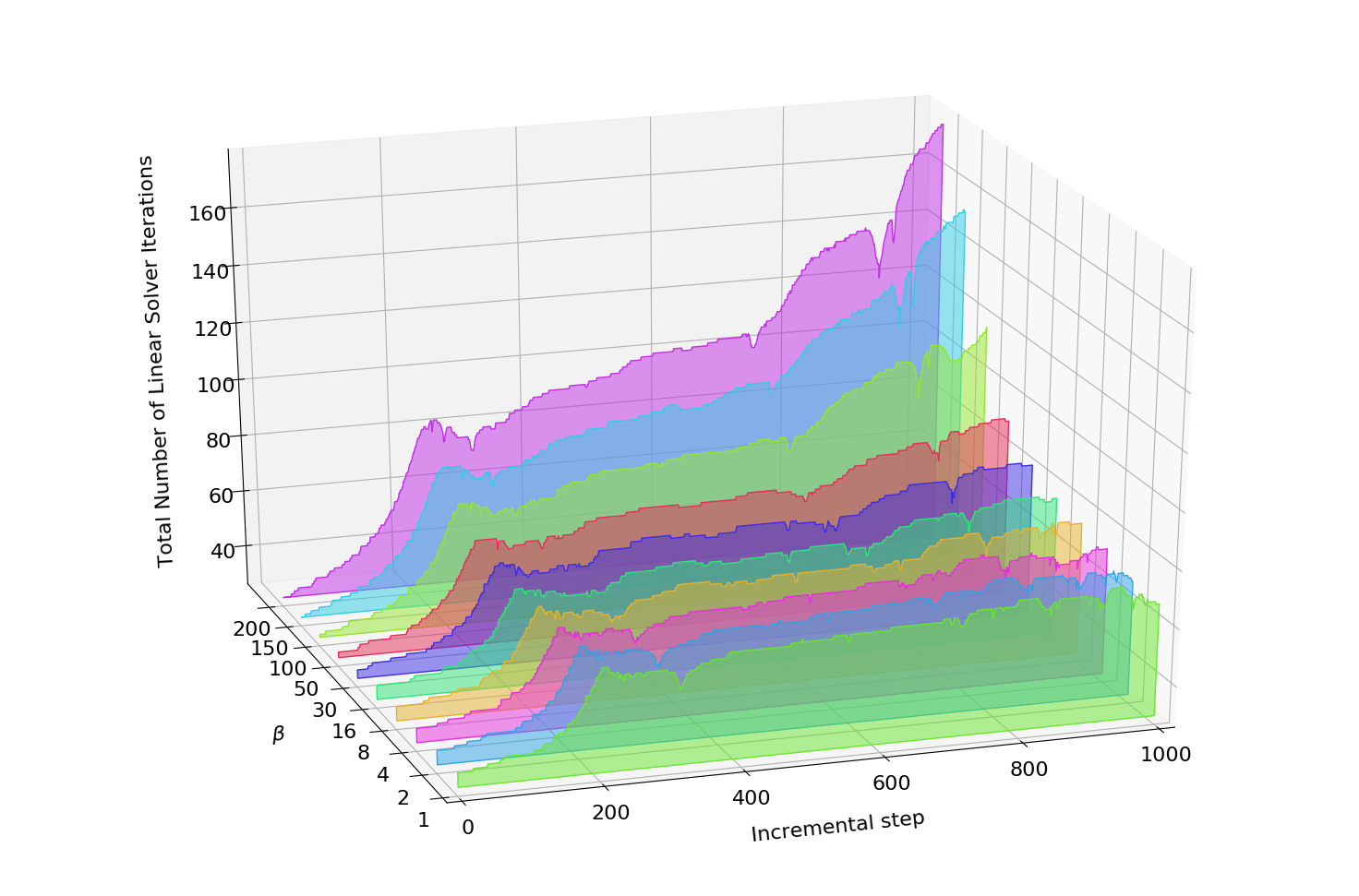} }
	\subfloat[$\epsilon = 20$]{ \includegraphics[width=0.5\linewidth]{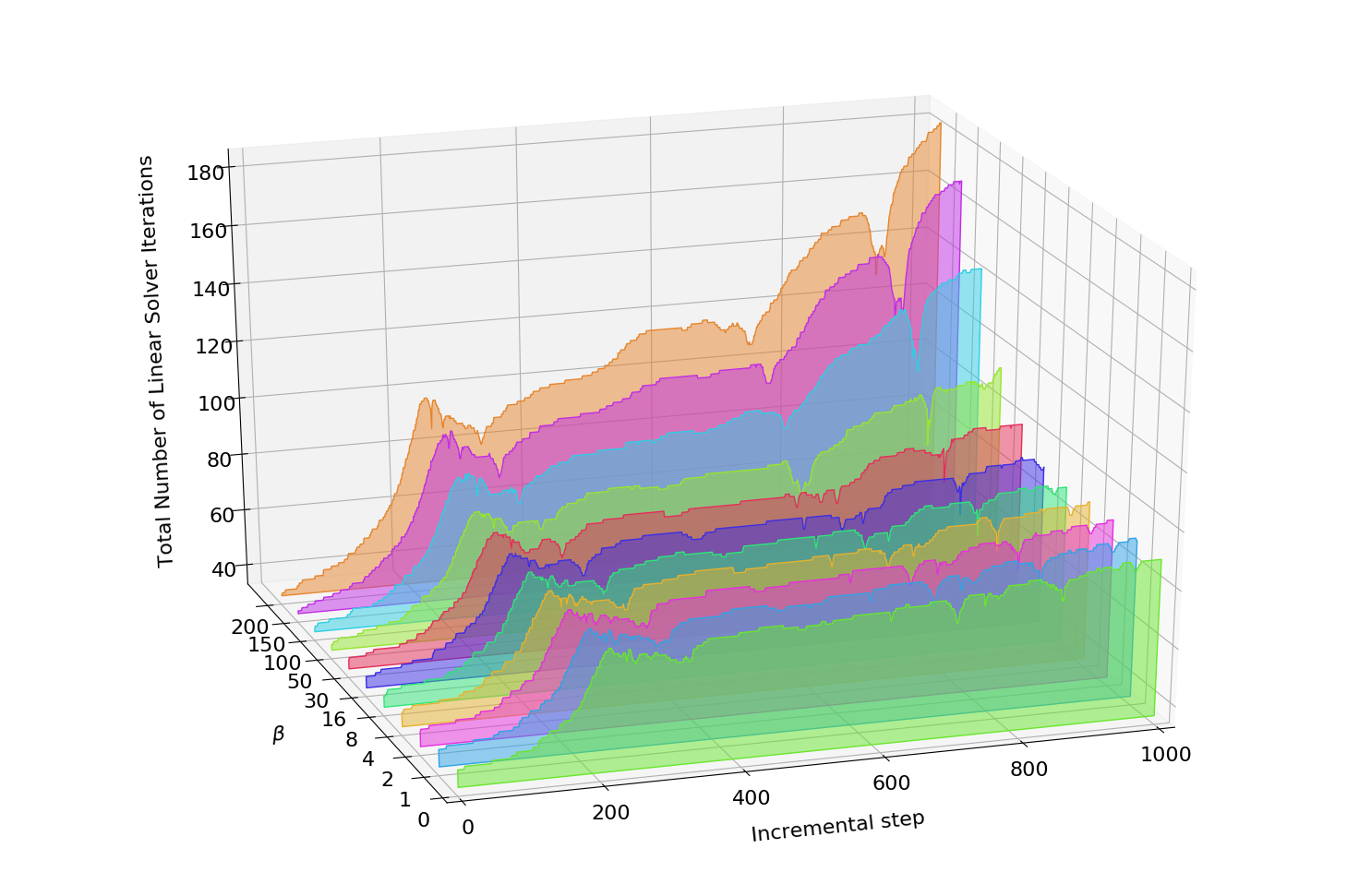} }
	\caption{NHK-C cylinder: total number of linear solver iterations recorded along the loading path (1k incremental steps). 
             Results are obtained varying the stabilization parameters $\beta$ and $\epsilon$ in order to show their influence on the performance of the solution strategy.
	\label{fig:cylinder_rot_beta_tests}}
\end{figure}

As we did for the beam deformation of \Sec~\ref{sec:beam_deformation}, 
we analyse the influence of the stabilization parameter on the performance of the $h$-multigrid solution strategy.
We consider a 1k increments loading path and choose $\beta \in [0,200]$ and $\epsilon \in \left\lbrace 0, 1, 10, 20 \right\rbrace $. 
In \Fig~\ref{fig:cylinder_rot_beta_tests}, for each incremental step solved by Newton's method,
we report the total number of linear solver iterations obtained varying $\beta$ and $\epsilon$. The average and maximum number of Newton iteration as well as the average and maximum number of linear solver iterations are tabulated in \Tab~\ref{tab:cylinder_rot_beta_epsilon_tests}.

The iteration spike observed at around one fifth of the loading path is due to buckling of the cylinder. 
We remark that for $\beta < 4$ and $\epsilon = 0$ the computation fails due to an insufficient amount of stabilisation. 
In the range $4 \leq \beta \leq 50$, the runs are successful and we observe a mild increase of the number of linear solver iterations. 
An excessive amount of stabilization ($\beta > 50$) deteriorates the solver efficiency leading to a significant increase of the computational time. 

\begin{table}[H]
	\centering
	\begin{tabular}{|c|l|c|c|c|c|c|c|c|c|c|c|c|c|}
		\hline 
		\multirow{2}{*}{\begin{tabular}[l]{@{}l@{}}$\epsilon$ \end{tabular}}  & \multicolumn{2}{|l|}{\multirow{2}{*}{}} & \multicolumn{11}{c|}{$\beta$}\\ \cline{4-14}
		&  \multicolumn{2}{|l|}{}  &  \multicolumn{1}{l|}{\textbf{0}}  &  \multicolumn{1}{l|}{\textbf{1}}  &  \multicolumn{1}{l|}{\textbf{2}}  &  \multicolumn{1}{l|}{\textbf{4}}  &  \multicolumn{1}{l|}{\textbf{8}}  &  \multicolumn{1}{l|}{\textbf{16}}  &  \multicolumn{1}{l|}{\textbf{30}}  &  \multicolumn{1}{l|}{\textbf{50}}  &  \multicolumn{1}{l|}{\textbf{100}}  &  \multicolumn{1}{l|}{\textbf{150}}  &  \multicolumn{1}{l|}{\textbf{200}}\\ 
		\hline
		\hline
		\multirow{4}{*}{\begin{tabular}[l]{@{}l@{}} 0 \end{tabular}} & \multirow{2}{*}{\begin{tabular}[l]{@{}l@{}}\textbf{Newton} \\  \textbf{iterations}\end{tabular}}   & \textbf{mean} & -& -& -  &  5  &  5  &  5  &  4  &  4  &  4  &  4  &  4\\ \cline{3-14}
		&  & \textbf{max} & -& -& -  &  5  &  5  &  5  &  4  &  4  &  4  &  4  &  4\\ \cline{2-14}
		&  \multirow{2}{*}{\begin{tabular}[l]{@{}l@{}}\textbf{Linear Solver} \\ \textbf{iterations} \end{tabular}}   & \textbf{mean} & -& -& - &  11  &  11  &  11  &  12  &  13  &  16  &  22  &  26\\ \cline{3-14}
		&  & \textbf{max}  & -& -& -  &  \multicolumn{1}{c|}{20}    &  \multicolumn{1}{c|}{20}    &  \multicolumn{1}{c|}{22}    &  \multicolumn{1}{c|}{25}    &  \multicolumn{1}{c|}{28}    &  \multicolumn{1}{c|}{46}    &  \multicolumn{1}{c|}{58}    &  \multicolumn{1}{c|}{64}  \\ 
		\hline
		\hline
		\multirow{4}{*}{\begin{tabular}[l]{@{}l@{}} 1 \end{tabular}} & \multirow{2}{*}{\begin{tabular}[l]{@{}l@{}}\textbf{Newton} \\  \textbf{iterations}\end{tabular}}   & \textbf{mean} & - &  5  &  5  &  5  &  5  &  5  &  4  &  4  &  4  &  4  &  4\\ \cline{3-14}
		&  & \textbf{max}  & - &  5  &  5  &  5  &  5  &  5  &  4  &  4  &  4  &  4  &  4\\ \cline{2-14}
		&  \multirow{2}{*}{\begin{tabular}[l]{@{}l@{}}\textbf{Linear Solver} \\ \textbf{iterations} \end{tabular}}   & \textbf{mean} &-  &  10  &  11  &  12  &  11  &  11  &  12  &  13  &  17  &  22  &  26\\ \cline{3-14}
		&  & \textbf{max}  & - &  \multicolumn{1}{c|}{19}    &  \multicolumn{1}{c|}{20}    &  \multicolumn{1}{c|}{20}    &  \multicolumn{1}{c|}{21}    &  \multicolumn{1}{c|}{22}    &  \multicolumn{1}{c|}{25}    &  \multicolumn{1}{c|}{29}    &  \multicolumn{1}{c|}{41}    &  \multicolumn{1}{c|}{54}    &  \multicolumn{1}{c|}{65}  \\ 
		\hline
		\hline
		\multirow{4}{*}{\begin{tabular}[l]{@{}l@{}} 10 \end{tabular}} & \multirow{2}{*}{\begin{tabular}[l]{@{}l@{}}\textbf{Newton} \\  \textbf{iterations}\end{tabular}}   & \textbf{mean}   &  5  &  5  &  5  &  5  &  5  &  4  &  4  &  4  &  4  &  4  &  4\\ \cline{3-14}
		&  & \textbf{max}   &  5  &  5  &  5  &  5  &  5  &  5  &  4  &  4  &  4  &  4  &  4\\ \cline{2-14}
		&  \multirow{2}{*}{\begin{tabular}[l]{@{}l@{}}\textbf{Linear Solver} \\ \textbf{iterations} \end{tabular}}   & \textbf{mean}  &  16  &  16  &  16  &  15  &  16  &  16  &  16  &  17  &  21  &  26  &  30\\ \cline{3-14}
		&  & \textbf{max}    &  \multicolumn{1}{c|}{26}    &  \multicolumn{1}{c|}{26}    &  \multicolumn{1}{c|}{26}    &  \multicolumn{1}{c|}{26}    &  \multicolumn{1}{c|}{27}    &  \multicolumn{1}{c|}{28}    &  \multicolumn{1}{c|}{30}    &  \multicolumn{1}{c|}{35}    &  \multicolumn{1}{c|}{49}    &  \multicolumn{1}{c|}{60}    &  \multicolumn{1}{c|}{72}  \\ 
		\hline
		\hline
		\multirow{4}{*}{\begin{tabular}[l]{@{}l@{}} 20 \end{tabular}} & \multirow{2}{*}{\begin{tabular}[l]{@{}l@{}}\textbf{Newton} \\  \textbf{iterations}\end{tabular}}   & \textbf{mean}   &  5  &  5  &  5  &  5  &  4  &  4  &  4  &  4  &  4  &  4  &  4\\ \cline{3-14}
		&  & \textbf{max}   &  5  &  5  &  5  &  5  &  5  &  4  &  4  &  4  &  4  &  4  &  4\\ \cline{2-14}
		&  \multirow{2}{*}{\begin{tabular}[l]{@{}l@{}}\textbf{Linear Solver} \\ \textbf{iterations} \end{tabular}}   & \textbf{mean}  &  19  &  19  &  19  &  19  &  19  &  19  &  20  &  21  &  26  &  30  &  33\\ \cline{3-14}
		&  & \textbf{max}    &  \multicolumn{1}{c|}{33}    &  \multicolumn{1}{c|}{33}    &  \multicolumn{1}{c|}{33}    &  \multicolumn{1}{c|}{33}    &  \multicolumn{1}{c|}{34}    &  \multicolumn{1}{c|}{36}    &  \multicolumn{1}{c|}{40}    &  \multicolumn{1}{c|}{47}    &  \multicolumn{1}{c|}{58}    &  \multicolumn{1}{c|}{68}    &  \multicolumn{1}{c|}{73}  \\ 
		\hline
	\end{tabular}
	\caption{NHK-C cylinder: average and maximum number of Newton and linear solver iterations recorded along the loading path. 
             Results are obtained varying the stabilization parameters $\beta$ and $\epsilon$ in order to show their influence on the performance of the solution strategy.}
	\label{tab:cylinder_rot_beta_epsilon_tests}
\end{table}

\section{Conclusions}
We developed and numerically validated a framework for the simulation of finite deformations
based on compressible and incompressible hyperelastic material models.
The framework relies on BR2 dG discretizations and allows to impose 
Dirichlet boundary conditions by means of Nitsche method and Lagrange multipliers. 
State of the art agglomeration based $h$-multigrid solution strategies have been successfully employed
to improve efficiency of the solution strategy. 

The proposed BR2 formulation provides the same attractive features of BR1 dG dicretizations 
for a reduced computational cost thanks to a more compact stencil: each cell is coupled solely 
with its neighbouring elements instead of neighbours plus neighbours of neighbours.
In order to better control the amount of stabilization the BR2 stabilization term 
relies on lifting operators defined in a polynomial space that is one degree higher than the 
polynomial space employed for test and trial functions. The approach has demonstrated effective over 
computational meshes composed of elements of standardized shape (triangles and quadrilaterals in 2D, tetrahedrals and hexahedrals in 3D)
and allows to get rid of non-local stabilization parameters based on the number of faces.

We demonstrated that the Lagrange multiplier method for imposing Dirichlet boundary conditions is 
more effective than Nitsche method in the sense that the number of incremental step can be reduced by orders of magnitude. 
Moreover the number of increments is insensitive to mesh density and polynomial degree. 

In order to achieve stability in a broader range of test case configurations, in particular is case of compression solicitations,
the proposed BR2 implementation requires an adaptive stabilization strategy featuring user dependent stabilisation parameters. 
Nevertheless, since the performance of the multigrid solution strategy is pretty insensitive to those stabilization parameters,
the computational expense is not affected by tuning the stabilization, as might be required when dealing with challenging applications.

Future research efforts will consider the possibility to utilize the proposed implementation within an unified high-order accurate framework for fluid-structure interaction where dG methods are employed both for computational fluid-dynamics and computational solid mechanics. 

\section*{Acknowledgements}
The authors acknowledge the support of \textit{Serioplast} (\url{www.serioplast.com}). This research is in partnership with the International Research Training Group (IRTG): DROPIT (Droplet Interaction Technologies-GRK-2160).


\end{document}